\documentclass{cpamart1}
\received{Month 200X}
\volume{000}
\startingpage{1}

\authorheadline{Chen G.-Q. et al}
\titleheadline{Global Solutions of the Compressible Euler-Poisson Equations}

\usepackage{amsfonts}
\usepackage{graphicx}
\usepackage{amssymb,mathrsfs}
\usepackage{marginnote}
\usepackage{color}
\usepackage[leqno]{amsmath}
\usepackage{mathrsfs}
\usepackage{mathptmx}

\numberwithin{equation}{section}
\newcommand{\dis}{\displaystyle}

\newcommand{\R}{\mathbb{R}}

\newcommand{\fb}{\mathfrak{b}}
\newcommand{\longrightharpoonup}{\rm - \!\!\! \rightharpoonup}

\def\v{\varepsilon}

\def\t{\theta}

\def\g{\gamma}
\def\d{\delta}

\def\r{\rho}
\def\s{\sigma}

\def\f{\frac}
\def\dd{{\rm d}}
\def\M{{\mathcal{M}}}
\def\bom{\boldsymbol{\omega}}
\def\bp{\boldsymbol{\psi}}

\newtheorem{theorem}{Theorem}[section]
\newtheorem{lemma}[theorem]{Lemma}
\newtheorem{corollary}[theorem]{Corollary}
\newtheorem{proposition}[theorem]{Proposition}

\theoremstyle{definition}
\newtheorem{definition}[theorem]{Definition}

\theoremstyle{remark}
\newtheorem{remark}[theorem]{Remark}

\begin{document}
\title{Global Solutions of \\ the Compressible Euler-Poisson Equations \\ with Large Initial Data of Spherical Symmetry}

\author{Gui-Qiang G. Chen}{
University of Oxford, Oxford UK}

\author{Lin He}{School of Mathematics, Sichuan University, Chengdu, China}

\author{Yong Wang}{
Academy of Mathematics and Systems Science
   Chinese Academy of Sciences, Beijing, China}

\author{Difan Yuan}{School of Mathematical Sciences, Beijing Normal University, Beijing, China}

\begin{abstract}
We are concerned with a global existence theory for finite-energy solutions of the multidimensional Euler-Poisson equations
for both compressible gaseous stars and plasmas with large initial data of spherical symmetry.
One of the main challenges is the strengthening of waves as they move radially inward towards the origin,
especially under the self-consistent gravitational field for gaseous stars.
A fundamental unsolved problem is whether the density of the global solution
forms a delta measure ({\it i.e.}, concentration) at the origin.
To solve this problem, we develop a new approach for the construction of approximate solutions as the solutions of an appropriately
formulated free boundary problem for the compressible Navier-Stokes-Poisson equations
with a carefully adapted class of degenerate density-dependent viscosity
terms, so that a rigorous convergence proof of the approximate solutions
to the corresponding global solution of the compressible Euler-Poisson equations with large initial data of spherical symmetry can be obtained.
Even though the density may blow up near the origin at a certain time,
it is proved that no delta measure ({\it i.e.}, concentration) in space-time is formed in the vanishing viscosity limit
for the finite-energy solutions of the compressible Euler-Poisson
equations for both gaseous stars and plasmas in the physical regimes under consideration.
\end{abstract}

\maketitle
\tableofcontents

\section{Introduction}
We are concerned with the global existence theory for spherically symmetric solutions of the
multidimensional (M-D) compressible Euler-Poisson equations (CEPEs) with large initial data.
CEPEs govern the motion of compressible gaseous stars or plasmas under a self-consistent
gravitational field
or an electric field,
which take the form:
\begin{align}\label{1.1-1}
	\begin{cases}
		\partial_t \rho+\mbox{div}\M=0,\\[1mm]
		\partial_t \M+\mbox{div}\big(\frac{\M\otimes\M}{\rho}\big)+\nabla p +\rho \nabla\Phi=\mathbf{0},\\
		\Delta \Phi = \kappa \rho,
	\end{cases}
\end{align}
for $t>0$, $\textbf{x}\in\mathbb{R}^n$, and $n\geq3$,
where $\rho$ is the density, $p$ is the pressure, $\M \in\mathbb{R}^n$ represents the momentum,
and $\Phi$ represents the gravitational potential of gaseous stars if $\kappa>0$
and the plasma electric field potential if  $\kappa<0$.
When $\rho>0$, $U=\frac{\M}{\rho}\in\mathbb{R}^n$ is the velocity.
By scaling, we always fix $\kappa=\pm1$ throughout this paper; that is,
$\kappa=1$ for the gaseous star and $\kappa=-1$ for the plasma.
The pressure-density relation is
\begin{equation}\nonumber
	p=p(\r)=a_0\r^{\gamma},
\end{equation}
where $\gamma>1$ is the adiabatic exponent.
Again, by scaling, constant $a_0>0$ may be chosen to be $a_0=\frac{(\gamma-1)^2}{4\gamma}$.

We consider the Cauchy problem for \eqref{1.1-1} with the Cauchy data:
\begin{equation}\label{1.1-2}
	(\rho, \M)|_{t=0}=(\rho_0, \M_0)(\mathbf{x})\longrightarrow (0, \mathbf{0})
	\qquad\,\,\,\,\mbox{as $|\mathbf{x}|\to \infty$},
\end{equation}
subject to the asymptotic condition:
\begin{equation}\label{1.1-3}
	\Phi(t, \mathbf{x})\longrightarrow 0 \qquad \,\,\,\, \mbox{as $|\mathbf{x}|\to \infty$}.
\end{equation}
In \eqref{1.1-2}, the initial far-field velocity has been assumed to be zero in \eqref{1.1-2}
without loss of generality, owing to the Galilean invariance of system \eqref{1.1-1}.
Since a global solution of CEPEs \eqref{1.1-1} normally contains the vacuum states
$\{(t,\mathbf{x})\,:\, \rho(t,\mathbf{x})=0\}$ where
the fluid velocity $U(t,\mathbf{x})$ is not well-defined,
we use the physical variables such as the momentum
$\M(t,\mathbf{x})$, or $\frac{\M(t,\mathbf{x})}{\sqrt{\rho(t,\mathbf{x})}}$
(which will be shown to be always well-defined globally), instead of $U(t,\mathbf{x})$,
when the vacuum states are involved.

The global existence for problem \eqref{1.1-1}--\eqref{1.1-3} is
challenging, mainly owing to the possible appearance of cavitation and concentration,
besides the formation of shock waves, in the solutions,
which leads to the lack of higher order regularity of the solutions, so
that our main focus has to be finite-energy solutions for CEPEs \eqref{1.1-1}.
To solve this existence problem,
we consider the vanishing viscosity limit of the solutions of
the compressible Navier-Stokes-Poisson equations (CNSPEs)
with carefully adapted density-dependent viscosity terms in $\R^n$:
\begin{align}\label{1.1}
	\begin{cases}
		\partial_t \rho+\mbox{div} \M=0,\\[1mm]
		\partial_t \M+\mbox{div}\big(\frac{\M\otimes\M}{\rho}\big)+\nabla p +\rho \nabla \Phi
		=\v \mbox{div}\big(\mu(\r)D(\frac{\M}{\rho})\big)+\v\nabla\big(\lambda(\r)\,\mbox{div}(\frac{\M}{\rho})\big),\\
		\Delta \Phi = \kappa \rho,
	\end{cases}
\end{align}
where $D(\frac{\M}{\rho})=\frac{1}{2}\big(\nabla (\frac{\M}{\rho})+(\nabla (\frac{\M}{\rho}))^\top\big)$ is the stress tensor,
the Lam\'{e} (shear and bulk) viscosity coefficients $\mu(\rho)$ and $\lambda(\rho)$ depend on the density
(that may vanish on the vacuum) and satisfy
\begin{equation*}
	\mu(\rho)\geq0,\quad \mu(\rho)+ n \lambda(\rho)\geq0 \qquad\,\,\, \mbox{for $\rho\geq0$},
\end{equation*}
and parameter $\v>0$ is the inverse of the Reynolds number.
Formally, as $\v\rightarrow 0+$, CNSPEs \eqref{1.1} converge to
CEPEs \eqref{1.1-1}. However, its rigorous mathematical proof has been one of the most challenging open
problems in mathematical fluid dynamics;
see Chen-Feldman \cite{Chen-Feldman2018}, Dafermos \cite{Dafermos}, and the references cited therein.

Many efforts have been made in the analysis of CEPEs \eqref{1.1-1}. We focus mainly on
some relevant time-dependent problems.
Some important progress has been made on the M-D CEPEs with $\kappa=-1$ (plasmas)
in Guo \cite{Guo-1998}, Guo-Ionescu-Pausader \cite{Guo-Ioescu-Pausader-2016}, Guo-Pausader \cite{Guo-Pausader-2011},
and Ionescu-Pausader \cite{Ionescu-Pausader-2013},
in which they proved the global existence of smooth solutions around a constant neutral background
under irrotational, smooth, and localized perturbations of the background with small amplitude.
For the $3$-D gaseous stars problem $(\kappa=1)$,
a compactly supported expanding
classical solution was discovered by Goldreich-Weber \cite{Goldreich-Weber-1980} in 1980;
see also \cite{Fu-Lin-1998,Makino-1992}.
Hadzic-Jang \cite{Hadzic-Jang-2018} proved the nonlinear stability of the Goldreich-Weber
solution under small spherically symmetric perturbations for the adiabatic exponent $\gamma=\frac43$,
while the problem for $\gamma\neq\frac43$ is still widely open.
When the initial density is small and has compact support,
Hadzic-Jang \cite{Hadzic-Jang-2017} constructed a class of global-in-time solutions
of the 3-D CEPEs in the Lagrangian coordinates for $\gamma=1+\frac1k$ with $k\in \mathbb{N}\backslash \{1\}$,
or  $\gamma\in(1,\frac{14}{13})$.
More recently, Guo-Hadzic-Jang \cite{Guo-H-J} constructed an infinite-D
	family of spherically symmetric collapsing solutions of the 3-D CEPEs \eqref{1.1-1} for $\gamma\in(1,\frac{4}{3})$,
	that is, the gas star continuously shrinks to be one point ({\it i.e.}, delta measure);
	see \cite{Guo-H-J-1} for the case that $\gamma=1$.
We also refer \cite{Lei-Gu,Luo-Xin-Zeng-2014,Makino-1986} for the local
well-posedness of smooth solutions.

On the other hand, owing to the strong hyperbolicity and nonlinearity,
the smooth solutions of \eqref{1.1-1} may break down in finite time, especially
when the initial data are large  ({\it cf.} \cite{Chen-Wang,Makino-1992}).
Therefore, the weak solutions have to be considered for the Cauchy problem with large initial data.
For gaseous stars ({\it i.e.}, $\kappa=1$) surrounding a solid ball,
Makino \cite{Makino-1997} proved the local existence of
weak solutions for $\gamma\in(1,\frac53]$  with spherical symmetry;
also see Xiao \cite{Xiao-2016} for global weak solutions
for a class of initial data.
For the compressible Euler equations,
we refer to
\cite{Chen2,Perepelitsa,Chen-Schrecker-2018,Chen-Wang-2018,CS-2012,Chen5,Chen5b,R.J.DiPerna2,Huang-Wang,JM-2015,Ph. LeFloch,Lions P.-L.1,Lions P.-L.2}
and the references therein.

For CNSPEs \eqref{1.1}, many efforts have also been made regarding the global existence of solutions.
For CNSPEs \eqref{1.1} with constant viscosity,
some global existence results for weak solutions for viscous gaseous stars ({\it i.e.}, $\kappa=1$)
have been obtained;
see also \cite{Ducomet-Feireisl-P-2004,Jang-2010,Kobayashi-Suzuki-2008,Kong-Li-2018,Okada-Makino-1993}
and the references cited therein.
For CNSPEs with density-dependent viscosity terms,
Zhang-Fang \cite{Zhang-Fang-2009} obtained a unique global weak solution for
a spherically symmetric vacuum free boundary problem with $\gamma>1$
for a small perturbation around some steady solution;
the global existence of spherically symmetric weak solutions was proved by Duan-Li \cite{Duan-Li}
for the $3$-D problem for $\kappa=1$ and $\gamma\in(\frac{6}{5},\frac{4}{3}]$
with stress free boundary condition and nonzero initial density
for arbitrarily large initial data.
Recently, Luo-Xin-Zeng \cite{Luo-Xin-Zeng-2016-Adv,Luo-Xin-Zeng-2016-CMP} proved
the existence and large-time stability of spherically symmetric smooth solutions of
the $3$-D viscous  problem (with $\kappa=1$) for a small perturbation around the
Lane-Emden solution for $\gamma\in(\frac43,2)$.
For the global existence of solutions of the compressible Navier-Stokes equations,
we refer to \cite{Feireisl-2001,Hoff-1995,Jiang-Zhang-2001,Lions-CNS-1998,Matsumura-1980}
for the case with constant viscosity,
\cite{Li-Xin-2015,VY} for the case with density-dependent viscosity,
and the references cited therein.
In particular, we remark that
the BD entropy estimate
developed in \cite{BD-2003-CPDE} for the derivative estimate of the density
plays a key role in \cite{Li-Xin-2015,VY}.
Such an estimate is based on the new mathematical entropy -- the BD entropy,
first discovered by Bresch-Desjardins \cite{BD-2003-CRMASP} for the particular case $(\mu,\lambda)=(\rho,0)$,
and later
generalized by Bresch-Desjardins  \cite{BD-2004-CRMASP} to include any viscosity
coefficients $(\mu, \lambda)$ satisfying the BD relation:
$\lambda(\rho)=\rho\mu'(\rho)-\mu(\rho)$;
also see \cite{BD-2007-JMPA}.
The BD-type entropy will also be used in this paper.

The idea of regarding inviscid gases as viscous gases with vanishing physical
viscosity can date back to the seminal paper by Stokes \cite{Stokes};
see also the important contributions in \cite{Dafermos,Hugoniot,Rankine,Lord Rayleigh}.
Most of the known results are for the vanishing viscosity limit from the compressible Navier-Stokes
to Euler equations.
The first rigorous convergence analysis of the vanishing physical viscosity limit
from the barotropic Navier-Stokes to Euler equations
was made by Gilbarg \cite{D. Gilbarg}, in which he established the mathematical existence
and vanishing viscosity limit of the Navier-Stokes shock layers.
For the convergence analysis confined in the framework of piecewise smooth solutions;
see \cite{GMWZ,D. Hoff-Liu,Xin-1993} and the references cited therein.
For general data, due to the lack of $L^\infty$ uniform estimate,
the $L^\infty$ compensated compactness framework \cite{Chen5,Chen5b,R.J.DiPerna2,Lions P.-L.1,Lions P.-L.2}
does not apply directly for the vanishing viscosity limit of the compressible Navier-Stokes equations.
LeFloch-Westdickenberg \cite{Ph. LeFloch} first developed an $L^p$ compensated compactness framework
for approximate solutions of the isentropic Euler equations for the adiabatic exponent $\gamma\in (1, \frac{5}{3})$.
In order to establish the vanishing viscosity limit as discussed above,
Chen-Perepelitsa \cite{Perepelitsa} generalized the $L^p$ compensated compactness framework, especially including the whole physical range of adiabatic exponent $\gamma>1$, by further developing/simplifying the proof arguments
and then applied it to establish rigorously the vanishing viscosity limit
of the solutions of the $1$-D compressible Navier-Stokes equations to
the corresponding relative finite-energy solutions of the Euler equations
for large initial data.
Most recently, Chen-Wang \cite{Chen-Wang-2018} established the vanishing viscosity limit of the compressible Navier-Stokes
equations with general data of spherical symmetry and obtained the global existence
of spherically symmetric solutions of the compressible Euler equations with large initial data,
in which it was proved that no delta measure is formed for the density function at the origin.

For problem \eqref{1.1-1}--\eqref{1.1-3}, owing to
the additional difficulties arisen from the possible appearance of concentration and cavitation,
besides the involvement of shock waves,
it has been a longstanding open problem to construct global finite-energy solutions
with large initial data of spherical symmetry.
The key objective of this paper is to solve this problem and establish the global existence of
spherically symmetric finite-energy solutions of \eqref{1.1-1}:
\begin{align}\label{1.2}
	\r(t,\mathbf{x})=\r(t,r),\,\,\,\, \M(t,\mathbf{x})=m(t,r)\f{\mathbf{x}}{r}, \,\,\,\, \Phi(t,\mathbf{x})=\Phi(t,r)
 \qquad \mbox{for $r=|\mathbf{x}|$},
\end{align}
subject to the initial condition:
\begin{align}\label{initial-1}
	(\rho, \M)(0,\mathbf{x})=
	(\rho_0,\M_0)(\mathbf{x})
	=(\rho_0(r),m_0(r)\frac{\mathbf{x}}{r}) \longrightarrow (0,\mathbf{0})
	\qquad \mbox{as $r\to \infty$},
\end{align}
and the asymptotic boundary condition:
\begin{align}\label{initial-2}
	\Phi(t,\mathbf{x})=\Phi(t,r)\longrightarrow 0
	\qquad\mbox{as $r\to \infty$}.
\end{align}
Since $\Phi(0,\mathbf{x})$ can be determined by the initial density in \eqref{initial-1}
and the boundary condition in \eqref{initial-2},
there is no need to impose initial data for $\Phi$.

To achieve this, we establish the vanishing viscosity limit of the corresponding spherically symmetric
solutions of CNSPEs \eqref{1.1} with the
adapted class of degenerate density-dependent viscosity terms and
approximate initial data of similar form to
\eqref{initial-1}.
For spherically symmetric solutions of form \eqref{1.2},
systems \eqref{1.1-1} and \eqref{1.1} become
\begin{align}\label{1.4}
	\begin{cases}
		\dis \r_t+ m_r+\f{n-1}r m=0,\\[1mm]
		\dis m_t+\big(\frac{m^2}{\rho}+p\big)_r+\f{n-1}r \frac{m^2}{\rho}+\rho \Phi_r=0,\\[1mm]
		\dis \Phi_{rr}+\frac{n-1}{r} \Phi_r=\kappa \rho,
	\end{cases}
\end{align}
and
\begin{align}\label{1.3}
	\begin{cases}
		\dis \r_t+ m_r+\f{n-1}r m=0,\\[1mm]
		\dis m_t+\big(\frac{m^2}{\rho}+p\big)_r+\f{n-1}r \frac{m^2}{\rho}
		+\rho \Phi_r = \v\Big((\mu+\lambda)\big((\frac{m}{\rho})_r+\frac{n-1}{r}\frac{m}{\rho}\big)\Big)_r-\v\frac{n-1}{r}\frac{m}{\rho}\mu_r,\\[1mm]
		\dis \Phi_{rr}+\frac{n-1}{r} \Phi_r=\kappa \rho,
	\end{cases}
\end{align}
respectively.

The study of spherically symmetric solutions can date back to the 1950s
and has been motivated by many important physical problems such as stellar dynamics including gaseous stars
and supernova formation \cite{Chandrasekhar,Rosseland,Whitham-1974}.
In fact, the most famous solutions of CEPEs \eqref{1.1-1} are the Lane-Emden steady solutions \cite{Chandrasekhar,Lin},
which describe spherically symmetric gaseous stars in equilibrium and minimize the energy among all possible configurations.
More precisely, for the $3$-D case, there exists a compactly supported and spherically symmetric steady solution
with finite mass for $\gamma\in(\frac65,2)$.
For the time-dependent system, the central feature is the strengthening of waves as they move radially inward near the origin,
especially under the self-gravitational force for gaseous stars.
The spherically symmetric solutions of the compressible Euler equations may blow up near
the origin \cite{Courant-Friedrichs-1962,Li-Wang-2006,MRRS-2020,Whitham-1974} at a certain time in some situations.
However, it has not been well understood how the spherically symmetric solutions of CEPEs \eqref{1.1-1} with
self-gravitational force (which drags the gas particles to the origin) blow up when the initial total-energy is finite.
A fundamental unsolved problem is whether a concentration is formed at the origin;
that is, the density becomes a delta measure at the origin, especially when a focusing spherical shock
is moving inward towards the origin under self-consistent gravitational field.

In this paper, we establish the global existence of finite-energy solutions of problem \eqref{1.1-1}--\eqref{1.1-3} for CEPEs
with spherical symmetry as the vanishing viscosity limits of global weak solutions of
CNSPEs \eqref{1.1} with corresponding initial and asymptotic conditions,
which indicates that no delta measure is formed for the density of the solution of problem \eqref{1.1-1}--\eqref{1.1-3}
in the limit indeed.
To achieve these,
the main point is to establish appropriate uniform estimates in $L^p$ and the $H^{-1}_{\rm loc}$--compactness of the entropy dissipation measures
for the solutions of CNSPEs \eqref{1.3} subject to the corresponding initial and asymptotic conditions.
Owing to the possible appearance of cavitation, the singularity of geometric source terms at the origin, as well as the gravitational force for the gaseous star case,
the global solutions of CNSPEs \eqref{1.3} with large initial data are not smooth in general.
Thus, we start with
the construction of approximate smooth solutions of
the truncated approximate problem \eqref{2.1}--\eqref{initial} for CNSPEs \eqref{1.1},
where the origin and the far-field are cut off, and a stress-free boundary condition is imposed.

In general, we have two basic candidates for the boundary conditions of the approximate problem \eqref{2.1}--\eqref{initial}:
One is to use the Dirichlet boundary conditions: $u(t,a)=u(t,b)=0$ as in \cite{Chen-Wang-2018}, in which case it is difficult
to obtain the higher integrability on the velocity (see Lemma \ref{lem2.6}) due to the far-field vacuum (since the total mass is finite).
Another choice is to use the vacuum free boundary condition; however the global existence of smooth solutions
with the vacuum free boundary condition and large initial data is still an important open problem.
One of our main observations is that the stress-free boundary conditions \eqref{2.6}--\eqref{2.5} we have adapted in \S3
serve our purpose to avoid the difficulties mentioned above.
Even though, we still have to overcome the following additional difficulties:

\begin{enumerate}
	\item[\rm (i)]  Owing to the effect of self-gravitational force for $\kappa=1$,
	we need condition \eqref{1.16} in \S2 to close the basic energy estimate
	for $\gamma\in(\frac{2n}{n+2},\frac{2(n-1)}{n}]$, which implies that
	the initial total mass
	can not be too large when the total initial-energy is fixed.
	The lower bound condition $\gamma>\frac{2n}{2+n}$ is essentially used when we deal with the gravitational potential.
	
	\item[\rm (ii)] To obtain the derivative estimate of the density, we use the BD entropy introduced
	in \cite{BD-2003-CRMASP,BD-2003-CPDE}; also see \cite{BD-2004-CRMASP,BD-2007-JMPA}.
	To close the bound, we need to control the boundary
	term involving $p(\rho_0^{\v,b}(b))b^n$ for the approximate initial data;
	see  \eqref{2.54}.
	To solve this problem, we construct the approximate initial data $(\rho_0^{\v,b},\rho_0^{\v,b}u_0^{\v,b})$
	so that $p(\rho_0^{\v,b}(b)) b^n$
	are uniformly bounded; see \eqref{2.60-1}.
	The details of constructing the approximate initial data are given in Appendix A; see Lemmas \ref{lemA.2}--\ref{propA.1}.

	\item[\rm (iii)] For the
	free boundary problem \eqref{2.1}--\eqref{initial} below,
	a follow-up point is whether the free boundary domain $\Omega_T$ (see \eqref{omega})
	will expand to the whole space as $b\rightarrow\infty$;
	otherwise, it would not be a good approximation to the original Cauchy problem.
	We deal with this difficulty by proving that
	\begin{equation}\label{0.18}
		b(t)\geq \frac{b}{2}\qquad \mbox{for $\, t\in[0,T]$},
	\end{equation}
	provided $b\gg 1$ for any given $T$.
	Property \eqref{2.60-1} of the constructed approximate initial density
	is crucial to prove \eqref{0.18}; see Lemma \ref{lem2.3} for details.

	\item[\rm (iv)] To utilize the $L^p$ compensated compactness framework \cite{Perepelitsa},
	we still need to have the higher velocity integrability.
	We use the entropy pairs $(\eta^{\#}, q^{\#})$  generated by $\psi(s)=\frac12s|s|$.
	Then we have to deal with the boundary term $(q^{\#}-u \eta^{\#})(t,b(t))b(t)^{n-1}$.
	In general, it is impossible to have a uniform bound for both $q^{\#}(t,b(t))b(t)^{n-1}$ and $(u \eta^{\#})(t,b(t))b(t)^{n-1}$.
	One of our key observations is the cancellation between $q^{\#}(t,b(t))$ and $(u \eta^{\#})(t,b(t))$ via observing that
	\begin{align}\label{0.19}
		|q^{\#}-u\eta^{\#}|\leq C\big(\rho^{\gamma} |u|+ \rho^{\gamma+\theta} \big).
	\end{align}
	With the help of the trace estimates in the basic estimates and the BD entropy estimate,
	it serves perfectly to obtain the uniform trace estimate for the terms on the right-hand side of \eqref{0.19}.
	On the other hand, the trace of $u_r$ can be handled by using \eqref{2.5}; see \eqref{2.101} for details.
\end{enumerate}

This paper is organized as follows:
In \S 2, we first introduce the notion of finite-energy solutions of
problem \eqref{1.1-1}--\eqref{1.1-3} for CEPEs
and then state the main theorems of this paper and several remarks.
In \S 3, we first derive some uniform estimates of
the solutions of the free boundary problem \eqref{2.1}--\eqref{initial}
for the approximate CNSPEs.
In \S 4, we establish the global existence of weak finite-energy solutions of \eqref{1.1} with large initial data of spherical symmetry and finite-energy.
Moreover, some uniform estimates in $L^p$ and the $H_{\rm loc}^{-1}$--compactness of entropy dissipation measures for the weak solutions
of CNSPEs \eqref{1.3} are also obtained.
In \S 5, the vanishing viscosity limit of weak solutions of CNSPEs \eqref{1.3} is proved by using the compensated compactness
framework \cite{Perepelitsa}, which leads to a global finite-energy solution of CEPEs \eqref{1.1-1}.
In Appendix A, we construct the approximate initial data with desired properties.

Throughout this paper,
we denote $L^p(\Omega), W^{k,p}(\Omega)$, and $H^k(\Omega)$ as the standard Sobolev
space on domain $\Omega$ for $p\in [1,\infty]$.
We also use $L^p(I;\, r^{n-1}\dd r)$ or $L^p([0, T)\times I;\,r^{n-1}\dd r \dd t)$
for open interval $I \subset \R_+$ with measure $r^{n-1}\dd r\,$ or $r^{n-1}\dd r \dd t\,$ correspondingly, and
$L^p_{\rm loc}([0,\infty); r^{n-1}\dd r)$ to represent $L^p([0,R);\, r^{n-1}\dd r)$ for any fixed $R>0$.

\section{Mathematical Problem and Main Theorems}

In this section, we first introduce the notion of  finite-energy solutions
of problem \eqref{1.1-1}--\eqref{1.1-3} for CEPEs
in $\mathbb{R}_+^{n+1}:=\mathbb{R}_+\times \mathbb{R}^n=[0, \infty)\times \mathbb{R}^n$ for $n\ge 3$.

We assume that the initial data $(\rho_0,\M_0)(\mathbf{x})$ and corresponding initial potential function $\Phi_0(\mathbf{x})$
have both finite initial total-energy:
\begin{align}\label{1.20a}
	E_0:=
	\begin{cases}
		\dis
		\int_{\mathbb{R}^n}\Big(\frac{1}{2}\big|\frac{\M_0}{\sqrt{\rho_0}}\big|^2
		+\rho_0 e(\rho_0)+\frac{1}{2}|\nabla_{\mathbf{x}}\Phi_0|^2\Big)(\mathbf{x})\,\dd\mathbf{x}<\infty
		 &\mbox{for $\kappa=-1$ $($plasmas$)$},\\[4mm]
		\dis
		\int_{\mathbb{R}^n}\Big(\frac{1}{2}\big|\frac{\M_0}{\sqrt{\rho_0}}\big|^2
		+\rho_0 e(\rho_0)\Big)(\mathbf{x})\,\dd\mathbf{x}<\infty\,\, &\hspace{-5mm}\mbox{for $\kappa=1$ $($gaseous stars$)$},
	\end{cases}
\end{align}
and finite initial total-mass:
\begin{equation}\label{1.6-1}
	M:=\int_{\mathbb{R}^n} \rho_0(\mathbf{x})\,\dd\mathbf{x}=\omega_{n} \int_0^\infty \rho_0(r)\,r^{n-1}\dd r <\infty,
\end{equation}
where $e(\rho):=\frac{a_0}{\gamma-1}\rho^{\gamma-1}$ represents the internal energy,
and
$\omega_n:=\frac{2\pi^{\frac{n}{2}}}{\Gamma(\frac{n}{2})}$ denotes the surface area
of the unit sphere in $\mathbb{R}^n$.

\begin{definition}\label{definition-Euler}
	A measurable vector function
	$(\rho, \M, \Phi)$ is said to be a finite-energy solution of the Cauchy problem \eqref{1.1-1}--\eqref{1.1-3}
	if the following conditions hold{\rm :}
	
\begin{enumerate}
\item[\rm (i)]	
		$\rho(t,\mathbf{x})\geq0\,\,\,$ {\it a.e.}, \,\, and $\,\,\, (\mathcal{M}, \frac{\mathcal{M}}{\sqrt{\rho}})(t,\mathbf{x})=\mathbf{0}\,\,$ \, {\it a.e.} \, on
		\, the \, vacuum \, states $\, \{(t,\mathbf{x})\,:\, \rho(t,\mathbf{x})=0\}$.
	
\item[\rm (ii)] For a.e. $t>0$, the total energy is finite{\rm :}

\begin{itemize}
\item For $\kappa=-1$ $($plasmas$)$,
			\begin{align} \label{1.20-1}
				\int_{\mathbb{R}^n}\Big(\frac{1}{2}\big|\frac{\M}{\sqrt{\rho}}\big|^2
				+\rho e(\rho)+\frac{1}{2}|\nabla_{\mathbf{x}}\Phi|^2\Big)(t,\mathbf{x})\, \dd \mathbf{x}
				\leq E_0;
			\end{align}
			
\item For $\kappa=1$ $($gaseous stars$)$,
			\begin{align}\label{1.20-2}
				\quad\quad \begin{cases}
					\dis
					\int_{\mathbb{R}^n}\Big(\frac{1}{2}\big|\frac{\M}{\sqrt{\rho}}\big|^2
					+\rho e(\rho)+\frac{1}{2}|\nabla_{\mathbf{x}}\Phi|^2\Big)(t,\mathbf{x})\, \dd \mathbf{x}
					\leq C(E_0,M),\\[3mm]
					\dis\int_{\mathbb{R}^n}\Big(\frac{1}{2}\big|\frac{\M}{\sqrt{\rho}}\big|^2
					+\rho e(\rho)-\frac{1}{2}|\nabla_{\mathbf{x}}\Phi|^2\Big)(t,\mathbf{x})\, \dd \mathbf{x}\\
					\qquad\leq \dis\int_{\mathbb{R}^n}\Big(\frac{1}{2}\big|\frac{\M_0}{\sqrt{\rho_0}}\big|^2
					+\rho_0 e(\rho_0)-\frac{1}{2}|\nabla_{\mathbf{x}}\Phi_0|^2\Big)\, \dd \mathbf{x}.
				\end{cases}
			\end{align}
		\end{itemize}

\item[\rm (iii)] For any $\zeta(t,\mathbf{x})\in C^{1}_0(\mathbb{R}_+^{n+1})$,
		\begin{align}\label{1.18-1}
			\int_{\mathbb{R}^{n+1}_+} \big(\rho \zeta_t + \M\cdot\nabla\zeta\big)\,\dd\mathbf{x}\dd t
			+\int_{\mathbb{R}^n} \rho_0(\mathbf{x}) \zeta(0,\mathbf{x})\,\dd\mathbf{x}=0.
		\end{align}
		
\item[\rm (iv)] For any $\bp(t,\mathbf{x})=(\psi_1,\cdots,\psi_n)(t,\mathbf{x})\in \big(C^1_0(\mathbb{R}_+^{n+1})\big)^n$,
		\begin{align}\label{1.19-1}
		&\int_{\mathbb{R}^{n+1}_+} \Big(\M\cdot\partial_t\bp +\frac{\M}{\sqrt{\rho}}\cdot \big(\frac{\M}{\sqrt{\rho}}\cdot \nabla\big)\bp
		+p(\rho)\,\mbox{\rm div}\,\bp\Big)(t,\mathbf{x})\, \dd\mathbf{x}\dd t\\
		&+\int_{\mathbb{R}^n} \M_0(\mathbf{x})\cdot \bp(0,\mathbf{x})\,\dd\mathbf{x}
		=\int_{\mathbb{R}^{n+1}_+} (\rho\nabla_{\mathbf{x}} \Phi\cdot \bp)(t,\mathbf{x})\,\dd\mathbf{x}\dd t.\nonumber
\end{align}

\item[\rm (v)] For all $\xi(\mathbf{x})\in C_0^1(\mathbb{R}^n)$,
		\begin{align}\label{1.19-2}
			\int_{\R^n} \nabla_{\mathbf{x}} \Phi(t,\mathbf{x})\cdot \nabla_{\mathbf{x}} \xi(\mathbf{x})\,\dd \mathbf{x}
			=-\kappa \int_{\R^n} \rho(t,\mathbf{x}) \, \xi(\mathbf{x})\,\dd\mathbf{x}\qquad \mbox{for {\it a.e.} $t\geq 0$}.
		\end{align}
	\end{enumerate}
\end{definition}

For the case that $\kappa=1$ (gaseous stars), denote $M_{\rm c}(\gamma)$ as the critical mass given by
\begin{align}\label{1.16}
	M_{\rm c}(\gamma):=
	\begin{cases}
		B_{n,\g}^{-\frac{n}{2}}\qquad\qquad\qquad\qquad &\mbox{for $\gamma=\frac{2(n-1)}{n}$},\\[2mm]
		\left(\frac{(n-2)B_{n,\g}}{n(\gamma-1)}\right)^{-\frac{n(\gamma-1)}{(n+2)\gamma-2n}}\,
		\left(\frac{(n-2)E_0}{(2(n-1)-n\gamma)\omega_n}\right)^{-\frac{2(n-1)-n\gamma}{(n+2)\gamma-2n}}\quad &\mbox{for $\gamma\in(\frac{2n}{n+2}, \frac{2(n-1)}{n})$},
	\end{cases}
\end{align}
with $\displaystyle B_{n,\g}:=\frac{2\kappa}{n(n-2)}\Big(\frac{a_0}{\gamma-1}\Big)^{-\frac{n-2}{n(\gamma-1)}}
\,\omega_n^{\frac{2(n-1)-n\gamma}{n(\gamma-1)}}\, \omega_{n+1}^{-\frac2n}$
depending only on $(n, \gamma)$.

We now state the main theorem of this paper.

\begin{theorem}[Main Theorem I:  Existence of spherically symmetric solutions of CEPEs]\label{thm:2.1}
	Consider problem \eqref{1.1-1}--\eqref{1.1-3} for CEPEs with
	large initial data of spherical symmetry of form \eqref{initial-1}--\eqref{initial-2}.
	Let $(\rho_0, \M_0, \Phi_0)(\mathbf{x})$ satisfy \eqref{1.20a}--\eqref{1.6-1}.
	In addition,

	\begin{enumerate}
		\item[\rm (a)] when $\kappa=-1$ $($plasmas$)$, assume that $\gamma>1$, and
          $\rho_0\in L^{\frac{2n}{n+2}}(\R^n)$ when $\gamma\in (1, \frac{2n}{n+2});$

		\item[\rm (b)] when $\kappa=1$ $($gaseous stars$)$,
		assume that $\gamma>\frac{2(n-1)}{n}$, or $M<M_{\rm c}(\gamma)$ when $\gamma\in (\frac{2n}{n+2}, \frac{2(n-1)}{n}]$.
	\end{enumerate}
	Then there exists a global finite-energy solution $(\rho, \M, \Phi)(t,\mathbf{x})$ of problem \eqref{1.1-1}--\eqref{1.1-3}
	and \eqref{initial-1}--\eqref{initial-2} with
	spherical symmetry of form \eqref{1.2} in the sense of  Definition {\rm \ref{definition-Euler}},
	where $(\rho, m, \Phi)(t,r)$ is determined by
	the corresponding
	system \eqref{1.4} with initial data
	$(\rho_0, m_0, \Phi_0)(r)$ given in \eqref{initial-1} subject to the asymptotic condition \eqref{initial-2}.
\end{theorem}

\begin{remark}
	By the Poisson equation, the initial condition on $\nabla_{\mathbf{x}} \Phi_0$ is indeed a condition on the initial density $\rho_0$.
	In fact, to make the Poisson equation solvable and sense,
	we need the additional condition $\rho_0\in (L^{\frac{2n}{n+2}}\cap L^1)(\R^n)$ for case $\kappa=-1$ {\rm (}plasmas{\rm )}
	if $\gamma>1$
    as required in Theorem {\rm \ref{thm:2.1}}.
	However, for case $\kappa=1$ {\rm (}gaseous stars{\rm )},
	such an additional condition is not required for $\gamma>\frac{2n}{n+2}$.
\end{remark}

\begin{remark}
	To the best of our knowledge, Theorem {\rm \ref{thm:2.1}} provides the first global-in-time solution of
	the M-D CEPEs \eqref{1.1-1} with large initial data.
	For $\kappa=1$ $($gaseous stars$)$, condition $\gamma>\frac{2n}{n+2}${\rm (}{\it e.g.} $\gamma>\frac{6}{5}$ for $n=3${\rm )} is necessary
	to ensure the global existence of finite-energy solutions with finite total mass,
	which
	corresponds to the one for the Lane-Emden solutions $(${\it cf}. {\rm  \cite{Chandrasekhar,Lin}}$)$.
	Moreover, it has been shown in {\rm  \cite{Chandrasekhar,Lin}} that there is no spherically symmetric steady solution
	of gas-stars for the $3$-D CEPEs
	when $\gamma\in (1,\frac{6}{5})$ with finite total mass{\rm ;}
	thus it has been conjectured that there is no global-in-time solution even in the weak sense in general.
\end{remark}

\begin{remark}
For the steady gaseous star problem, it is well-known that there exists no steady white dwarf star with total mass
larger than the Chandrasekhar limit $M_{\rm ch}$ when $\gamma\in (\frac{6}{5}, \frac{4}{3}]${\rm ;} see {\rm \cite{Chandrasekhar}}.
In this paper, for the $3$-D time-dependent gaseous star problem with $\gamma \in(\frac{6}{5},\frac{4}{3}]$,
we also need the restriction on the total mass of the gaseous star{\rm :} $M<M_{\rm c}(\gamma)$ with $M_c(\gamma)$ defined in \eqref{1.16},
which is consistent with the phenomenon observed by Chandrasekhar {\rm \cite{Chandrasekhar}}.
It is interesting to clarify whether the delta measure could be formed at some time when $M<M_{\rm c}(\gamma)$ is violated.
Indeed, when $\gamma\in(1,\frac{4}{3})$, Guo-Hadzic-Jang {\rm \cite{Guo-H-J}} recently constructed an infinite-D family of spherically
symmetric collapsing solutions of the $3$-D CEPEs \eqref{1.1-1}{\rm ;}
that is, the gas star continuously shrinks to be one point $(${\it i.e.}, delta measure$)${\rm ;} see {\rm \cite{Guo-H-J-1}}
for the case that $\gamma=1$.
\end{remark}

To establish Theorem \ref{thm:2.1},
we first construct global weak solutions for CNSPEs \eqref{1.1}
with appropriately adapted degenerate density-dependent viscosity terms and approximate initial data:
\begin{equation}\label{initial-data}
	(\rho, \M, \Phi)|_{t=0}=(\rho_0^\v, \M_0^\v, \Phi_0^\v)(\mathbf{x})\longrightarrow (\rho_0, \M_0, \Phi_0)(\mathbf{x}) \qquad\,\,\,\mbox{as $\v\to 0$},
\end{equation}
constructed as in Appendix A satisfying Lemmas \ref{lemA.2}--\ref{propA.1}, subject to the asymptotic boundary condition:
\begin{equation}\label{initial-data-2}
	\Phi^\v (t,\mathbf{x}) \longrightarrow 0 \qquad\,\,\,\mbox{as $|\mathbf{x}|\to \infty$}.
\end{equation}

For clarity, we adapt the viscosity terms with
$(\mu,\lambda)=(\rho,0)$
in \eqref{1.1} and focus on the case when $\v\in (0,1]$
without loss of generality throughout this paper.

\begin{definition}\label{definition-NSP}
A vector function
$(\rho^\v,\M^\v, \Phi^\v)$  is said to be a weak solution of problem \eqref{1.1} and \eqref{initial-data}--\eqref{initial-data-2}
with $(\mu, \lambda)=(\rho, 0)$
if the following conditions hold{\rm :}

\begin{enumerate}
\item[\rm (i)]	
$\rho^\v(t,\mathbf{x})\geq0\,\,$ {\it a.e.}, and $(\mathcal{M}^\v,\frac{\mathcal{M}^\v}{\sqrt{\rho^\v}})(t,\mathbf{x})=\mathbf{0}\,\,$
{\it a.e.} on the vacuum states $\{(t,\mathbf{x})\,:\,\rho^\v(t,\mathbf{x})=0\}$,
\begin{align*}
&\rho^\v\in L^{\infty}\big([0,T]; L^\gamma(\mathbb{R}^n)\big),
 \quad   \nabla\sqrt{\rho^\v}\in  L^{\infty}\big([0,T]; L^2(\mathbb{R}^n)\big),\\
&\frac{\M^\v}{\sqrt{\rho^\v}}\in L^{\infty}\big([0,T]; L^2(\mathbb{R}^n)\big), \quad \Phi^\v\in L^\infty([0,T];L^{\f{2n}{n-2}}(\mathbb{R}^n) ),\\
&\nabla\Phi^\v\in L^\infty([0,T];L^2(\mathbb{R}^n));
\end{align*}
		
\item[\rm (ii)] For any $t_2\geq t_1\geq0$ and any $\zeta(t,\mathbf{x})\in C^1_0(\mathbb{R}_+^{n+1})$,
the mass equation $\eqref{1.1}_1$ holds in the sense{\rm :}
\begin{align*}
&\int_{\mathbb{R}^n}(\rho^\v \zeta )(t_2,\mathbf{x})\, \dd\mathbf{x}
 -\int_{\mathbb{R}^n}(\rho^\v \zeta )(t_1,\mathbf{x})\, \dd\mathbf{x}\\
&=\int^{t_2}_{t_1}\int_{\mathbb{R}^n}\big(\rho^\v\zeta_t+\M^\v\cdot \nabla \zeta\big)(t,\mathbf{x})\, \dd\mathbf{x} \dd t;
\end{align*}
		
\item[\rm (iii)] For any $\bp=(\psi_1,\cdots,\psi_n)\in \big( C^2_0(\mathbb{R}_+^{n+1})\big)^n$,
the momentum equations $\eqref{1.1}_2$ hold in the sense{\rm :}
\begin{align*}
&\int_{\mathbb{R}^{n+1}_+} \Big(\M^\v\cdot \bp_t
	+\frac{\M^\v}{\sqrt{\rho^\v}}\cdot \big(\frac{\M^\v}{\sqrt{\rho^\v}}\cdot \nabla\big)\bp
	+ p(\rho^\v)\, \mbox{\rm div}\,\bp\Big)\,\dd\mathbf{x}\dd t\\
&\quad +\int_{\mathbb{R}^n} \M_0^\v(\mathbf{x})
 \cdot\bp(0,\mathbf{x})\, \dd\mathbf{x}\nonumber\\
&=-\v\int_{\mathbb{R}^{n+1}_+}
 \Big\{\frac{1}{2}\M^\v\cdot \big(\Delta \bp
 +\nabla\mbox{\rm div}\,\bp\big)
 +\frac{\M^\v}{\sqrt{\rho^\v}}\cdot \big(\nabla\sqrt{\rho^\v}\cdot \nabla\big)\bp \\
&\qquad\qquad\quad\,\,\,+\nabla\sqrt{\rho^\v}\cdot \big(\frac{\M^\v}{\sqrt{\rho^\v}}\cdot \nabla\big)\bp\Big\}\,  \dd\mathbf{x}\dd t
  + \int_{\mathbb{R}^{n+1}_+} (\rho^\v\nabla_{\mathbf{x}} \Phi^\v\cdot \bp )(t,\mathbf{x})\,\dd\mathbf{x}\dd t;
		\end{align*}
		
		\item[\rm (iv)] For any $t\geq 0$ and $\phi(\mathbf{x})\in C_0^1(\mathbb{R}^n)$,
		\begin{align}\nonumber
			\int_{\mathbb{R}^n} \nabla \Phi^\v(t,\mathbf{x})\cdot \nabla \phi(\mathbf{x})\,\dd\mathbf{x}
			=-\kappa\int_{\mathbb{R}^n} \rho^\v(t,\mathbf{x})\, \phi(\mathbf{x})\, \dd\mathbf{x}.
		\end{align}
	\end{enumerate}
\end{definition}

Consider spherically symmetric solutions of form \eqref{1.2}.
Then systems \eqref{1.1-1} and \eqref{1.1} for such solutions
become \eqref{1.4} and \eqref{1.3}, respectively.

A pair
$(\eta(\rho,m),q(\rho,m))$ of functions of $(\rho,m)$
is called an entropy pair of the $1$-D Euler system ({\it i.e.}, consisting of the first two equations of system \eqref{1.4} with $n=1$ and $\Phi\equiv0$)
if the pair $(\eta(\rho,m),q(\rho,m))$ satisfies
$$
\partial_t\eta(\rho(t,r), m(t,r))+ \partial_r q(\rho (t,r),m(t,r))=0
$$
for any smooth solution $(\rho, m)(t,r)$ of the 1-D Euler system;
see Lax \cite{P. D. Lax}.
Furthermore, $\eta(\rho,m)$ is called a weak entropy if
\begin{align}\nonumber
	\eta|_{\rho=0}=0\qquad \text{for any fixed $u=\frac{m}{\rho}$}.
\end{align}
From \cite{Lions P.-L.2},
it is known ({\it cf}. \cite{Perepelitsa,Chen6,Lions P.-L.2})
that any weak entropy $(\eta,q)$ can be represented by
\begin{align}\label{weakentropy}
	\begin{cases}
		\displaystyle \eta^\psi(\rho,m)=\eta(\rho,\rho u)=\int_{\R}\chi(\rho;s-u)\psi(s)\,\dd s,\\[3mm]
		\displaystyle q^\psi(\rho,m)=q(\rho,\rho u)=\int_{\R}(\theta s+(1-\theta)u)\chi(\rho;s-u)\psi(s)\, \dd s,
	\end{cases}
\end{align}
where  $\chi(\rho;s-u)=[\rho^{2\theta}-(s-u)^2]_{+}^{\fb}$ is the kernel
with $\fb=\frac{3-\gamma}{2(\gamma-1)}>-\frac{1}{2}$
and $\theta=\frac{\gamma-1}{2}.$
In particular, when $\psi(s)=\frac{1}{2}s^2,$
the entropy pair is the pair of the mechanical energy and the associated energy flux:
\begin{align}\label{m-entropy}
	\eta^{*}(\rho,m)=\frac{m^2}{2\rho}+\rho e(\rho),\quad q^{*}(\rho,m)=\frac{m^3}{2\rho^2}+m (\rho e(\rho))'.
\end{align}

\begin{theorem}[Main Theorem II: Existence and inviscid limit for CNSPEs]\label{thm1.2}
	Consider CNSPEs \eqref{1.1} with $n\ge 3$
	and the spherically symmetric approximate initial data  \eqref{initial-data} satisfying that
	\begin{align}
		&(\rho^\v_0,m_0^\v)(r) \longrightarrow (\rho_0,m_0)(r)
		\quad \mbox{in $L^{q}([0,\infty);r^{n-1}\dd r)\times L^1([0,\infty);r^{n-1}\dd r)$},\label{1.13-2}\\[2mm]
		&(E_0^\v, E_1^\v)\longrightarrow (E_0,0),  \label{1.13}
			\end{align}
as $\v\to 0+$ for $q=\max\{\gamma,\frac{2n}{n+2}\}$,	and
\begin{align}		
\int_0^{\infty}\rho^\v_0(r)\, r^{n-1}\dd r= \frac{M}{\omega_n},\qquad E_1^\v\leq C(1+M)\v, \label{1.13-5}
\end{align}
\begin{align}
E_0^\v:=
\begin{cases}
\dis \omega_n\int_0^\infty\Big(\frac{1}{2}\Big|\frac{m_0^\v}{\sqrt{\rho_0^\v}}\Big|^2
+ \rho_0^\v e(\rho_0^\v) +\frac{1}{2}|\Phi_{0r}^\v|^2\Big)(r)\, r^{n-1}\dd r<\infty \, &\mbox{for $\kappa=-1$},\\[3.5mm]
\dis \omega_n\int_0^\infty\Big(\frac{1}{2}\Big|\frac{m_0^\v}{\sqrt{\rho_0^\v}}\Big|^2+\rho_0^\v e(\rho_0^\v)\Big)(r)\,r^{n-1}\dd r<\infty \,\, &\mbox{for $\kappa=1$},
\end{cases}\label{1.13-3}\\[1mm]
&\hspace{-10cm}E_1^\v:=\omega_n\v^2\int_0^\infty \big|(\sqrt{\rho_0^\v})_r\big|^2\, r^{n-1} \dd r. \label{1.13-1}
\end{align}
	
	In addition,
	\begin{enumerate}
		\item[\rm (a)]
		when $\kappa=-1$ $($plasmas$)$, assume that $\gamma>1${\rm ;}
		
		\item[\rm (b)] when $\kappa=1$ $($gaseous stars$)$, assume that $\gamma>\frac{2(n-1)}{n}$ or that there exists $\v_0\in (0,1]$ such that, for $\v\in(0,\v_0]$,
		\begin{equation}\label{1.13-4}
			M<M_{\rm c}^\v(\gamma)\qquad \mbox{ when $\gamma\in (\frac{2n}{n+2}, \frac{2(n-1)}{n}]$},
		\end{equation}
		where $M_{\rm c}^\v(\gamma)$ is the critical mass defined by replacing $E_0$ in \eqref{1.16} with $E_0^\v$.
	\end{enumerate}

	\noindent Then the following results hold{\rm :}

	\noindent{\bf Part I (Existence of global solutions of CNSPEs).}
	For each fixed  $\v\in(0,\v_0],$ there exists a globally-defined spherically symmetric weak solution{\rm :}
\begin{align*}
(\rho^\v, \M^\v, \Phi^\v)(t,\mathbf{x})
&=(\rho^\v(t,r),m^\v(t,r) \frac{\mathbf{x}}{r}, \Phi^\v(t,r))\\
&=(\rho^\v(t,r),\rho^\v(t,r)u^\v(t,r) \frac{\mathbf{x}}{r}, \Phi^\v(t,r))
\end{align*}
	of problem \eqref{1.1} and \eqref{initial-data}--\eqref{initial-data-2} for CNSPEs
	in the sense of Definition {\rm \ref{definition-NSP}}, where
\begin{align*}
&u^\v(t,r)=\frac{m^\v(t,r)}{\rho^\v(t,r)}\quad \mbox{a.e. on $\{(t,r)\,:\, \rho^\v(t,r)\ne 0\}$}, \\
& u^\v(t,r)=0 \quad \mbox{a.e. on $\{(t,r)\,:\, \rho^\v(t,r)=0\}$}.
\end{align*}
Moreover,  $(\rho^\v,m^\v, \Phi^\v)(t,r)$ satisfies the following{\rm :}	
	
\begin{enumerate}
\item[\rm (i)] For any fixed $T\in (0,\infty)$, the following uniform bounds hold for any $t\in [0,T]${\rm :}
\begin{align}
			&\int_0^\infty \rho^{\v}(t,r)\,r^{n-1}\dd r=\int_0^{\infty}\rho^\v_0(r)\,r^{n-1}\dd r=\frac{M}{\omega_n},\label{1.17}\\
			& \Phi^\v_r(t,r)=\frac{\kappa}{r^{n-1}} \int_0^r \rho^\v(t,z)\, z^{n-1}\dd z,\label{1.17-1}\\
			&\int_0^\infty \eta^{*}(\rho^\v,m^\v)(t,r)\,r^{n-1}\dd r+\v\int_{\R^2_+} (\rho^\v |u^\v|^2)(t,r)\,r^{n-3} \dd r\dd t \label{1.18}\\
			&\quad
			+\int_0^{\infty} \Big(\int_0^r  \rho^\v(t,z)\, z^{n-1}\dd z\Big)\rho^{\v}(t,r)\,r \dd r +\|\Phi^\v(t)\|_{L^{\frac{2n}{n-2}}(\mathbb{R}^n)} \nonumber\\[1.5mm]
			&\quad +\|\nabla \Phi^\v(t)\|_{L^2(\mathbb{R}^n)} \leq C(M,E_0),\nonumber\\[3mm]
	&\v^2\int_0^{\infty}\big|\big(\sqrt{\r^{\v}(t,r)}\big)_r\big|^2 r^{n-1}\dd r+\v\int_0^T\int_{0}^\infty
             \big|\big((\rho^{\v})^{\frac{\gamma}{2}}\big)_r\big|^2r^{n-1}\dd r\dd t\label{1.19}\\[1.5mm]
		&\quad \leq C(M,E_0, T), \nonumber\\[3mm]
			&\int_0^T\int_d^{D} \big(\rho^{\v}(t,r)|u^{\v}(t,r)|^3+(\rho^{\v}(t,r))^{\gamma+\theta}\big)\,\dd r \dd t\leq C(d,D,M,E_0,T),\label{1.20}\\[2mm]
			&\int_0^T\int_d^{D}(\rho^{\v}(t,r))^{\gamma+1}\,\dd r\dd t\leq C(d,D,M,E_0,T),\label{1.21}
\end{align}
for any compact subset $[d,D]\Subset (0,\infty)$,
where and whereafter $C(M,E_0)>0$, $C(M,E_0,T)>0$,
and $C(d, D, M, E_0,T)>0$ are three universal constants independent of $\v$, but may depend on $(\gamma,n)$
and $(d, D, M, E_0,T)$, respectively.

\item[\rm (ii)] The following energy inequality holds for both $\kappa=\pm 1${\rm :}
\begin{align}\label{1.12-1}
&\int_{\mathbb{R}^n}\Big(\frac{1}{2}\Big|\frac{\M^\v}{\sqrt{\rho^\v}}\Big|^2
+\rho^\v e(\rho^\v)-\frac{\kappa}{2}|\nabla_{\mathbf{x}}\Phi^\v|^2\Big)(t,\mathbf{x})\, \dd \mathbf{x}\\
&\leq \int_{\mathbb{R}^n}\Big(\frac{1}{2}\Big|\frac{\M_0^\v}{\sqrt{\rho_0^\v}}\Big|^2
+\rho_0^\v e(\rho_0^\v)-\frac{\kappa}{2}|\nabla_{\mathbf{x}}\Phi_0^\v|^2\Big)\, \dd \mathbf{x}
\qquad\,\, \mbox{for $t\geq 0$}.\nonumber
\end{align}
		
\item[\rm (iii)]
Let $(\eta^\psi, q^\psi)$ be an entropy pair defined in \eqref{weakentropy}
for a smooth  function $\psi(s)$ of compact support on $\mathbb{R}$.
Then, for $\v\in (0,\v_0]$,
\begin{align}\label{1.22}
\partial_t\eta^\psi(\rho^\v,m^\v)+\partial_rq^\psi(\rho^\v,m^\v) \qquad \mbox{is compact in }\, H^{-1}_{\rm loc}(\mathbb{R}^2_+),
\end{align}
where $H^{-1}_{\rm loc}(\mathbb{R}^2_+)$ represents $H^{-1}((0,T]\times I)$ for any $T>0$ and bounded open subset $I\Subset (0,\infty)$,
and $\mathbb{R}^2_+:=\{(t,r)\ :\ t\in (0,\infty),\ r\in(0,\infty) \}$.
\end{enumerate}

\noindent{\bf Part II (Inviscid limit and existence of global solutions of CEPEs).} For the global weak
	solutions $(\rho^{\v},\M^{\v}, \Phi^\v)(t,\mathbf{x})=(\rho^{\v}(t,r),m^{\v}(t,r) \frac{\mathbf{x}}{r}, \Phi^\v(t,r))$
	of problem  \eqref{1.1} and \eqref{initial-data}--\eqref{initial-data-2} for CNSPEs
	established in  Part {\rm I},
	there exist both a subsequence {\rm (}still denoted{\rm )} $(\rho^{\v},m^{\v}, \Phi^\v)(t,r)$ and a vector function $(\rho, m, \Phi)(t,r)$
	such that, as $\v\rightarrow0$,
	\begin{align*}
	(\rho^{\v},m^{\v})(t,r)\longrightarrow (\rho,m)(t,r)\qquad
		\mbox{in $L^{q_1}_{\rm loc}(\R^2_+)\times L^{q_2}_{\rm loc}(\R^2_+)$},
	\end{align*}
 with $q_1\in[1,\gamma+1)$ and $q_2\in [1,\frac{3(\gamma+1)}{\gamma+3})$, and
	\begin{align}\nonumber
		\begin{split}
			&\Phi^{\v} \rightharpoonup \Phi \quad  \mbox{weakly in $L^2(0,T; H^1_{\rm loc}(\mathbb{R}^n))$},\\
			& \Phi_r^{\v}(t,r)r^{n-1}=\kappa\int_0^r  \rho^{\v}(t,z)\, z^{n-1}\dd z\nonumber\\
			&\qquad\qquad\quad\longrightarrow \, \Phi_r(t,r)r^{n-1}
			=\kappa\int_0^r  \rho(t,z)\,z^{n-1}\dd z \quad  \mbox{a.e. $(t,r)\in \R_+^2$},\\
			& \int_0^\infty|\Phi_r^{\v}(t,r)-\Phi_r(t,r)|^2\, r^{n-1} \dd r \to 0 \qquad \mbox{if $\gamma>\frac{2n}{n+2}$},
		\end{split}
	\end{align}
	with $\displaystyle(\rho,\M, \Phi)(t,\mathbf{x}):=(\rho(t,r),\ m(t,r)\frac{\mathbf{x}}{r},\ \Phi(t,r))$
	to be a global spherically symmetric finite-energy
	solution of problem \eqref{1.1-1}--\eqref{1.1-3}
	for CEPEs in the sense of Definition {\rm \ref{definition-Euler}}.
\end{theorem}

\begin{remark}\label{remark:2.5a}
	In Theorem {\rm \ref{thm1.2}}, the approximate initial data functions $(\rho_0^\v, \M_0^\v,\Phi_0^\v)$
	satisfying conditions {\rm \eqref{1.13-2}}--{\rm \eqref{1.13-4}} are constructed
	in Lemmas {\rm \ref{lemA.2}--\ref{propA.1}} in Appendix A.
	The restriction, $\v\in(0,\v_0]$,
	for the gaseous star case $\kappa=1$ with $\gamma\in (\frac{2n}{n+2}, \frac{2(n-1)}{n}]$
	is mainly due to the construction of approximate initial data in Appendix A.
	Then Theorem {\rm \ref{thm:2.1}} is a direct corollary of Theorem {\rm \ref{thm1.2}}.
	From now on, we always assume $\v\in(0,\v_0]$
	without loss of our main objectives for the inviscid limit.
\end{remark}

\section{Construction and Uniform Estimates of Approximate Solutions}\label{section3}

In order to deal with the difficulties for the appearance of cavitation and the singularity at the origin,
besides shock waves,
as well as uniform estimates of approximate solutions,
we construct our approximate solutions via the following approximate free boundary
problem for CNSPEs:
\begin{equation}\label{2.1}
	\begin{cases}
		\displaystyle\rho_t+(\r u)_r+\frac{n-1}{r}\rho u=0,\\[2mm]
		\displaystyle(\rho u)_t+(\rho u^2+p(\rho))_r+\frac{n-1}{r}\rho u^2 +\frac{\kappa\rho}{r^{n-1}}\int_a^r \rho(t,z)z^{n-1} \dd z\\[2mm]
		\qquad\qquad\qquad\qquad\qquad\quad =\v\big(\r(u_r+\frac{n-1}{r}u)\big)_r-\v\frac{n-1}{r}u\rho_r,
	\end{cases}
\end{equation}
for $(t,r)\in \Omega_T$ with moving domain:
\begin{equation}\label{omega}
	\Omega_T=\{(t,r)\ :\,  a\leq r\leq b(t),\,0\leq t\leq T \},
\end{equation}
where $\{r=b(t):0<t\leq T\}$ is a free boundary determined by
\begin{equation}\label{2.6}
	\begin{cases}
		b'(t)=u(t,b(t)) \qquad \mbox{for $t>0$},\\
		b(0)=b,
	\end{cases}
\end{equation}
and $a=b^{-1}$ with $b\gg1$.
On the free boundary $r=b(t)$, the stress-free boundary condition is chosen:
\begin{equation}\label{2.5}
	\big(p(\rho)-\v\rho(u_r+\frac{n-1}{r}u)\big)(t,b(t))=0 \qquad \mbox{for $t>0$}.
\end{equation}
On the fixed boundary $r=a=b^{-1}$, we impose the Dirichlet boundary condition:
\begin{equation}\label{2.2}
	u|_{r=a}=0\qquad\, \mbox{for $t>0$}.
\end{equation}
The initial condition is
\begin{equation}\label{initial}
	(\rho,\rho u)|_{t=0}=(\r_0^{\v,b},\rho_0^{\v,b}u_0^{\v,b})(r)\qquad\,\mbox{for $ r\in [a,b]$}.
\end{equation}
We always assume that the initial data functions $(\rho_0^{\v,b}, u_0^{\v,b})(r)$ are smooth
and compatible with the boundary conditions \eqref{2.5}--\eqref{2.2},
and $0<C_{\v,b}^{-1}\leq \rho_0^{\v,b}(r)\leq C_{\v,b}<\infty$.

For later use, we define
\begin{align}\label{2.7-1}
	E_{0}^{\v,b}:=
	\begin{cases}
		\dis \int_{a}^{b} \Big\{\rho_0^{\v,b}\big(\frac{1}{2}|u_0^{\v,b}|^2+ e(\rho_0^{\v,b})\big)
		+\frac{1}{2r^{2(n-1)}} \Big(\int_a^r \rho_0^{\v,b}(z)\,z^{n-1} \dd z\Big)^2\Big\}\,\omega_n r^{n-1}\dd r\\
		\hspace{9cm} \mbox{for $\kappa=-1$},\\[3mm]
		\dis \int_{a}^{b} \rho_0^{\v,b}\big(\frac{1}{2} |u_0^{\v,b}|^2+ e(\rho_0^{\v,b})\big)\,\omega_n r^{n-1}\dd r \hspace{3.26cm}\mbox{for $\kappa=1$},
	\end{cases}\\[1.5mm]
	 E_1^{\v,b}:= \v^2 \int_a^b \big|\big(\sqrt{\rho_0^{\v,b}}\big)_r\big|^2\,\omega_n r^{n-1} \dd r.\hspace{5cm}&\label{2.7-2}
\end{align}
When $\kappa=1$, for the given total energy $E_0^{\v,b}>0$, similar to \eqref{1.16},  we  define the critical mass:
\begin{align}\label{2.36}
	M_{\rm c}^{\v,b}(\gamma):=
	\begin{cases}
		B_{n,\g}^{-\frac{n}{2}}\qquad\qquad\qquad\quad &\mbox{for $\gamma=\frac{2(n-1)}{n}$},\\[2mm]
		\left(\frac{(n-2)B_{n,\g}}{n(\gamma-1)}\right)^{-\frac{n(\gamma-1)}{(n+2)\gamma-2n}}\,
		\left(\frac{(n-2)E_0^{\v,b}}{(2(n-1)-n\gamma)\omega_n}\right)^{-\frac{2(n-1)-n\gamma}{(n+2)\gamma-2n}}\, &\mbox{for $\gamma\in(\frac{2n}{n+2}, \frac{2(n-1)}{n})$}.
	\end{cases}
\end{align}

For the initial data $(\rho^\v_0,m^\v_0)$ imposed in \eqref{initial-data} satisfying \eqref{1.13-2}--\eqref{1.13-4},
it follows from Lemma \ref{propA.1} in Appendix A that
there exists a sequence of smooth functions $(\rho^{\v,b}_0, u_0^{\v,b})$ defined on $[a,b]$ such that,
as $ b\rightarrow \infty$,
\begin{align}\label{5.11}
\begin{cases}
(\rho_0^{\v,b}, \rho_0^{\v,b}u_0^{\v,b})\longrightarrow (\rho_0^\v, m_0^\v) \quad\mbox{in $L^{\hat{q}}([a,b];r^{n-1}\dd r)\times L^1([a,b];r^{n-1}\dd r)$},\\[1.5mm]
(E_0^{\v,b}, E_1^{\v,b}) \longrightarrow (E_0^\v, E_1^\v),
\end{cases}
\end{align}
with $\hat{q}\in \{1,\gamma\}$ when $\kappa=1$ (gaseous stars) and $\hat{q}\in \{1,\gamma,\frac{2n}{n+2}\}$ when $\kappa=-1$ (plasmas), and
\begin{align}
	&\int_a^b\rho_0^{\v,b}(r)\,r^{n-1}\dd r=\frac{M}{\omega_n}>0,\qquad  E_1^{\v,b}+E_1^\v\leq C(1+M)\v, \label{2.7}\\
	&\,\, \rho_0^{\v,b}(b) \cong b^{-(n-\alpha)} \qquad \mbox{with $\alpha:=\min\{\frac{1}{2},(1-\frac{1}{\gamma})n\}$}.\label{2.60-1}
\end{align}
Moreover, for each fixed $\v\in(0,\v_0]$, there exists a large constant $\mathfrak{B}(\v)>0$ such that, when $\kappa=1$ with $\gamma\in(\frac{2n}{n+2}, \frac{2(n-1)}{n}]$,
\begin{equation}
	M<M_{\rm c}^{\v,b}(\gamma)\qquad\,\, \mbox{for $b\geq \mathfrak{B}(\v)$}. \label{2.60-3}
\end{equation}
Property \eqref{2.60-1} is important for us to close the BD-type entropy estimate in Lemma \ref{lem2.2} below.

Once problem \eqref{2.1}--\eqref{initial} is solved, we define the potential function $\Phi$  to be the solution of the Poisson equation:
\begin{equation}\label{2.24}
	\displaystyle\Delta \Phi=\kappa\rho\,\mathbf{1}_{\Omega_t},\qquad \lim_{|\textbf{x}|\rightarrow\infty} \Phi=0,
\end{equation}
with  $\Omega_t=\{\textbf{x}\in\mathbb{R}^n\,:\,a\leq |\textbf{x}|\leq b(t) \}$,
for which we have extended $\rho$ to be zero outside $\Omega_t$, where $\mathbf{1}_{\Omega_t}$ is the indicator function
of $\Omega_t$ (which is $1$ when $\textbf{x}\in \Omega_t$ and $0$ otherwise).
In fact, we can show that $\Phi(t,\mathbf{x})=\Phi(t,r)$ with
\begin{equation}\label{2.24-1}
	\Phi_r(t,r)=
	\begin{cases}
		\displaystyle 0 \qquad &\mbox{for $0\leq r\leq a$},\\[2mm]
		\displaystyle \frac{\kappa}{r^{n-1}} \int_a^r  \rho(t,z)\,z^{n-1}\dd z\quad &\mbox{for $a\leq r\leq b(t)$},\\[3mm]
		\displaystyle \frac{\kappa M}{\omega_n r^{n-1}} \quad &\mbox{for $r\geq b(t)$}.
	\end{cases}
\end{equation}

In this section, parameters $(\v, b)$ are fixed such that $\v\in(0,\v_0]$ and $b\geq \mathfrak{B}(\v)$.
For $n=3$, the existence of global smooth solutions $(\r^{\v,b}, u^{\v,b})$ of problem \eqref{2.1}--\eqref{initial}
has been proved by Duan-Li \cite{Duan-Li} for $\gamma\in(\frac65,\frac43]$ and $\kappa=1$ with $0<\r^{\v,b}(t,r)<\infty$.
In fact, for $n\geq 3$, the global existence of smooth solutions of our approximate problem \eqref{2.1}--\eqref{initial}
can be obtained by using similar arguments as in \cite[\S 3]{Duan-Li} for $\kappa=-1$ with $\gamma\in(1,\infty)$,
and for $\kappa=1$ with $\gamma\in(\frac{2(n-1)}{n},\infty)$, or with $\gamma\in (\frac{2n}{n+2},\frac{2(n-1)}{n}]$
and $M<M_{\rm c}^{\v,b}(\gamma)$,
so we omit the details here.

Notice that the upper and lower bound of $\r^{\v,b}$ in \cite{Duan-Li} depend on parameters $(\v, b)$.
Therefore, some careful uniform estimates, independent of $b$, are required so that we can take the limit: $b\rightarrow\infty$
to obtain the global weak solutions of  problem \eqref{1.1}
and \eqref{initial-data}--\eqref{initial-data-2}
in \S 4 below
as approximate solutions of problem \eqref{1.1-1}--\eqref{1.1-3}.
Throughout this section, for simplicity, we drop the superscript in both the approximate solutions
$(\rho^{\v, b}, u^{\v,b})(t,r)$ and the approximate initial data $(\rho^{\v,b}_0, u^{\v,b}_0)$
when no confusion arises.

For strong solutions, it is convenient to deal with IBVP \eqref{2.1}--\eqref{initial} in the Lagrangian coordinates.
It follows from \eqref{2.6} that
\begin{equation*}
	\frac{\dd}{\dd t}\int_a^{b(t)}\rho(t,r)\,r^{n-1}\dd r= (\rho u)(t,b(t))b(t)^{n-1}-\int_a^{b(t)}(\rho u r^{n-1})_r(t,r)\,\dd r=0,
\end{equation*}
which yields that
\begin{equation}\label{2.8}
	\int_a^{b(t)}\rho(t,r)\,r^{n-1}\dd r=\int_a^b\rho_0(r)\,r^{n-1}\dd r=\frac{M}{\omega_n}\qquad \mbox{for any $t\geq0$}.
\end{equation}
For $r\in[a,b(t)]$ and $t\in[0,T]$, we define the Lagrangian coordinates $(\tau,x)$ as
\begin{equation*}
	x(t,r)=\int_a^r \rho(t,y)\,y^{n-1}dy,\quad  \tau=t,
\end{equation*}
which translates domain $[0,T]\times[a,b(t)]$ into a fixed domain $[0,T]\times[0,\frac{M}{\omega_n}]$.
A direct calculation shows that
\begin{align*}\label{2.9}
&\nabla_{(t,r)} x =(-\rho ur^{n-1},\rho r^{n-1}),
\,\,\, \nabla_{(t,r)}\tau=(1, 0),\\
&\nabla_{(\tau,x)}r=(u, \rho^{-1}r^{1-n}), \,\,\, \nabla_{(\tau,x)}t=(1,0).
\end{align*}

Applying the Euler-Lagrange transformation, IBVP \eqref{2.1}--\eqref{2.6} becomes
\begin{equation}\label{2.10}
	\begin{cases}
		\displaystyle\rho_\tau+\rho^2(r^{n-1}u)_x=0,\\[2mm]
		\displaystyle u_\tau+r^{n-1} p_x=-\kappa \frac{x}{r^{n-1}}+\v r^{n-1}\big(\rho^2(r^{n-1}u)_x\big)_x-(n-1)\v r^{n-2}u\rho_x,
	\end{cases}
\end{equation}
for $(\tau,x)\in[0,T]\times[0,\frac{M}{\omega_n}]$, and
\begin{equation}\label{2.11}
	u(\tau,0)=0, \quad \big(p-\v \rho^2(r^{n-1}u)_x\big)(\tau,\frac{M}{\omega_n})=0 \qquad\,\mbox{for $\tau\in[0,T]$},
\end{equation}
where $r=r(\tau,x)$ is defined by
\begin{equation}\label{2.12}
	\frac{\dd}{\dd\tau}r(\tau,x)=u(\tau,x) \qquad \mbox{for $x\in[0,\frac{M}{\omega_n}]$ and $\tau\in[0,T]$},
\end{equation}
and the fixed boundary $x=\frac{M}{\omega_n}$ corresponds to the free boundary $b(\tau)=r(\tau,\frac{M}{\omega_n})$
in the Eulerian coordinates.

\begin{lemma}[Basic energy estimate]\label{lem2.1}
	Any smooth solution $(\rho, u)(t,r)$ of  problem \eqref{2.1}--\eqref{initial} satisfies  the following energy identity{\rm :}
	\begin{align}\label{2.13-1}
     &\int_a^{b(t)}\big(\frac{1}{2}\rho u^2 + \rho e(\rho)\big)(t,r)\,r^{n-1} \dd r
     -\frac{\kappa}{2} \int_a^\infty \frac{1}{r^{n-1}} \Big(\int_a^r \rho(t,z)\, z^{n-1} \dd z\Big)^2 \dd r\\
     &\quad+\v\int_0^t\int_a^{b(s)} \big(\rho u_r^2+(n-1)\rho \f{u^2}{r^2}\big)(s,r)\,r^{n-1} \dd r\dd s\nonumber\\
     &\quad +(n-1) \v\int_0^t (\r u^2) (s,b(s))b(s)^{n-2}\,\dd s\nonumber\\
     &=\int_a^{b}\big(\frac{1}{2}\rho_0u_0^2 +\rho_0e(\rho_0)\big)(r)\,r^{n-1} \dd r
     -\frac{\kappa}{2} \int_a^\infty \frac{1}{r^{n-1}} \Big(\int_a^r \rho_0(z)\, z^{n-1} \dd z\Big)^2\,\dd r,\nonumber
	\end{align}
	where $\rho(t,r)$ is understood to be $0$ for $r\in[0,a]\cup (b(t),\infty)$ in the second term of
	the right-hand side $($RHS$)$ and the second term of the left-hand side $($LHS$)$.  In particular, the following
	estimates hold{\rm :}
	
	\smallskip
	\noindent{Case} {\rm 1:} If $\kappa=-1$ $($plasmas$)$ with $\gamma>1$, then
\begin{align}\label{2.13}
	&\int_{a}^{b(t)}\big(\frac{1}{2}\rho u^2+ \rho e(\rho)\big)(t,r) \, r^{n-1} \dd r+\frac{1}{2}\int_a^{\infty} \frac{1}{r^{n-1}}\Big(\int_a^r \rho(t,z)\, z^{n-1} \dd z\Big)^2 \dd r\\
	&\quad +\v\int_0^t\int_a^{b(s)} \big(\rho u_r^2+(n-1)\rho \frac{u^2}{r^2} \big)(s,r)\,r^{n-1}\dd r \dd s\nonumber\\
	&\quad +(n-1)\v\int_{0}^t (\rho u^2)(s,b(s))b(s)^{n-2}\,\dd s\nonumber\\
	&=\int_{a}^{b} \Big(\big(\frac{1}{2}\rho_0 u_0^2+ \rho_0 e(\rho_0)\big)(r)
	+\frac{1}{2 r^{2(n-1)}}\big(\int_a^r \rho_0(z)\,z^{n-1} \dd z\big)^2\Big)\, r^{n-1} \dd r
	:=\frac{E_{0}^{\v,b}}{\omega_n}.\nonumber
\end{align}
	
	\noindent{Case} {\rm 2:} $\,$ If $\kappa=1$ $($gaseous stars$)$ with $\gamma\in(\frac{2n}{n+2}, \frac{2(n-1)}{n}]$ and $M<M_{\rm c}^{\v,b}(\gamma)$,
	then
	\begin{align}\label{2.14}
     &\int_{a}^{b(t)}\big(\frac{1}{2}\rho u^2+C_{\gamma}\rho e(\rho)\big)(t,r)\, r^{n-1} \dd r\\
     &\quad +\v\int_0^t\int_a^{b(s)} \big(\rho u_r^2+(n-1)\rho \frac{u^2}{r^2} \big)(s,r)\, r^{n-1} \dd r \dd s\nonumber\\
     &\quad+(n-1)\v\int_{0}^t (\rho u^2)(s,b(s))b(s)^{n-2}\,\dd s\nonumber\\
     &\leq \int_{a}^{b} \big(\frac{1}{2}\rho_0 u_0^2+\rho_0 e(\rho_0)\big)(r)\,r^{n-1}\dd r=:\frac{E_0^{\v,b}}{\omega_n},\nonumber
	\end{align}
	where the positive constant $C_{\gamma}>0$ is defined as
	\begin{align}\label{2.14-1}
		C_{\gamma}:=
		\begin{cases}
			\displaystyle 1-
			B_{n,\g} M^{\frac{2}{n}}
			\qquad &\mbox{for $\gamma=\frac{2(n-1)}{n}$},\\[3mm]
			\displaystyle\frac{2(n-1)-n\gamma}{n-2}\qquad &\mbox{for $\gamma\in (\frac{2n}{n+2}, \frac{2(n-1)}{n})$}.
		\end{cases}
	\end{align}
	
	\noindent{Case {\rm 3:}} If $\kappa=1$ $($gaseous stars$)$ with $\gamma> \frac{2(n-1)}{n}$,  then
	\begin{align}\label{2.15}
    \int_{a}^{b(t)}&\frac{1}{2}\big(\rho u^2+ \rho e(\rho)\big)(t,r)\,r^{n-1}\dd r\\
    &+\v\int_0^t\int_a^{b(s)} \big(\rho u_r^2+(n-1)\rho \frac{u^2}{r^2} \big)(s,r)\,r^{n-1} \dd r \dd s\nonumber\\
    &+(n-1)\v\int_{0}^t (\rho u^2)(s,b(s))b(s)^{n-2}\,\dd s\nonumber\\
    \leq &\,\, \frac{E_0^{\v,b}}{\omega_n}+C(M),\nonumber
	\end{align}
	where $C(M)>0$ is some positive constant depending only on the total initial-mass $M$.
\end{lemma}

\begin{proof} We divide the proof into four steps.

\smallskip
1. Multiplying $\eqref{2.10}_2$ by $u$ and then integrating the resultant equation over $x\in[0,\frac{M}{\omega_n}]$ yield that
\begin{align}\label{2.16}
&\frac{1}{2}\frac{\dd}{\dd\tau}\int_0^{\frac{M}{\omega_n}}u^2 \dd x+\int_0^{\frac{M}{\omega_n}}\big(p-\v\r^2(r^{n-1}u)_x\big)_xu\,r^{n-1}\dd x\\
&=-\v (n-1)\int_0^{\frac{M}{\omega_n}} u^2 \rho_x r^{n-2}\,\dd x
-\kappa \int_0^{\frac{M}{\omega_n}}\frac{x }{r^{n-1}} u \,\dd x.\nonumber
\end{align}
For the second term of \eqref{2.16}-LHS ({\it i.e.}, the left-hand side of \eqref{2.16}),
it follows from $\eqref{2.10}_1$ and \eqref{2.11}--\eqref{2.12} and integration by parts that
\begin{align}\label{2.17}
	&\int_0^{\frac{M}{\omega_n}}\big(p(\rho)-\v\r^2(r^{n-1}u)_x\big)_xu\,r^{n-1} \dd x\\
	&=-\int_0^{\frac{M}{\omega_n}}\big(p(\rho)-\v\r^2(r^{n-1}u)_x\big)(r^{n-1}u)_x\,\dd x\nonumber \\
	&=-\int_0^{\frac{M}{\omega_n}}p(\rho)(r^{n-1}u)_x\, \dd x
	+\v\int_0^{\frac{M}{\omega_n}}\r^2\big((r^{n-1}u)_x\big)^2\,\dd x\nonumber \\
	&=a_0\int_0^{\frac{M}{\omega_n}}\rho^{\gamma-2}\rho_\tau\,\dd x
	+\v\int_0^{\frac{M}{\omega_n}}\r^2\big(r^{n-1}u_x+(n-1)r^{n-2}r_xu\big)^2\,\dd x\nonumber\\
	&=\frac{\dd}{\dd\tau}\int_0^{\frac{M}{\omega_n}} e(\rho)\,\dd x\nonumber\\
    &\quad +\v\int_0^{\frac{M}{\omega_n}}\Big(\r^2(r^{n-1}u_x)^2
	+(n-1)^2\frac{u^2}{r^2}+2(n-1)r^{n-2}\r u u_x\Big)\dd x.\nonumber
\end{align}
For the first term of \eqref{2.16}-RHS ({\it i.e.}, the right-hand side of \eqref{2.16}),
a direct calculation shows that
\begin{align}\label{2.18}
	(n-1)\v\int_0^{\frac{M}{\omega_n}}\r_xu^2r^{n-2} \dd x
	&=-(n-1)\v\int_0^{\frac{M}{\omega_n}}\big(2 r^{n-2}\rho uu_x +(n-2)\frac{ u^2}{r^2}\big)\,\dd x\\
	&\quad\, +(n-1)\v(\r u^2 r^{n-2})(\tau,\frac{M}{\omega_n}).\nonumber
\end{align}
For the last term of \eqref{2.16}-RHS, it follows from \eqref{2.12} that
\begin{align}\label{2.19}
	-\kappa \int_0^{\frac{M}{\omega_n}}\frac{x }{r^{n-1}} u \,\dd x
	=-\kappa \int_0^{\frac{M}{\omega_n}}\frac{x }{r^{n-1}} r_\tau\,\dd x
	=\frac{\kappa}{n-2}\frac{\dd}{\dd\tau} \int_0^{\frac{M}{\omega_n}} \frac{x}{r^{n-2}}\,\dd x.
\end{align}
Substituting \eqref{2.17}--\eqref{2.19} into \eqref{2.16}, we have
\begin{align} \label{2.20}
	&\frac{\dd}{\dd\tau}\Big\{\int_0^{\frac{M}{\omega_n}}\big(\frac{1}{2}u^2+ e(\rho)\big)\,\dd x
	-\frac{\kappa}{n-2}\int_0^{\frac{M}{\omega_n}} \frac{x}{r^{n-2}}\,\dd x\Big\}\\
	&\quad +\v\int_0^{\frac{M}{\omega_n}}\big(\r^2(r^{n-1}u_x)^2 +(n-1)\frac{u^2}{r^2}\big)\,\dd x
	+(n-1)\v(\r u^2 r^{n-2})(\tau,\frac{M}{\omega_n})\nonumber\\
    &=0.\nonumber
\end{align}
Plugging \eqref{2.20} back to the Eulerian coordinates, we obtain
\begin{align} \label{2.20-1}
&\frac{\dd}{\dd t}\Big\{\int_a^{b(t)} \big(\frac{1}{2} \rho u^2 + \rho e(\rho)\big)\,r^{n-1}\dd r
-\frac{\kappa}{n-2}\int_a^{b(t)}  \big(\int_a^r \rho\, z^{n-1} \dd z\big)\rho r\dd r\Big\}\\
&\quad+\v\int_a^{b(t)} \big(\rho u_r^2+(n-1)\rho\f{u^2}{r^2}\big)\, r^{n-1}  \dd r
+(n-1)\v (\r u^2)(t,b(t))b(t)^{n-2}\nonumber\\
&=0. \nonumber
\end{align}
Then integrating \eqref{2.20-1} over $[0,t]$ leads to
\begin{align}\label{2.22}
&\int_a^{b(t)}\big(\frac{1}{2} \rho u^2 + \rho e(\rho)\big)\,r^{n-1} \dd r
-\frac{\kappa}{n-2}\int_a^{b(t)}  \big(\int_a^r \rho(t,z)\,z^{n-1} \dd z\big)\rho\, r \dd r\\
&\quad +\v\int_0^t\int_a^{b(s)} \big(\rho u_r^2+(n-1) \rho \f{u^2}{r^2}\big)\,r^{n-1}\dd r\dd s\nonumber\\
&\quad +(n-1) \v\int_0^t (\r u^2)(s,b(s)) b(s)^{n-2}\, \dd s\nonumber\\
&=\int_a^{b}\big(\frac{1}{2}\rho_0 u_0^2 + \rho_0 e(\rho_0)\big)\,r^{n-1} \dd r
-\frac{\kappa}{n-2}\int_a^{b}  \big(\int_a^r \rho_0(z)\,z^{n-1} \dd z\big)\rho_0(r)\, r \dd r.\nonumber	
\end{align}

\smallskip
2. To close the estimates, we need to control the terms involving potential $\Phi$.
Noting \eqref{2.24-1}, a direct calculation shows that
\begin{align}\label{2.23-1}
&\frac{\kappa}{n-2}\int_a^{b(t)} \big(\int_a^r \rho\, z^{n-1} \dd z\big) \rho \,r \dd r \\
& =\frac{1}{(n-2)\kappa}\int_a^{b(t)}  (r^{n-1}\Phi_r)_r \Phi_r \,r \dd r\nonumber\\
&=\frac{1}{2\kappa}\Big\{\int_a^{b(t)}  |\Phi_r|^2\,r^{n-1}\dd r+\frac{1}{(n-2)b(t)^{n-2}}\big(\int_a^{b(t)} \rho\,z^{n-1} \dd z\big)^2\Big\}\nonumber\\
&= \frac{1}{2\kappa}\Big\{\int_a^{b(t)}  |\Phi_r|^2\,r^{n-1}\dd r+\left(\frac{M}{\omega_n}\right)^2\frac{1}{(n-2)b(t)^{n-2}} \Big\}.\nonumber
\end{align}
On the other hand, it follows from \eqref{2.24-1} that
\begin{align}
	\|\nabla \Phi\|_{L^2(\mathbb{R}^n)}^2
	&=\omega_n \Big\{\int_a^{b(t)} |\Phi_r|^2 \, r^{n-1}\dd r+ \int_{b(t)}^\infty  |\Phi_r|^2\, r^{n-1}\dd r \Big\}\nonumber\\
	&=\omega_n \Big\{\int_a^{b(t)}  |\Phi_r|^2\, r^{n-1}\dd r+ \left(\frac{M}{\omega_n}\right)^2 \int_{b(t)}^\infty r^{-n+1}\dd r\Big\}\nonumber\\
	&=\omega_n \Big\{\int_a^{b(t)} |\Phi_r|^2 \,r^{n-1}\dd r+ \left(\frac{M}{\omega_n}\right)^2\frac{1}{(n-2)b(t)^{n-2}}\Big\},\nonumber
\end{align}
which, together with \eqref{2.23-1}, yields that
\begin{align}\label{2.23}
		\frac{\kappa}{n-2}\int_a^{b(t)} \big(\int_a^r \rho z^{n-1} \dd z\big) \rho\,r \dd r
	&= \frac{1}{2\kappa\omega_n} \|\nabla \Phi\|^2_{L^2(\mathbb{R}^n)}\\
	&=\frac{1}{2\kappa} \int_a^\infty \frac{1}{r^{n-1}} \big(\int_a^r \rho\, z^{n-1} \dd z\big)^2 \dd r,\nonumber
\end{align}
where we need to understand $\rho$ to be zero  for $r\in[0,a)\cup (b(t),\infty)$ in the last equality of \eqref{2.23}.

\smallskip
3. Substituting \eqref{2.23} into \eqref{2.22}, we conclude \eqref{2.13-1}.
When $\kappa=-1$ (plasmas), \eqref{2.13} follows directly from \eqref{2.13-1}.

\smallskip
4. When $\kappa=1$ (gaseous stars), from \eqref{2.22} and \eqref{2.23},
the gravitational potential part has to be carefully controlled.
Multiplying \eqref{2.24} by $\Phi$ and integrating by parts yield
\begin{align}\label{2.25}
	\|\nabla \Phi\|_{L^2(\mathbb{R}^n)}^2&\leq \|\Phi\|_{L^{\frac{2n}{n-2}}(\mathbb{R}^n)}\,\|\rho\|_{L^{\frac{2n}{n+2}}(\Omega_t)}
	\leq\sqrt{A_n} \|\nabla \Phi\|_{L^2(\mathbb{R}^n)}\, \|\rho\|_{L^{\frac{2n}{n+2}}(\Omega_t)},
\end{align}
where
the positive constant $A_n=\frac{4}{n(n-2)}\omega_{n+1}^{-\frac{2}{n}}>0$ is the sharp constant for the Sobolev inequality which is given in Lemma \ref{Sobolev}.
Then it follows directly from \eqref{2.25} that
\begin{align}\label{2.26}
\|\nabla \Phi\|_{L^2(\mathbb{R}^n)}^2&\leq A_n \|\rho\|_{L^{\frac{2n}{n+2}}(\Omega_t)}^2\\
	&\leq A_n \|\rho\|_{L^1(\Omega_t)}^{2(1-\vartheta)}\, \|\rho\|_{L^\gamma(\Omega_t)}^{2\vartheta}\nonumber\\
	&\leq A_n \omega_n^{\frac{n-2}{n(\gamma-1)}}M^{\frac{(n+2)\gamma-2n}{n(\gamma-1)}}
	\Big(\int_a^{b(t)} \rho^{\gamma}\, r^{n-1} \dd r\Big)^{\frac{n-2}{n(\gamma-1)}},\nonumber
\end{align}
where $\vartheta=\frac{(n-2)\gamma}{2n(\gamma-1)}$,
and we have used the condition that $\gamma\geq \frac{2n}{n+2}$.

Substituting \eqref{2.26} into \eqref{2.23} and using \eqref{A.2}, we have
\begin{align}\label{2.29}
	&\frac{\kappa}{n-2}\int_a^{b(t)}  \big(\int_a^r \rho\, z^{n-1} \dd z\big)\rho \,r \dd r\\
	&\leq \frac{2\kappa}{n(n-2)}\omega_{n+1}^{-\frac{2}{n}}\omega_n^{\frac{2(n-1)-n\gamma}{n(\gamma-1)}}
	\Big(\frac{\gamma-1}{a_0}\Big)^{\frac{n-2}{n(\gamma-1)}}M^{\frac{(n+2)\gamma-2n}{n(\gamma-1)}}
	\, \Big(\int_a^{b(t)} \rho e(\rho)\, r^{n-1} \dd r\Big)^{\frac{n-2}{n(\gamma-1)}}\nonumber\\
	&:=B_{n,\g} M^{\frac{(n+2)\gamma-2n}{n(\gamma-1)}} \Big(\int_a^{b(t)} \rho e(\rho)\,r^{n-1} \dd r\Big)^{\frac{n-2}{n(\gamma-1)}},\nonumber
\end{align}
where $B_{n,\g}$ is the constant defined in \eqref{1.16}.

Noting \eqref{2.29}, we use the internal energy to control the gravitational part.
It follows from \eqref{2.29} that
\begin{align}\label{2.30}
&\int_a^{b(t)}\rho e(\rho)\,r^{n-1} \dd r
	-\frac{\kappa}{n-2}\int_a^{b(t)}  \Big(\int_a^r \rho \,z^{n-1} \dd z\Big)\rho\, r\dd r\\
&\geq \int_a^{b(t)}\rho e(\rho)\,r^{n-1} \dd r
	-B_{n,\g} M^{\frac{(n+2)\gamma-2n}{n(\gamma-1)}} \Big(\int_a^{b(t)} \rho e(\rho)\, r^{n-1} \dd r\Big)^{\frac{n-2}{n(\gamma-1)}}.\nonumber
\end{align}

For the case that $\kappa=1$ with $\gamma>\frac{2(n-1)}{n}$, which implies that $\frac{n-2}{n(\gamma-1)}<1$.
Then it follows from \eqref{2.30} and the H\"older inequality that
\begin{align*}
	&\int_a^{b(t)}\rho e(\rho)\,r^{n-1} \dd r
	-\frac{\kappa}{n-2}\int_a^{b(t)}\Big(\int_a^r \rho \,z^{n-1} \dd z\Big)\rho\, r\dd r\\
	&\geq \frac{1}{2} \int_a^{b(t)}\rho e(\rho)\, r^{n-1} \dd r-C(n,M),\nonumber
\end{align*}
which, together with \eqref{2.22}, yields \eqref{2.15}.

For the case that $\kappa=1$ with $\gamma=\frac{2(n-1)}{n}$, {\it i.e.}, $\frac{n-2}{n(\gamma-1)}=1$, we use \eqref{2.22}  and  \eqref{2.30} to obtain
\begin{align}\label{2.31-1}
&\int_a^{b(t)}\rho e(\rho)\,r^{n-1}\dd r
	-\frac{\kappa}{n-2}\int_a^{b(t)}\Big(\int_a^r \rho\, z^{n-1} \dd z\Big)\rho\,r \dd r\\
&\geq \Big(1- B_{n,\g} M^{\frac{2}{n}}
	\Big)\int_a^{b(t)}\rho e(\rho)\,r^{n-1} \dd r
	:=C_{\gamma}\int_a^{b(t)}\rho e(\rho)\,r^{n-1}\dd r,\nonumber
\end{align}
provided that
$M<M_{\rm c}^{\v,b}(\gamma):=B_{n,\g}^{-\frac{n}2}$.

For the case that $\kappa=1$ with $\gamma\in(\frac{2n}{n+2}, \frac{2(n-1)}{n})$, we define
\begin{equation*}
	F(s)=s-B_{n,\g} M^{\frac{(n+2)\gamma-2n}{n(\gamma-1)}} \, s^{\frac{n-2}{n(\gamma-1)}}\qquad \mbox{for $s\geq0$}.
\end{equation*}
A direct calculation shows that
\begin{align*}
	\begin{cases}
		\displaystyle F'(s)=1-\frac{n-2}{n(\gamma-1)}B_{n,\g} M^{\frac{(n+2)\gamma-2n}{n(\gamma-1)}} \, s^{\frac{2(n-1)-n\gamma}{n(\gamma-1)}},\\[4mm]
		\displaystyle  F''(s)=-\frac{(2(n-1)-n\gamma)(n-2)}{n^2(\gamma-1)^2}\,
		B_{n,\g} M^{\frac{(n+2)\gamma-2n}{n(\gamma-1)}}\,s^{\frac{3n-2-2n\gamma}{n(\gamma-1)}},
	\end{cases}
\end{align*}
which yields that $F''(s)<0$ for any $s>0$ if $\gamma<\frac{2(n-1)}{n}$, so that $F(s)$ is concave for $s>0$.
We denote
\begin{align}\label{2.33}
	s_\ast = \Big(\frac{(n-2)B_{n,\g}}{n(\gamma-1)}\Big)^{-\frac{n(\gamma-1)}{2(n-1)-n\gamma}}\, M^{-\frac{(n+2)\gamma-2n}{2(n-1)-n\gamma}},
\end{align}
which is the critical point of $F(s)$ so that $F'(s_\ast)=0$.
The maximum of $F(s)$ for $s\geq0$ is
\begin{align}\nonumber
	F(s_\ast)=\frac{2(n-1)-n\gamma}{n-2} \Big(\frac{(n-2)B_{n,\g}}{n(\gamma-1)}\Big)^{-\frac{n(\gamma-1)}{2(n-1)-n\gamma}}
	\, M^{-\frac{(n+2)\gamma-2n}{2(n-1)-n\gamma}}>0  \qquad \mbox{for $\gamma<\frac{2(n-1)}{n}$}.
\end{align}

Now we claim that, under condition \eqref{2.60-3},
\begin{equation}\label{2.38}
	\int_a^{b(t)} \rho e(\rho)\, r^{n-1} \dd r<s_\ast.
\end{equation}
Noting that $\frac{2(n-1)-n\gamma}{(n+2)\gamma-2n}>0$ for $\gamma\in(\frac{2n}{n+2}, \frac{2(n-1)}{n})$,
it follows from \eqref{2.36}, \eqref{2.60-3}, and \eqref{2.33} that
\begin{equation}\label{2.40}
	F(s_\ast)>\frac{E_0^{\v,b}}{\omega_n},
\end{equation}
and
\begin{align}\label{2.41}
	s_\ast &>\Big(\frac{(n-2)B_{n,\g}}{n(\gamma-1)}\Big)^{-\frac{n(\gamma-1)}{2(n-1)-n\gamma}}\,
	M_c^{\v,b}(\gamma)^{-\frac{(n+2)\gamma-2n}{2(n-1)-n\gamma}} \\
	&=\frac{n-2}{2(n-1)-n\gamma} \frac{E_0^{\v,b}}{\omega_n}
	>\frac{2 E_0^{\v,b}}{\omega_n}, \nonumber
\end{align}
where we have used $\frac{n-2}{2(n-1)-n\gamma}>1+\frac{n}{2}>2$ for $\gamma\in (\frac{2n}{n+2}, \frac{2(n-1)}{n})$.
Then it follows from \eqref{2.22}, \eqref{2.30}, and \eqref{2.40}--\eqref{2.41}  that
\begin{align}
	&F(\int_a^{b(t)} \rho e(\rho)\, r^{n-1} \dd r)\leq \frac{E_0^{\v,b}}{\omega_n}<F(s_\ast),\label{2.42}\\
	&\int_a^{b} \rho_0 e(\rho_0)\, r^{n-1} \dd r\leq  \frac{E_0^{\v,b}}{\omega_n}<s_\ast.\label{2.43}
\end{align}
Thus, due to the continuity of $\int_a^{b(t)} (\rho e(\rho))(t,r)\,r^{n-1} \dd r$ with respect to $t$,
\eqref{2.38} must hold.
Otherwise,
there exists some time $t_0>0$ such that $\int_a^{b(t)} \rho e(\rho)\, r^{n-1} \dd r=s_\ast$, which yields
\begin{align}\nonumber
	F(\int_a^{b(t_0)} \rho e(\rho) \,r^{n-1} \dd r)=F(s_\ast)>\frac{E_0^{\v,b}}{\omega_n},
\end{align}
which contradicts \eqref{2.42}. Therefore, \eqref{2.38} always holds under condition \eqref{2.60-3}.

Now, under condition \eqref{2.60-3}, it follows from \eqref{2.38} that
\begin{align}\label{2.39}
	&F(\int_a^{b(t)} \rho e(\rho)\, r^{n-1} \dd r)\\
	&\geq\Big(1-B_{n,\g} M^{\frac{(n+2)\gamma-2n}{n(\gamma-1)}}\, s_\ast^{\frac{2(n-1)-n\gamma}{n(\gamma-1)}}\Big)
	\int_a^{b(t)} \rho e(\rho)\, r^{n-1} \dd r \nonumber\\
	&= \frac{2(n-1)-n\gamma}{n-2} \int_a^{b(t)} \rho e(\rho)\, r^{n-1} \dd r.\nonumber
\end{align}
Thus, \eqref{2.14} follows directly from \eqref{2.31-1} and \eqref{2.39}.
This completes the proof.
\end{proof}

Using \eqref{2.24-1},  \eqref{2.13}--\eqref{2.14}, \eqref{2.15}, \eqref{2.23}, and \eqref{A.1},
we have the following estimates for the potential function $\Phi$.

\begin{corollary}\label{cor3.1}
	Under the conditions of Lemma {\rm \ref{lem2.1}},
	\begin{align*}
		&|r^{n-1} \Phi_r(t,r)|\leq \frac{M}{\omega_n}\qquad\mbox{for $\, (t,r)\in [0,\infty)\times[0,\infty)$},\\
		&\int_a^{b(t)}\Big(\int_a^r \rho(t,z)\,z^{n-1} \dd z\Big)\rho(t,r)\,r \dd r
		+\|\Phi(t)\|_{L^{\frac{2n}{n-2}}(\mathbb{R}^n)}+\|\nabla \Phi(t)\|_{L^2(\mathbb{R}^n)}\\[1.5mm]
		&\quad\leq C(M,E_0)\quad\mbox{for $t\geq0$}.
	\end{align*}
\end{corollary}

For later use, we analyze the behavior of density $\rho$  on the free boundary. It follows from $\eqref{2.10}_1$ and \eqref{2.11} that
\begin{align}\label{2.44-1}
	\rho_{\tau}(\tau,\frac{M}{\omega_n})=-\frac{a_0}{\v}\r^{\g}(\tau,\frac{M}{\omega_n}).
\end{align}
Then we obtain
\begin{align}\nonumber
	\rho(\tau,\frac{M}{\omega_n})=\rho_0(\frac{M}{\omega_n})\,\Big(1+\frac{(\g-1)^2}{\v}e(\r_0(\frac{M}{\omega_n}))\tau\Big)^{-\frac{1}{\g-1}}.
\end{align}
In the Eulerian coordinates, it is equivalent to the form:
\begin{align}\label{2.44}
	\rho(t,b(t))=\rho_0(b)\Big(1+\frac{(\g-1)^2}{\v}e(\r_0(b))t\Big)^{-\frac{1}{\g-1}}\leq \rho_0(b).
\end{align}
The density behavior on the free boundary \eqref{2.44} is important, which will be used frequently later.

\begin{lemma}[BD-type entropy estimate]  \label{lem2.2}
	Under the conditions of Lemma {\rm \ref{lem2.1}}, for any given $T>0$, the following holds for any $t\in[0,T]${\rm :}
	\begin{align}\label{2.45}
     &\v^2\int_a^{b(t)}\big|\big(\sqrt{\rho(t,r)}\big)_r\big|^2\,r^{n-1} \dd r+\frac{3a_0\v}{\gamma} \int_0^t\int_a^{b(s)} \big|\big(\rho^{\frac{\gamma}{2}}\big)_r\big|^2\,r^{n-1}\dd r\dd s\\
     &\quad +\frac{1}{n}p(\rho(t,b(t)))b(t)^n+\frac{1}{n\v}\int_0^t p(\rho(s,b(s)))p'(\rho(s,b(s)))b(s)^n\,\dd s\nonumber\\[1.5mm]
     &\leq  C(E_0,M,b^n\rho_0^\gamma(b),T) \leq C(E_0,M,T).\nonumber
	\end{align}
\end{lemma}

\begin{proof}  We divide the proof into  four steps.

1. For convenience, we start with the solution in the Lagrangian coordinates $(\tau, x)$.
It follows from $\eqref{2.10}_{1}$ that
\begin{align*}
	\rho_{x\tau}=-\big(\rho^2(r^{n-1}u)_x\big)_x,
\end{align*}
which, together with  $\eqref{2.10}_{2}$, yields that
\begin{align}\label{2.47}
	u_{\tau}+r^{n-1}p_x=-\v r^{n-1}\rho_{x\tau}-\v(n-1) r^{n-2}u\rho_x -\kappa\frac{x}{r^{n-1}}.
\end{align}
Then \eqref{2.47} can be rewritten by using \eqref{2.12} as
\begin{align}\label{2.49}
	(u+\v r^{n-1}\r_x)_{\tau}+r^{n-1}p_x=-\kappa\frac{x}{r^{n-1}}.
\end{align}
Multiplying \eqref{2.49} by $u+\v r^{n-1}\rho_x$ yields
\begin{align}\label{2.48}
	&\frac{1}{2}\frac{\dd}{\dd\tau}\int_0^{\frac{M}{\omega_n}}(u+\v r^{n-1}\r_x)^2\,\dd x
	+\v\int_0^{\frac{M}{\omega_n}}p(\rho)_x\rho_x\, r^{2n-2}\dd x\\
    &\quad +\int_0^{\frac{M}{\omega_n}}p(\rho)_x u\,r^{n-1}\dd x\nonumber\\
	&=\frac{\kappa}{n-2}\frac{\dd}{\dd\tau} \int_0^{\frac{M}{\omega_n}} \frac{x}{r^{n-2}}\,\dd x
	-\kappa\v\int_0^{\frac{M}{\omega_n}}\rho_x x\,\dd x,\nonumber
\end{align}
where we have used \eqref{2.19}.
Using $\eqref{2.10}_1$, \eqref{2.12}, and \eqref{2.44-1}, we have
\begin{align}
&\int_0^{\frac{M}{\omega_n}}p(\rho)_xu\,r^{n-1}\dd x \label{2.50} \\
	&=-\int_0^{\frac{M}{\omega_n}}p(\rho)(r^{n-1}u)_x\,\dd x
	+( p(\rho)u r^{n-1})(\tau,\frac{M}{\omega_n})\nonumber\\
&=\frac{\dd}{\dd\tau}\int_0^{\frac{M}{\omega_n}}e(\rho)\,\dd x
	+( p(\rho)r^{n-1})(\tau,\frac{M}{\omega_n})\, r_\tau(\tau,\frac{M}{\omega_n})\nonumber\\
&=\frac{\dd}{\dd\tau}\int_0^{\frac{M}{\omega_n}}e(\rho)\,\dd x
	+\Big(\frac{1}{n}p(\rho(\tau,\frac{M}{\omega_n}))b(\tau)^n\Big)_\tau\nonumber\\
&\quad\,	-\frac{1}{n}p'(\rho(\tau,\frac{M}{\omega_n}))\rho_{\tau}(\tau,\frac{M}{\omega_n})b(\tau)^n\nonumber\\
&=\frac{\dd}{\dd\tau}\int_0^{\frac{M}{\omega_n}}e(\rho)\, \dd x
	+\Big(\frac{1}{n}p(\rho(\tau,\frac{M}{\omega_n}))b(\tau)^n\Big)_\tau\nonumber\\
&\quad\,	+\frac{1}{n\v}p(\rho(\tau,\frac{M}{\omega_n}))p'(\rho(\tau,\frac{M}{\omega_n}))b(\tau)^n,\nonumber\\[3mm]
&\int_0^{\frac{M}{\omega_n}}x\rho_x\, \dd x
	=-\int_0^{\frac{M}{\omega_n}} \rho\,\dd x +\frac{M}{\omega_n} \rho(\tau,\frac{M}{\omega_n}). \label{2.51}
\end{align}
Substituting \eqref{2.50}--\eqref{2.51} into \eqref{2.48} yields
\begin{align}\label{2.52}
	&\frac{\dd}{\dd\tau}\Big\{\int_0^{\frac{M}{\omega_n}}\frac{1}{2}(u+\v r^{n-1}\r_x)^2\,\dd x
	+\int_0^{\frac{M}{\omega_n}}e(\rho)\,\dd x-\frac{\kappa}{n-2}\int_0^{\frac{M}{\omega_n}} \frac{x}{r^{n-2}}\,\dd x\Big\}\\
	&\quad\, +\v \int_0^{\frac{M}{\omega_n}}p'(\r)\r_x^2\,r^{2n-2} \dd x
	+\Big(\frac{1}{n}p(\rho(\tau,\frac{M}{\omega_n}))b(\tau)^n\Big)_\tau\nonumber\\
	&\quad\, +\frac{1}{n\v}p(\rho(\tau,\frac{M}{\omega_n}))p'(\rho(\tau,\frac{M}{\omega_n})) b(\tau)^n\nonumber \\
	&=\v\kappa \int_0^{\frac{M}{\omega_n}} \rho\,\dd x-\v \kappa\frac{M}{\omega_n} \rho(\tau,\frac{M}{\omega_n}).\nonumber
\end{align}
Integrating \eqref{2.52} over $[0,\tau]$, we have
\begin{align}\label{2.53}
&\int_0^{\frac{M}{\omega_n}}\frac{1}{2}(u+\v r^{n-1}\r_x)^2\,\dd x
+\int_0^{\frac{M}{\omega_n}}e(\rho)\,\dd x-\frac{\kappa}{n-2}\int_0^{\frac{M}{\omega_n}} \frac{x}{r^{n-2}}\,\dd x\\
&\quad+ \v \int_0^\tau\int_0^{\frac{M}{\omega_n}} p(\r) \r_x^2\,r^{2n-2}\dd x\dd s
+\frac{1}{n}p(\rho(\tau,\frac{M}{\omega_n}))b(\tau)^n\nonumber\\
&\quad +\frac{1}{n\v}\int_0^\tau p(\rho(\tau,\frac{M}{\omega_n}))p'(\rho(\tau,\frac{M}{\omega_n})) b(s)^n\,\dd s\nonumber\\
&=\int_0^{\frac{M}{\omega_n}}\frac{1}{2}(u_0+\v r_0^{n-1}\rho_{0x})^2\,\dd x+\int_0^{\frac{M}{\omega_n}}e(\rho_0)\,\dd x
-\frac{\kappa}{n-2}\int_0^{\frac{M}{\omega_n}} \frac{x}{r_0^{n-2}}\, \dd x\nonumber\\
&\quad+\frac{1}{n}b^n p(\rho_0(\frac{M}{\omega_n}))+\v\kappa\int_0^\tau \int_0^{\frac{M}{\omega_n}} \rho\,\dd x\dd s
-\v\kappa\frac{M}{\omega_n} \int_0^\tau\rho(s,\frac{M}{\omega_n})\,\dd s.\nonumber	
\end{align}
Plugging \eqref{2.53} back to the Eulerian coordinates and using Lemma \ref{lem2.1}, we obtain
\begin{align}\label{2.54}
	&\frac{\v^2}{4}\int_a^{b(t)} \big|\big(\sqrt{\rho(t,r)}\big)_r\big|^2\,r^{n-1} \dd r
	+\frac{4 a_0\v}{\gamma} \int_0^t\int_a^{b(s)} \big|\big(\rho^{\frac{\gamma}{2}}\big)_r\big|^2\,r^{n-1}\dd r\dd s\\
	&\quad+\frac{1}{n}p(\rho(t,b(t)))b(t)^n
	+\frac{1}{n\v}\int_0^t p(\rho(s,b(s)))p'(\rho(s,b(s))) b(s)^n\,\dd s\nonumber\\
	&\leq C(E_0,M)+\frac{1}{n}p(\rho_0(b))b^n +\v\kappa\int_0^t \int_a^{b(s)} \rho^2\, r^{n-1} \dd r\dd s\nonumber\\
	&\quad -\v\kappa \frac{M}{\omega_n} \int_0^t\rho(s,b(s))\,\dd s.\nonumber
\end{align}

2. For the second term of \eqref{2.54}-RHS, it follows from \eqref{2.60-1} that
\begin{equation}\label{2.60-2}
	\frac{1}{n}p(\rho_0(b))b^n\leq C.
\end{equation}
For the last term of \eqref{2.54}-RHS, by using \eqref{2.44}, we have
\begin{align}\label{2.55}
	\left|\v\kappa \frac{M}{\omega_n} \int_0^t\rho(s,b(s))\,\dd s\right|\leq C(M) \rho_0(b) T\leq C(M,T).
\end{align}

3. For $\kappa=-1$ (plasmas), then \eqref{2.45} follows from \eqref{2.54} and \eqref{2.55}.

\smallskip
4. To close the estimates for $\kappa=1$ (gaseous stars), we still need to bound the third term of \eqref{2.54}-RHS:
\begin{equation}\nonumber
	\v\kappa\int_0^t \int_a^{b(s)} \rho^2\, r^{n-1} \dd r\dd s=\frac{\v\kappa}{\omega_n}\int_0^t \|\rho(s)\|^2_{L^2(\Omega_s)} \dd s.
\end{equation}
We estimate the above term in the following two cases:

\smallskip
\noindent{\it Case} 1. For $\gamma\geq 2$, then it is bounded as
\begin{align}\label{2.56}
	\v\kappa\int_0^t \int_a^{b(s)} \rho^2 \, r^{n-1}\dd r\dd s
	&\leq C\v\int_0^t \int_a^{b(s)} \rho\big(1+ e(\rho)\big)\,r^{n-1} \dd r\dd s\\[2mm]
    &\leq C(E_0,M,T). \nonumber
\end{align}

\smallskip
\noindent{\it Case} 2. For $\gamma\in (\frac{2n}{n+2}, 2)$, we notice that $\frac{n\gamma}{n-2}>2$ and use the interpolation inequality to obtain
\begin{align}\label{2.57}
	\|\rho(t)\|_{L^2(\Omega_t)}\leq \|\rho(t)\|_{L^{\frac{n\gamma}{n-2}}(\Omega_t)}^{\bar{\vartheta}}\, \|\rho(t)\|_{L^{\gamma}(\Omega_t)}^{1-\bar{\vartheta}}
	\qquad \mbox{with $\bar{\vartheta}=\frac{n(2-\gamma)}{4}$}.
\end{align}
For $B_R(\mathbf{0})\subset \mathbb{R}^n$, the following Sobolev inequality holds:
\begin{align} \label{2.27-1}
	\|f\|_{L^{\frac{2n}{n-2}}(B_R(\mathbf{0}))}
	\leq C\big(\|\nabla f\|_{L^2(B_R(\mathbf{0}))}+\frac1R \|f\|_{L^2(B_R(\mathbf{0}))}\big).
\end{align}
It follows from \eqref{2.8} that
\begin{align}
	\frac{M}{\omega_n}&=\int_a^{b(t)} \rho(t,r)\, r^{n-1} \dd r\nonumber\\
	&\leq  \Big(\int_a^{b(t)} \rho^{\gamma}\, r^{n-1} \dd r\Big)^{\frac1\gamma} \Big(\int_a^{b(t)} r^{n-1} \dd r\Big)^{1-\frac1\gamma}\nonumber\\
	&\leq n^{\frac1\gamma-1} b(t)^{n(1-\frac1\gamma)} \Big(\int_a^{b(t)} \rho^{\gamma}\,r^{n-1} \dd r\Big)^{\frac1\gamma},\nonumber
\end{align}
which yields that
\begin{align}\label{2.27}
	b(t)^{-1}\leq \Big(n^{\frac{1-\gamma}{\gamma}} \,\frac{\omega_n}{M}\Big)^{\frac{\gamma}{n(\gamma-1)}}
	\Big(\int_a^{b(t)} \rho^{\gamma}\, r^{n-1} \dd r\Big)^{\frac{1}{n(\gamma-1)}}\leq C(M,E_0).
\end{align}
This, together with Lemma \ref{lem2.1} and \eqref{2.27-1},  yields that
\begin{align}\label{2.58}
	\|\rho(t)\|_{L^{\frac{n\gamma}{n-2}}(\Omega_t)}
	&=\Big(\|\rho^{\frac{\gamma}{2}}(t)\|_{L^{\frac{2n}{n-2}}(\Omega_t)}\Big)^{\frac2{\gamma}}\\
	&\leq C\Big(\|\nabla (\rho^{\frac{\gamma}{2}})\|_{L^2(\Omega_t)}+ b(t)^{-1}\|\rho^{\frac{\gamma}{2}}\|_{L^2(\Omega_t)}\Big)^{\frac{2}{\gamma}}\nonumber\\
	&\leq C(M,E_0) \Big(1+\big(\int_a^{b(t)}|(\rho^{\frac\gamma2})_r|^2r^{n-1}\dd r\big)^{\frac{1}{\gamma}}\Big).\nonumber
\end{align}
Substituting \eqref{2.58} into \eqref{2.57} and using \eqref{lem2.1} and the H\"older inequality, we have
\begin{align}\label{2.59}
	&\v\kappa\int_0^t \int_a^{b(s)} \rho^2\, r^{n-1} \dd r\dd s\\
	&\leq C(M,E_0) \int_0^t\v\Big(1+\big(\int_a^{b(s)}|(\rho^{\frac\gamma2})_r|^2r^{n-1}\dd r\big)^{\frac{1}{\gamma}}\Big)^{2\bar{\vartheta}}\dd s\nonumber\\
	&\leq C(M,E_0,T) +\frac{\v a_0}{\gamma}\int_0^t\int_a^{b(s)}|(\rho^{\frac\gamma2})_r|^2\,r^{n-1}\dd r\dd s,\nonumber
\end{align}
where we have used $\dis \frac{2\bar{\vartheta}}{\gamma}\in(0,1)$ for $\gamma>\frac{2n}{n+2}$.
Finally, combining \eqref{2.54}--\eqref{2.60-2}, \eqref{2.56}, and \eqref{2.59}, we obtain \eqref{2.45}.
This completes the proof.
\end{proof}

From \eqref{2.27}, we know that $b(t)$ has a uniform positive lower bound.
However, to take limit $b\rightarrow\infty$,
we need to make sure that domain $\Omega_T$ can expand to the whole physical space for fixed $\v>0$; that is,
$\displaystyle\inf_{t\in[0,T]}b(t)\longrightarrow\infty$ as $b\rightarrow\infty$.

\begin{lemma}[{Expanding of domain $\Omega_T$}]\label{lem2.3}
	Given $T>0$ and $\v\in(0,\v_0]$, there exists a positive constant $C_1(M,E_0,T,\v)>0$ such that,
	if $\dis b\geq \max\{C_1(M,E_0,T,\v),\mathfrak{B}(\v)\}$,
	\begin{align}\label{2.60}
		b(t)\geq \frac{b}{2}\qquad\,\, \mbox{for $t\in[0,T]$},
	\end{align}
	where $\mathfrak{B}(\v)$ is defined for \eqref{A.51}.
\end{lemma}

\begin{proof}
Noting the continuity of $b(t)$, we first make the {\it a priori } assumption:
\begin{equation}\label{2.63-1}
	b(t)\geq \frac{b}{4}.
\end{equation}
Integrating \eqref{2.6} over $[0,t]$ yields
\begin{align}\label{3.28}
	b(t)=b+\int_0^t u(s,b(s))\,\dd s.
\end{align}
A direct calculation by using \eqref{2.44}, \eqref{2.63-1}, and Lemma \ref{lem2.1} yields that
\begin{align}\label{2.63}
&\left|\int_0^tu(s,b(s))\,\dd s \right|\\
&\leq \frac{C}{\sqrt{\v}}\Big(\int_0^t\v (\rho u^2 r^{n-2})(s,b(s))\,\dd s\Big)^{\frac12}
\Big( \int_0^t\frac{1}{\rho (s,b(s))b(s)^{n-2}}\dd s  \Big)^{\frac12}\nonumber\\
&\leq \frac{C(M,E_0)}{\sqrt\v}\Big( \int_0^t\frac{1}{\rho (s,b(s))b(s)^{n-2}}\dd s  \Big)^{\frac12}\nonumber\\
&=\frac{C(M,E_0)}{\sqrt\v}
\bigg( \int_0^t\frac{\big(1+\frac{(\gamma-1)^2}{\v} e(\rho_0(b))s\big)^{\frac{1}{\gamma-1}}}{\rho_0(b)b(s)^{n-2}}
\dd s \bigg)^{\frac12}\nonumber\\
&\leq C_0(M,E_0)\,\Big(\frac{1+T}{\v}\Big)^{\frac{\g}{2(\gamma-1)}} \rho_{0}(b)^{-\frac12}b^{-\frac{n-2}{2}}.\nonumber
\end{align}
We take $C_1(M,E_0,T,\v):=\big(4 C_0(M,E_0) \big)^{\frac2\alpha}\left(\frac{1+T}{\v}\right)^{\frac{\g}{\alpha(\gamma-1)}}$.
Then we use \eqref{2.60-1} and \eqref{2.63} to conclude that
\begin{align} \label{3.29}
	\left|\int_0^tu(s,b(s))\,\dd s \right|&\leq \frac{b}{4},
\end{align}
provided $b\geq C_1(M,E_0,T,\v)$. Then it follows from  \eqref{3.28} and \eqref{3.29}  that
\begin{equation}\label{2.66-1}
	b(t)\geq \frac{3b}{4}.
\end{equation}
Thus, we have closed our {\it a priori} assumption \eqref{2.63-1}.
Then, using \eqref{2.66-1} and the continuity arguments, we conclude \eqref{2.60}.
\end{proof}

\begin{lemma}[Higher integrability on the density]\label{lem2.4}
	Let $(\r,u)$ be the smooth solution of \eqref{2.1}--\eqref{initial}.
	Then,  under the assumption of Lemma {\rm \ref{lem2.1}},
	\begin{align}\label{2.65}
		\int_0^T\int_K\r^{\g+1}(t,r)\,\dd r \dd t\leq C(K,M,E_0,T)
	\end{align}
for any $K\Subset (a,b(t))$ and any $t \in[0,T]$.
\end{lemma}

\begin{proof}  We divide the proof into five steps.

1.  For given $K\Subset (a,b(t))$ for any $t\in[0,T]$,
there exist $d$ and $D$ such that $K\Subset (d,D)\Subset[a,b(t)]$.
Let $w(r)$ be a smooth compact support function with $\text{supp}\,w\subseteq(d,D)$ and $w(r)=1$ for $r\in K$.
Multiplying $\eqref{2.1}_2$ by $w(r)$, we have
\begin{align}\label{2.66}
	&(\r u w)_t+\big((\r u^2+p(\rho))w\big)_r
	+\frac{n-1}{r} \r u^2 w\\
	&=\v\big(\rho(u_r+\frac{n-1}{r}u) w\big)_r-\v \frac{n-1}{r}u \rho_r w\nonumber\\
	&\quad +\big(\r u^2+p(\rho)-\v\rho(u_r+\frac{n-1}{r} u)\big)w_r
	-\frac{\kappa \rho w}{r^{n-1}} \int_a^r \rho\,z^{n-1}\dd z.\nonumber
\end{align}
Integrating \eqref{2.66} over $[d,r)$ to obtain
	\begin{align}\label{2.67}
		p(\rho)w
		&=\v \rho\big(u_r+\frac{n-1}{r}u\big) w -\v\int_d^r \frac{n-1}{z} u\rho_z w\, \dd z\\
		&\quad -\Big(\big(\int_d^r\rho u w\,\dd z\big)_t+ \rho u^2 w\Big) -\int_d^r \frac{n-1}{z} \rho u^2 w\,\dd z\nonumber \\
		&\quad+\int_d^r\big(\rho u^2+p(\rho)-\v\rho(u_z+\frac{n-1}{z} u)\big)w_z\, \dd z\nonumber\\
		&\quad -\kappa \int_d^r  \big(\int_a^z \rho\, y^{n-1}dy\big)\frac{\rho w}{z^{n-1}} \dd z.\nonumber
	\end{align}
	Multiplying \eqref{2.67} by $\rho w$, we have
	\begin{align}\label{2.68}
     \rho\,p(\rho)w^2
     	&=\v \rho^2\big(u_r+\frac{n-1}{r}u\big) w^2 -\v\rho w\int_d^r \frac{n-1}{z} u\rho_z w\, \dd z\\
     	&\quad-\Big(\rho w\big(\int_d^r\rho u w\,\dd z\big)_t+ \rho^2 u^2 w^2\Big) -\rho w\int_d^r \frac{n-1}{z} \rho u^2 w\,\dd z\nonumber\\
     	&\quad +\rho w\int_d^r\big(\rho u^2+p(\rho)-\v\rho(u_z+\frac{n-1}{z} u)\big)w_z\, \dd z\nonumber\\
     	&\quad-\kappa \rho w\int_d^r  \big(\int_a^z \rho\, y^{n-1}dy\big)\frac{\rho w}{z^{n-1}} \dd z.\nonumber
	\end{align}
	Notice that
	\begin{align}
		\r w\Big(\int_d^r\r u w\, \dd z\Big)_t+ \r^2 u^2 w^2
		&=\Big(\r w\int_d^r\r u w \,\dd z\Big)_t+\Big(\r  uw\int_d^r\r u w \,\dd z\Big)_r\nonumber\\
		&\quad-\r uw_r \int_d^r \r u w \, \dd z+\frac{n-1}{r}\r u w \int_d^r \r u w \,\dd z.\nonumber
	\end{align}
	Then it following from
	\eqref{2.68} that
\begin{align}\label{2.69}
	\r p(\rho) w^2&=\v \r^2\big(u_r+\frac{n-1}{r} u\big) w^2 -\v\r w\int_d^r \frac{n-1}{z} u\rho_z w\,\dd z\\
	&\quad -\Big(\r w\int_d^r\r u w\,\dd z\Big)_t -\Big(\r  uw\int_d^r\r u w \,\dd z\Big)_r\nonumber\\
	&\quad +\r uw_r \int_d^r \r u w\, \dd z-\frac{n-1}{r} \r u w \int_d^r \r u w \,\dd z\nonumber\\
	&\quad +\r w\int_d^r\Big(\r u^2+p(\rho)-\v\rho(u_z+\frac{n-1}{z}u)\Big)w_z\, \dd z\nonumber\\
	 &\quad -\r w\int_d^r \frac{n-1}{z}\r u^2 w\, \dd z -\kappa \rho w \int_d^r  \big(\int_a^z \rho\, y^{n-1}dy\big)\frac{\rho w}{z^{n-1}}\, \dd z \nonumber\\
    &:= \sum_{i=1}^9 I_i.\nonumber
\end{align}

2. To estimate $I_i, i=1,\cdots, 9$, in \eqref{2.69}, we first notice that
\begin{align}\label{2.70}
	\int_d^D(\rho+\rho|u|)\, \dd r\leq \f{C}{d^{n-1}}\int_d^D (\rho+\rho u^2)\,r^{n-1}\dd r\leq C(d, M, E_0).
\end{align}
Then it follows from \eqref{2.70} that
	\begin{align}\label{2.71}
     	\left| \int_0^T\int_d^D I_3\, \dd r\dd t \right|
     	&\leq \left| \int_d^D (\rho w)(T,r)\,\Big(\int_d^{r}(\rho u w)(T,z) \, \dd z\Big)\, \dd r \right|\\
     	&\quad + \left| \int_d^D (\rho w)(0,r)\, \Big(\int_d^{r}(\rho u w)(0,z) \, \dd z\Big)\, \dd r \right|\nonumber\\
     	&\leq C\,\sup_{t\in[0,T]}\left\{ \int_d^D \rho(t,r)\, \dd r\, \int_d^D (\rho |u|)(t,z)\, \dd z \right\}\nonumber\\
     	&\leq C(d,M,E_0,T),\nonumber
	\end{align}
	\begin{align}\label{2.71-2}
		&\left| \int_0^T\int_d^D \big(I_5+I_6+I_8\big)\, \dd r\dd t \right|\\
		&\leq C(d)\int_0^T\left(\int_d^D (\rho+\rho |u|^2)(t,r)\, \dd r\right)
		\left(\int_d^r (\rho+\rho |u|^2)(t,z)\, \dd z\right)\, \dd t\nonumber\\
		&\leq C(d,M,E_0,T),\nonumber
	\end{align}
	and
	\begin{align}\label{2.71-3}
     	&\left| \int_0^T\int_d^D I_9\, \dd r\dd t \right|\\
     	&\leq C(d)\int_0^T \left(\int_d^D\rho(t,r)\, \dd r\right)\left(\int_d^r\rho(t,z)\, \dd z\right)
     	\left(\int_d^z \rho(t,y)\, y^{n-1}\, \dd y\right) \dd t\nonumber\\
     	&\leq C(d)\int_0^T \left(\int_d^D \rho(t,r)\, r^{n-1}\, \dd r\right)^3
     	\dd t\nonumber\\
     	&\leq C(d,M,E_0,T).\nonumber
	\end{align}
	Since $\mbox{supp}\,w\subseteq (d,D)$, it is clear that
	\begin{align}\label{2.71-1}
		\left| \int_0^T\int_d^D I_4\, \dd r\dd t \right|=\left| \int_0^T\int_d^D  \Big(\r  uw\int_d^r\r u w \,\dd z\Big)_r\, \dd r\dd t \right|=0.
	\end{align}

3. Next, we estimate $I_7$. Noting that
\begin{align*}
	&\Big|\int_0^T\int_d^D \r w \int_d^r \big(\r u^2+p(\rho)\big) w_z\,\dd z \dd r \dd t\Big| \leq C(d,M,E_0,T),\\[3mm]
	&\v\Big|\int_0^T\int_d^D \r w \int_d^r \rho(u_z+\frac{n-1}{z}u) w_z\,\dd z \dd r \dd t \Big|\nonumber\\
	&\quad\leq C(d,M,E_0) \Big(\v\int_0^T\int_d^D \big(z^{n-1}\rho|u_z|^2+z^{n-1}\rho u^2\big)\,\dd z\dd t \Big)\nonumber\\[1.5mm]
	&\quad \leq C(d,M,E_0,T),
\end{align*}
we have
\begin{align}\label{2.72}
	\Big|\int_0^T\int_d^D I_7\, \dd r \dd t\Big| \leq C(d,M,E_0,T).
\end{align}

4. For $I_2$, integrating by parts, we have
\begin{align}
	&\Big|\int^r_d \frac{n-1}{z} u\rho_z w \,\dd z\Big|\nonumber\\
	&\leq \frac{n-1}{r} |(\rho u w)(t,r)|
	+C\Big|\int_d^r\frac{1}{z}\Big( \rho u_z w +\rho u w_z-\frac{1}{z}\rho uw\Big)(t,z)\dd z\Big|\nonumber\\
	&\leq \frac{n-1}{r} |(\rho u w)(t,r)|+C(d)\int_d^D \rho u_r^2\,r^{n-1}\dd r
	+C(d,M, E_0,T),\nonumber
\end{align}
which yields that
\begin{align}\label{2.73}
	&\Big|\int_0^T\int_d^DI_2 \,\dd r\dd t\Big|\\
    &=\Big|\int_0^T\int_d^D\v\r w\int_d^r \frac{n-1}{z} u\rho_z w\, \dd z \dd r\dd t\Big|\nonumber\\
	&\leq C(d,M,E_0,T)\Big(1+\int_0^T\int_d^D\v  \rho u_r^2\,r^{n-1} \dd r \dd t\Big)
	+\v\int_0^T\int_d^D \rho^3 w^2\, \dd r\dd t\nonumber\\
	&\leq C(d,M,E_0,T)+\v\int_0^T\int_d^D \rho^3 w^2\, \dd r\dd t.\nonumber
\end{align}

For $I_1$, notice that
\begin{align}\label{2.74}
	&\Big|\int_0^T\int_d^D I_1 \,\dd r \dd t\Big|\\
	&\leq \v\int_0^T\int_d^D \r^2\big|u_r+\frac{n-1}{r} u\big| w^2\,\dd r \dd t\nonumber\\
	&\leq \v\int_0^T\int_d^D \rho \Big(\r^2 +C\big(u_r+\frac{n-1}{r}u\big)^2\Big) w^2\, \dd r \dd t\nonumber\\
	&\leq C(d,M,E_0,T)+\v\int_0^T\int_d^D \r^3 w^2\, \dd r \dd t.\nonumber
\end{align}
To close the estimates for $I_1$ and $I_2$, we need to bound the last term of both \eqref{2.73}-RHS and \eqref{2.74}-RHS.
There are three cases:

\smallskip
{\it Case} 1:  $\gamma\in(1,2]$. Notice that
\begin{align}\label{2.75}
	&\v\int_0^T\int_d^D \r^3 w^2\, \dd r \dd t\\
	&\leq \v \int_0^T \big(\int_d^D \r^{\g} \dd r\big) \sup_{r\in[d,D]}\big(\r^{3-\g} w^2\big)\,\dd t\nonumber\\
	& \leq C(d,M,E_0)\int_0^T\v \sup_{r\in[d,D]}\big(\r^{3-\g} w^2\big)\, \dd t\nonumber\\
	&\leq C(d,M,E_0)\int_0^T\int_d^D\v\big|\big(\r^{3-\g} w^2\big)_r(t,r)\big|\,\dd r\dd t\nonumber\\
	&\leq C_2(d,M,E_0)\int_0^T\int_d^D\v \big( \r^{2-\g}|\r_r| w^2+ \r^{3-\g} w |w_r|\big)\, \dd r\dd t.\nonumber
\end{align}
A direct calculation shows that
\begin{align}
	&\int_0^T\int_d^D \v \r^{2-\g}|\r_r| w^2\, \dd r \dd t\label{2.76}\\
	&\leq \int_0^T\int_d^D \v\r^{\g-2} \r_r^2\,\dd r \dd t+\f\v2\int_{0}^T\int_d^D \r^{3(2-\g)} w^2\, \dd r \dd t\nonumber\\
	&\leq C(d,M,E_0,T)+\f{1}{4C_2(d,M,E_0)}\int_{0}^T\int_d^D \v \rho^{3} w^2\, \dd r \dd t,\nonumber
\end{align}
and
\begin{align}
	&\int_0^T\int_d^D \v \r^{3-\g} w |w_r|\, \dd r \dd t\label{2.77}\\
	& \leq \int_0^T\Big(\v \sup_{r\in[d,D]}(\r w)(t,r) \int_d^D\r^{2-\g}|w _r|\,\dd r\Big) \dd t\nonumber\\
	&\leq C(d,M,E_0)\int_0^T\v \sup_{r\in[d,D]}(\r w)(t,r)\,  \dd t\nonumber\\
	&\leq C(d,M,E_0)\int_0^T \int_d^D \v \big(|\r_r| w+ \r |w_r|\big) \,\dd r \dd t\nonumber\\
	&\leq C(d,M,E_0)\Big(\int_0^T\int_d^D \v\r^{\g-2} \r_r^2\,\dd r \dd t
	+\int_0^T\int_d^D \big(\rho+\rho^{2-\g} w\big)\, \dd r \dd t\Big)\nonumber\\
	&\leq C(d,M,E_0, T).\nonumber
\end{align}
Combining \eqref{2.75}--\eqref{2.77}, we have
\begin{align}\label{2.78}
	\v\int_0^T\int_d^D \r^3 w^2\, \dd r \dd t\leq C(d,M,E_0,T)\qquad\mbox{for $\g\in(1,2]$}.
\end{align}

\smallskip
{\it Case} 2: $\g\in[2,3]$. We have
\begin{align}\label{2.79}
&\v \int_0^T\int_d^D \r^3 w^2\, \dd r \dd t\\
&\leq \v \int_0^T\Big(\sup_{r\in[d,D]}(\r^2 w)(t,r)\int_d^D \r w \,\dd r\Big) \dd t\nonumber\\
&\leq C(d,M,E_0) \int_0^T \int_d^D \v\big(\r|\r_r| w+\r^2 |w_r|\big)\, \dd r \dd t\nonumber\\
&\leq C(d,M,E_0) \int_0^T \int_d^D \big(\v^2\r^{\g-2}|\r_r|^2 w+\r^2 |w_r|+\r^{4-\g} w\big)\,\dd r \dd t\nonumber\\
&\leq C(d,M,E_0,T).\nonumber
\end{align}

\smallskip
{\it Case} 3:  $\gamma\in(3,\infty)$. It is direct to see that
\begin{align}\label{2.80}
	\v \int_0^T\int_d^D \r^3 w^2 \,\dd r \dd t\leq C\int_0^T\int_d^D (\rho+\rho^\gamma)\,\dd r \dd t\leq C(d,M,E_0,T).
\end{align}
Now substituting \eqref{2.78}--\eqref{2.80} into \eqref{2.73}--\eqref{2.74} yields that
\begin{align}\label{2.81}
	\Big|\int_0^T\int_d^D (I_1+I_2)\, \dd r \dd t\Big|\leq C(d,M,E_0,T).
\end{align}
5. Integrating \eqref{2.69} over $[0,T]\times [d,D]$ and
then using \eqref{2.71}--\eqref{2.72} and \eqref{2.81},
we conclude \eqref{2.65}.
\end{proof}

To use the compensated compactness framework in \cite{Perepelitsa},
we still need to obtain the higher integrability on the velocity.
For this, we require to exploit several important properties of some special entropy pairs.

First, taking $\psi(s)=\frac12 {s}{|s|}$ in \eqref{weakentropy}, then the corresponding entropy and entropy flux are represented as
\begin{align}\label{2.82}
	\begin{cases}
		\displaystyle\eta^{\#}(\r,\r u)=\f12 \r \int_{-1}^1 (u+\r^{\t} s) |u+\r^{\t}s| [1-s^2]_+^{\fb} \dd s,\\[3mm]
		\displaystyle q^{\#}(\r, \r u)=\f12 \r \int_{-1}^1 (u+\theta\r^{\t}s)(u+\r^{\t} s) |u+\r^{\t}s| [1-s^2]_+^{\fb} \dd s.
	\end{cases}
\end{align}
A direct calculation shows that
\begin{align}\label{2.83}
	|\eta^{\#}(\r,\r u)|\leq C_{\g} \big(\r |u|^2+\r^{\g}\big), \qquad\,\, q^{\#}(\r,\r u)\geq C_{\g}^{-1} \big(\r |u|^3+\r^{\g+\t}\big),
\end{align}
where and whereafter $C_\g>0$ is a universal constant depending only on $\g>1$.
We regard $\eta^{\#}$ as a function of $(\rho,m)$ to obtain
\begin{align}\nonumber
	\begin{cases}
		\displaystyle\eta^{\#}_\r=\int_{-1}^1 \big(-\f12 u+(\t+\f12)\r^{\t} s\big)\,|u+\r^{\t}s| [1-s^2]_+^{\fb} \dd s,\\[3mm]
		\displaystyle \eta^{\#}_m=\int_{-1}^1|u+\r^{\t}s| [1-s^2]_+^{\fb} \dd s.
	\end{cases}
\end{align}
It is direct to check that
\begin{align}\label{2.85}
	|\eta^{\#}_m|\leq C_\g \big(|u|+\rho^\theta\big),\qquad\,\, |\eta^{\#}_\rho|\leq C_\g\big(|u|^2+\rho^{2\theta}\big).
\end{align}
From \cite{Perepelitsa,Chen6}, we know that
\begin{align}\label{2.86}
	\r u \partial_\r \eta^{\#}+\r u^2 \partial_m\eta^{\#}-q^{\#}
	&=\f{\t}{2} \rho^{1+\theta}\int_{-1}^1  (u-\r^{\t}s)s |u+\r^{\t}s|\, [1-s^2]_+^{\fb} \dd s\leq 0.
\end{align}

The following lemma is crucial to control the trace estimates for the higher integrability on the velocity.
In fact, we have the boundary parts $(u \eta^{\#})(t,b(t))$ and $q^{\#}(t,b(t))$,
and it is impossible to have the uniform trace bound (independent of $\v$) for each of them.
Our key point is to identify the cancelation between these two boundary parts.

\begin{lemma}\label{lem2.5}
	For the entropy pair defined in \eqref{2.82},
	\begin{align}\label{2.87}
		|q^{\#}-u\eta^{\#}|\leq C_\g \big(\rho^{\gamma} |u|+ \rho^{\gamma+\theta} \big).
	\end{align}
\end{lemma}

\begin{proof} It follows from \eqref{2.82} that
\begin{align}\label{2.88}
	q^{\#}-u\eta^{\#}&=\frac{1}{2}\theta\rho^{1+2\theta}\int_{-1}^1s^2|u+\rho^\theta |[1-s^2]_+^{\fb} \dd s\\
	&\quad+\frac{1}{2}\theta\rho^{1+\theta}u\int_{-1}^1s|u+\rho^\theta s|[1-s^2]_+^{\fb} \dd s\nonumber\\
	&:=I_1+I_2.\nonumber
\end{align}
A direct calculation shows that
\begin{equation}\label{2.89}
	|I_1|\leq C_{\g}\big(\rho^\gamma|u|+\rho^{\gamma+\theta}\big).
\end{equation}

For $I_2$, we note that $I_2=0$ if $u=0$. Thus, it suffices to consider $u\neq0$. We divide the proof into three cases.

{\it Case} 1.  If $u>0$ and $u+\rho^\theta s\geq0$ for all $s\in[-1,1]$, then it follows that
\begin{align}
	\int_{-1}^1s|u+\rho^\theta s|[1-s^2]_{+}^{\fb}\,\dd s
	&=u \int_{-1}^1 s[1-s^2]_{+}^{\fb}\,\dd s+\r^\theta\int_{-1}^{1}s^2[1-s^2]_{+}^{\fb}\,\dd s\nonumber\\
	&=\r^\theta\int_{-1}^{1}s^2[1-s^2]_{+}^{\fb}\,\dd s,\nonumber
\end{align}
which
yields that
\begin{equation}\label{2.90}
	|I_2|\leq C_{\g}\rho^{1+2\theta}|u|=C_{\g}\rho^{\g}|u|.
\end{equation}

{\it Case} 2.
If $u>0$ and $u+\rho^\theta s_0=0$ for some $s_0\in[-1,1]$, then
$s_0=-\frac{u}{\rho^\theta}\in[-1,1]$ so that
\begin{equation}\nonumber
	|u|\leq \rho^{\theta},
\end{equation}
which yields that
\begin{align}\label{2.91}
	|I_2|\leq C_\g \rho^{1+\theta}|u|\big(|u|+\rho^\theta\big)
	\leq C_\g \rho^{1+3\theta}=C_\g \rho^{\gamma+\theta}.
\end{align}

{\it Case} 3. If $u<0$, by similar arguments as in Cases 1--2,
\begin{equation*}
	|I_2|\leq C_\g \big(\rho^{\g}|u|+\rho^{\gamma+\theta}\big),
\end{equation*}
which, together with \eqref{2.88}--\eqref{2.91} yields \eqref{2.87}.
\end{proof}

\begin{lemma}[Higher integrability on the velocity]\label{lem2.6}
	Let $(\r,u)$ be the smooth solution of \eqref{2.1}--\eqref{2.6}.
	Then, under the assumption of Lemma {\rm \ref{lem2.1}},
	\begin{align}\label{2.95}
		\int_0^T\int_d^D  \big(\rho|u|^3+\rho^{\gamma+\theta}\big)(t,r)\, r^{n-1} \dd r \dd t\leq C(d,D,M,E_0, T)
	\end{align}
	for any $(d,D)\Subset [a, b(t)]$.
\end{lemma}

\begin{proof}  Multiplying $\eqref{2.1}_1$ by $r^{n-1} \eta^{\#}_\r$ and $\eqref{2.1}_2$ by $r^{n-1} \eta^{\#}_m$,
we have
\begin{align}\label{2.96}
& ( \eta^{\#}r^{n-1})_t+( q^{\#}r^{n-1})_r+(n-1) \big(-q^{\#}+\r u \eta^{\#}_\r+\r u^2 \eta^{\#}_m\big)r^{n-2}\\
& =\v  \eta^{\#}_m\, \Big((\rho u_r)_r+(n-1)\rho \Big(\frac{u}{r}\Big)_r\Big)r^{n-1}
-\kappa \eta^{\#}_m\, \rho \int_a^r \rho\,z^{n-1} \dd z.\nonumber
\end{align}
Using \eqref{2.6}, a direct calculation shows that
\begin{align}\label{2.97}
\frac{\dd}{\dd t}\int_r^{b(t)}\eta^{\#}\,z^{n-1}\dd z
& =\eta^{\#}(t,b(t))b(t)^{n-1} b'(t)+\int_r^{b(t)}\eta^{\#}_t(t,z)\,z^{n-1}\dd z\\
&= (u\eta^{\#})(t, b(t))b(t)^{n-1} +\int_r^{b(t)} \eta^{\#}_t(t,z)\,z^{n-1}\dd z.\nonumber
\end{align}
Integrating \eqref{2.96} over $[r,b(t))$, then using \eqref{2.86} and  \eqref{2.97},  we have
\begin{align}\label{2.98}
	&q^{\#}(t,r)r^{n-1} \\
   &\leq -\v \int_r^{b(t)}  \eta^{\#}_m(t,z)(\rho u_z)_z\, z^{n-1}\dd z
	-(n-1)\v\int_r^{b(t)}\eta^{\#}_m(t,z)\rho\big(\frac{u}{z}\big)_z\, z^{n-1}\dd z\nonumber\\
	&\quad\, +\Big( \int_r^{b(t)}\eta^{\#}(t,z)\,z^{n-1}\dd z\Big)_t
	+\big(q^{\#}-u\eta^{\#}\big)(t,b(t))b(t)^{n-1}\nonumber\\
	&\quad\,+\kappa \int_r^{b(t)} \big(\int_a^y \rho \,z^{n-1} \dd z\big)\rho\,\eta^{\#}_m\, \dd y.\nonumber
\end{align}

We now bound each term of \eqref{2.98}-RHS. First,
for the term involving the trace estimates in \eqref{2.98}, it follows from \eqref{2.44} and  Lemmas \ref{lem2.1}, \ref{lem2.2}, and \ref{lem2.5}
that
\begin{align}\label{2.99}
&\int_0^T \big|(q^{\#}-u\eta^{\#})(t,b(t))\big|b(t)^{n-1}\, \dd t\\
	&\leq C\int_0^T\big(\rho^{\gamma+\theta}(t,b(t))+(\rho^{\gamma}|u|)(t, b(t))\big)b(t)^{n-1}\,\dd t\nonumber\\
	&\leq  C\Big(\int_0^T\v  (\r|u|^2)(t,b(t))b(t)^{n-2}\,\dd t\Big)^{\frac{1}{2}}
	\Big(\int_0^T\frac{1}{\v}\rho^{2\gamma-1}(t,b(t)) b(t)^n\,\dd t\Big)^{\frac{1}{2}}\nonumber\\
	&\,\quad+C(M,E_0,T)\int_0^T\rho^{{\frac{\gamma}{n}+\theta}}(t,b(t))\,\dd t\leq C(M,E_0,T).\nonumber
\end{align}
Observe from \eqref{2.99} that the free boundary approximation is
ideal for the existence of solutions of CEPEs with  finite mass.

To estimate the first term of \eqref{2.98}-RHS, we integrate by parts to obtain
\begin{align}\label{2.100}
	&\v\int_r^{b(t)}\eta^{\#}_m(\rho u_z)_z\,z^{n-1}\dd z\\
	&=\v\Big(\eta^{\#}_m(t,b(t))\, (\rho u_r)(t,b(t))b(t)^{n-1}-\eta^{\#}_m(t,r)(\rho u_r)(t,r)r^{n-1} \Big)\nonumber\\
	&\quad\, -(n-1)\v\int_r^{b(t)}\eta^{\#}_m\, \rho u_z\,z^{n-2}\dd z
	-\v\int_r^{b(t)}\rho u_z(\eta^{\#}_{mu}u_z+\eta^{\#}_{m\rho}\rho_z)\,z^{n-1}\dd z,\nonumber
\end{align}
where we have regarded $\eta^{\#}_m$ as a function of $(\rho, u)$. Using  \eqref{2.5} and $\eqref{2.85}$, we have
\begin{align}\label{2.101}
&\big|\v \eta^{\#}_m(t,b(t))\, (\rho u_r)(t,b(t))b(t)^{n-1}\big|\\
	&=\Big|\eta^{\#}_m(t,b(t))\,\Big(\v\rho\big(u_r+\frac{n-1}{r}u\big)(t,b(t))-(n-1)\v b(t)^{-1}\,(\rho u)(t,b(t))\Big)b(t)^{n-1}\Big|\nonumber\\
	&=\left|a_0 \eta^{\#}_m(t,b(t))\, \rho^{\g}(t,b(t))b(t)^{n-1}-(n-1)\v \eta^{\#}_m(t,b(t))\,(\rho u)(t,b(t))b(t)^{n-2}\right|\nonumber\\[1mm]
	&\leq C\big\{(\rho^{\g}|u|)(t,b(t))+\rho^{\gamma+\theta}(t,b(t))\big\}b(t)^{n-1}\nonumber\\[0.5mm]
	&\quad\, +C\v \big\{(\rho|u|^2)(t,b(t))+(\rho^{1+\theta}|u|)(t,b(t))\big\}b(t)^{n-2}  \nonumber\\[1mm]
	&\leq C\big\{(\rho^{\g}|u|)(t,b(t))+\rho^{\gamma+\theta}(t,b(t))\big\}b(t)^{n-1} \nonumber\\[0.5mm]
	&\quad\, +C\v\big\{(\rho|u|^2)(t,b(t))+ \rho^{\gamma}(t,b(t))\big\} b(t)^{n-2}.\nonumber
\end{align}
Thus, using \eqref{2.101}, by similar arguments as in \eqref{2.99},  we have
\begin{align}\label{2.102}
	\int_0^T
	\big|\v b^{n-1}(t)\eta^{\#}_m(t,b(t))(\rho u_r)(t,b(t))\big|\,
	\dd t\leq C(M,E_0, T).
\end{align}

A direct calculation shows that
\begin{align}\label{2.103}
	|\eta^{\#}_{mu}| + \r^{1-\t} |\eta^{\#}_{m\r}|\leq C.
\end{align}
Integrating \eqref{2.100} over $[0,T]\times [d,D]$ and then using \eqref{2.85} and \eqref{2.102}--\eqref{2.103} lead to
\begin{align}\label{2.104}
		&\int_0^T\int_d^D\Big|\v\int_r^{b(t)}  \eta^{\#}_m\,  (\rho u_z)_z \,z^{n-1}\dd z\Big| \dd r \dd t\\
	&\leq C(D, M,E_0,T)+C\v\int_0^T\int_d^D  \rho |u_r| \big(|u|+\rho^\theta\big)\,r^{n-1}\dd r\dd t\nonumber\\
	&\quad\, +C\int_0^T\int_d^D\int_r^{b(t)}\v \rho|u_z|\big(|u|+\rho^\theta\big)\,z^{n-2}\dd z\dd r\dd t\nonumber\\
	&\quad\, +C\int_0^T\int_d^D\int_r^{b(t)}\v \rho|u_z|\big(|u_z|+\rho^{\theta-1}|\rho_z|\big)\, z^{n-1}\dd z\dd r\dd t\nonumber\\
	&\leq C(D, M,E_0,T)+C(d,D)\int_0^T\int_d^{b(t)}\v\big( \rho|u_z|^2+\rho^{\gamma-2}|\rho_z|^2\big)\,z^{n-1}\dd z \dd t\nonumber\\
	&\quad\, +C(d,D)\int_0^T\int_d^{b(t)}\v \big(\rho |u|^2+ \rho^{\g}\big)\,z^{n-1} \dd z \dd t \nonumber\\
	&\leq  C(d,D, M,E_0,T),\nonumber
\end{align}
where we have used Lemmas \ref{lem2.1}--\ref{2.2}.

For the second term of \eqref{2.98}-RHS, we have
\begin{align}\label{2.105}
		&\int_0^T\int_d^D\Big|\int_r^{b(t)}\v \eta^{\#}_m(t,z)\, \rho(\frac{u}{z})_z\, z^{n-1}\dd z\Big|\dd r\dd t\\
	&\leq C(D)\int_0^T\int_d^{b(t)}\v \big(\rho |u|+\rho^{1+\theta}\big)\big(\frac{|u_z|}{z}+\frac{|u|}{z^2}\big)\,z^{n-1}\dd z\dd t,\nonumber\\
	&\leq C(d,D)\int_0^T\int_d^{b(t)}\v \big(\rho |u|^2+\rho |u_z|^2+\rho^{\g}\big)\,z^{n-1}\dd z\dd t\nonumber\\
	&\leq C(d,D, M,E_0,T).\nonumber
\end{align}
For the third term of \eqref{2.98}-RHS,
\begin{align}\label{2.106}
	&\Big|\int_0^T\int_d^D\Big(\int_r^{b(t)}\eta^{\#}(\rho,\rho u)\, z^{n-1}\dd z\Big)_t\dd r\dd t\Big|\\
	&\leq \Big|\int_d^D\int_r^{b(T)}\eta^{\#}(\rho,\rho u)(T,z)\, z^{n-1}\dd z\dd r\Big|\nonumber\\
	&\quad\, +\Big|\int_d^D\int_r^{b}\eta^{\#}(\rho_0,\rho_0 u_0)\,z^{n-1} \dd z\dd r\Big|\nonumber\\[2mm]
	&\leq C(D,M,E_0).\nonumber
\end{align}
For the last term of \eqref{2.98}-RHS, it follows from \eqref{2.83} that
\begin{align}\label{2.107}
	&\int_0^T\int_d^D\Big|\kappa \int_r^{b(t)}\eta^{\#}_m\, \rho \big(\int_a^y \rho\, z^{n-1} \dd z\big) dy\Big|\dd r\dd t\\
	&\leq C(D,M)\int_0^T\int_d^{b(t)} \big(\rho|u|+\rho^{1+\theta}\big)\, \dd r \dd t\nonumber\\
	&\leq C(D,M)\int_0^T\int_d^{b(t)} \big(\rho|u|^2+\rho+\rho^{\gamma}\big)\, \dd r \dd t
	\leq C(d,D,M,E_0).\nonumber
\end{align}

Integrating \eqref{2.98} over $[0,T]\times [d,D]$, then using \eqref{2.83}, \eqref{2.99}, and \eqref{2.104}--\eqref{2.107},
we conclude that
\begin{align}
	\int_0^T \int_d^D \big(\r |u|^3+\r^{\g+\theta}\big)\, r^{n-1}\dd r \dd t
	\leq C \int_0^T \int_d^D  q^{\#}(t,r)\, r^{n-1} \dd r \dd t
	\leq C(d,D, M,E_0,T).\nonumber
\end{align}
\end{proof}

\section{Existence of Global Finite-Energy Solutions}

In this section, for fixed $\v>0$, we take limit $b\rightarrow\infty$ to obtain
global weak solutions of CNSPEs with some uniform bounds,
which are essential for applying the compensated compactness framework in \S 5 below.
We often denote the solutions of \eqref{2.1}--\eqref{initial}
as $(\rho^{\v,b}, u^{\v,b})$ for simplicity of presentation
in this section, since $\rho^{\v,b}>0$ on $[0,T]\times [a,b(t)]$ for fixed $b>0$.

To take the limit, we have to be careful, since the weak solutions  may involve the vacuum.
We use the similar compactness arguments as in  \cite{Mellet,Zhenhua Guo2} to handle
the limit: $b\rightarrow\infty$.
First of all, we understand our solutions $(\rho^{\v,b}, u^{\v,b})$ to be the zero extension
of $(\rho^{\v,b}, u^{\v,b})$ on  $\big([0,T]\times[0,\infty)\big)\backslash \Omega_T$.
It follows from  Lemma \ref{lem2.3} that
\begin{equation}\label{3.3-1}
	\lim_{b\rightarrow\infty} \sup_{t\in[0,T]} b(t)= \infty,
\end{equation}
which implies that domain $[0,T]\times [a,b(t)]$ expands to $[0,T]\times (0,\infty)$
as $b\to \infty$.
That is, for any compact set $K\Subset (0,\infty)$, when $b\gg1$,
$K\Subset (a,b(t))$ for all $t\in[0,T]$.

Now we define
\begin{equation}\label{5.5-1}
	\M^{\v,b}(t,\textbf{x}):=m^{\v,b}(t,r) \frac{\textbf{x}}{r}=(\rho^{\v,b}u^{\v,b})(t,r)  \frac{\textbf{x}}{r}.
\end{equation}
Then it is direct to check that the corresponding vector function $(\rho^{\v,b}, \M^{\v,b}, \Phi^{\v,b})$
is a classical solution of CNSPEs
for $(t,\textbf{x})\in [0,\infty)\times \Omega_t $:
\begin{align}\label{5.5-2}
	\begin{cases}
		\partial_t \rho^{\v,b}+\mbox{div}\,\M^{\v,b}=0,\\[2mm]
		\partial_t \M^{\v,b}+\mbox{div}\big(\frac{\M^{\v,b}\otimes \M^{\v,b}}{\rho^{\v,b}}\big)+\nabla p(\rho^{\v,b})
		= -\rho^{\v,b} \nabla\Phi^{\v,b} +\v \mbox{div}\big(\rho^{\v,b} D(\frac{ \M^{\v,b}}{\rho^{\v,b}})\big),\\[2mm]
		\Delta \Phi^{\v,b} =\kappa \rho^{\v,b},
	\end{cases}
\end{align}
with $\M^{\v,b}\Big|_{\textbf{x}\in \partial B_a(\mathbf{0})}={\bf 0}$.

\subsection{Taking limit $b\rightarrow \infty$}

\begin{lemma}\label{lem5.2}
	For fixed $\v>0$, there exists a function $\rho^{\v}(t,r)$ such that, as $b\rightarrow\infty$ $($up to a subsequence$)$,
	\begin{equation} \label{5.12}
		(\sqrt{\rho^{\v,b}}, \,\rho^{\v,b})\longrightarrow (\sqrt{\rho^{\v}},\, \rho^{\v})  \quad  \mbox{ {\it a.e.} and strongly in} \  C(0,T; L^q_{\rm loc})
	\end{equation}
	for any $q\in [1,\infty)$, where $L^q_{\rm loc}$ denotes  $L^q(K)$ for $K\Subset (0,\infty)$.
\end{lemma}

\begin{proof}
 It follows from Lemmas \ref{lem2.1} and \ref{lem2.2}  that
\begin{equation}\label{3.6-1}
	\sqrt{\rho^{\v,b}}\in L^\infty(0,T; H^{1}_{\rm loc})\hookrightarrow L^\infty(0,T; L^{\infty}_{\rm loc})
	\qquad\mbox{uniformly in $b>0$}.
\end{equation}
Using the mass equation $\eqref{2.1}_1$ and Lemma \ref{lem2.1}, we have
\begin{align}
	-\partial_t\sqrt{\rho^{\v,b}}&= \Big(\sqrt{\rho^{\v,b}}\Big)_r u^{\v,b}+ \frac12 \sqrt{\rho^{\v,b}} u_r^{\v,b}+\frac{n-1}{2r}\sqrt{\rho^{\v,b}} u^{\v,b}\nonumber\\
	&=\Big(\sqrt{\rho^{\v,b}} u^{\v,b}\Big)_r-\frac12 \sqrt{\rho^{\v,b}} u_r^{\v,b}
	+\frac{n-1}{2r}\sqrt{\rho^{\v,b}} u^{\v,b} \in L^{2}(0,T; H^{-1}_{\rm loc})\nonumber
\end{align}
uniformly in $b>0$,
which, together with the Aubin-Lions lemma, yields that
\begin{align}
	\sqrt{\rho^{\v,b}} \quad \mbox{is compact in $C(0,T; L^q_{\rm loc})$ for any $q\in [1,\infty)$}.\nonumber
\end{align}

Notice that, for any $K\Subset (0,\infty)$ and $b_1,b_2\in (1,\infty)$,
$$
|\rho^{\v,b_1}(t,r)-\rho^{\v,b_2}(t,r)|\le C_{T,K} \Big|\sqrt{\rho^{\v,b_1}}-\sqrt{\rho^{\v,b_2}}\Big|
\qquad \mbox{for any $(t,r)\in [0, T]\times K$},
$$
where $C_{T,K}>0$ is a constant independent of $b_1$ and $b_2$.
Then there exists a function $\rho^{\v}(t,r)$ such that, as $b\rightarrow\infty$ $($up to a subsequence$)$,
$(\sqrt{\rho^{\v,b}}, \rho^{\v,b})\longrightarrow (\sqrt{\rho^{\v}}, \rho^{\v})$ a.e. and strongly in $C(0,T; L^q_{\rm loc})$ for any $q\in [1,\infty)$.
Then \eqref{5.12} follows.
\end{proof}

\begin{corollary}\label{cor5.3}
	For fixed $\v>0$, the pressure function sequence $p(\rho^{\v,b})$ is uniformly bounded in $L^\infty(0,T; L^q_{\rm loc})$ for all $q\in[1,\infty]$
	and, as $b\rightarrow\infty$ $($up to a subsequence$)$,
	\begin{align*}
		p(\rho^{\v,b})\longrightarrow p(\rho^{\v}) \qquad  \mbox{strongly in $L^q(0,T; L^q_{\rm loc})$ for all $q\in[1,\infty)$}.
	\end{align*}
\end{corollary}

\begin{lemma}\label{lem5.5}
	For fixed $\v>0$, as $b\rightarrow\infty$ $($up to a subsequence$)$,
	the momentum function sequence $m^{\v,b}:=\rho^{\v,b} u^{\v,b}$ converges strongly in $L^2(0,T; L^q_{\rm loc})$ to
	some function $m^{\v}(t,r)$ for all $q\in[1,\infty)$.
	In particular, we have
	\begin{align*}
		m^{\v,b}=\rho^{\v,b} u^{\v,b}\longrightarrow m^{\v}(t,r) \qquad  \mbox{a.e. in $[0,T]\times(0,\infty)$}.
	\end{align*}
\end{lemma}

\begin{proof}
Notice that $\sqrt{\rho^{\v,b}}$ is uniformly bounded in  $L^\infty(0,T; L^{\infty}_{\rm loc})$
and $\sqrt{\rho^{\v,b}} u^{\v,b}$ is uniformly bounded in $L^\infty(0,T; L^{2}_{\rm loc})$ in $b>0$,
which imply that
\begin{align}\label{5.18}
	\rho^{\v,b} u^{\v,b}=\sqrt{\rho^{\v,b}} \,\Big(\sqrt{\rho^{\v,b}} u^{\v,b}\Big)
	\quad \mbox{is uniformly bounded in} \  L^\infty(0,T; L^{2}_{\rm loc}).
\end{align}
A direct calculation shows that
\begin{align}\label{5.19}
	(\rho^{\v,b} u^{\v,b})_r&=\rho^{\v,b}_r u^{\v,b}+\rho^{\v,b} u_r^{\v,b}\\
	&=2\Big(\sqrt{\rho^{\v,b}}\Big)_r\, \Big(\sqrt{\rho^{\v,b}}u^{\v,b}\Big)
      +\sqrt{\rho^{\v,b}}\,\Big(\sqrt{\rho^{\v,b}} u_r^{\v,b}\Big)\nonumber
\end{align}
is uniformly bounded in $L^{2}(0,T; L^1_{\rm loc})$.
Then it follows from \eqref{5.18}--\eqref{5.19} that
\begin{align}\label{5.20}
	\rho^{\v,b} u^{\v,b} \qquad  \mbox{is uniformly bounded in $L^{2}(0,T; W^{1,1}_{\rm loc})$}.
\end{align}

It follows from \eqref{3.6-1} and Lemma \ref{lem2.1} that
\begin{align}\label{5.21}
\begin{cases}
\displaystyle\partial_r\Big(\Big(\sqrt{\rho^{\v,b}} u^{\v,b}\Big)^2\Big)\in L^{\infty}(0,T;W_{\rm loc}^{-1,1}),\\[2.5mm]
\displaystyle\frac{n-1}{r} \Big(\sqrt{\rho^{\v,b}} u^{\v,b}\Big)^2 \in L^{\infty}(0,T; L_{\rm loc}^{1}),\\[2.5mm]
\displaystyle \partial_r p(\rho^{\v,b}) \in L^{2}(0,T;H_{\rm loc}^{-1}),\\[2mm]
\displaystyle	\kappa \frac{\rho^{\v,b}}{r^{n-1}} \int_a^r \rho^{\v,b}(t,z)\,z^{n-1} \dd z \in L^\infty(0,T; L^\infty_{\rm loc}),
\end{cases}
\end{align}
and
\begin{align}\label{5.22}
	\sqrt{\rho^{\v,b}}\,\Big(\sqrt{\rho^{\v,b}} u^{\v,b}_r\Big)
	+\frac{n-1}{r} \sqrt{\rho^{\v,b}}\,\Big(\sqrt{\rho^{\v,b}}u^{\v,b}\Big) \in L^2(0,T; L^2_{\rm loc})
\end{align}
uniformly in $b>0$.

Therefore, it follows from \eqref{5.22} that
\begin{align}\label{5.23}
	\partial_r\big(\rho^{\v,b}(u_r^{\v,b}+\frac{n-1}{r} u^{\v,b}) \big)\in L^{2}(0,T; H^{-1}_{\rm loc})
	\qquad \mbox{uniformly in $b$}.
\end{align}
Also, using Lemmas \ref{lem2.1} and \ref{lem2.2}, we have
\begin{align}\label{5.24}
	&\frac{n-1}{r} u^{\v,b} \partial_r \rho^{\v,b}\\
	&=\frac{2(n-1)}{r} \Big(\sqrt{\rho^{\v,b}}\Big)_r \,\Big(\sqrt{\rho^{\v,b}}u^{\v,b}\Big)\in L^2(0,T; L^1_{\rm loc})
	\,\,\, \mbox{uniformly in $b$}.\nonumber
\end{align}
Thus, substituting \eqref{5.21} and \eqref{5.23}--\eqref{5.24} into $\eqref{2.1}_2$ yields that
\begin{align}
	\partial_t(\rho^{\v,b} u^{\v,b})\in L^{2}(0,T; W^{-1,1}_{\rm loc})
	\qquad \mbox{uniformly in $b>0$}, \nonumber
\end{align}
which, together with \eqref{5.20} and the Aubin-Lions lemma,
implies that
\begin{align*}
	\rho^{\v,b} u^{\v,b} \qquad  \mbox{is compact in $L^{2}(0,T; L^q_{\rm loc})$ for all $q\in[1,\infty)$}.
\end{align*}
This completes the proof.
\end{proof}

\begin{lemma}\label{lem5.6}
	The limit function $m^\v(t,r)$ in Lemma {\rm \ref{lem5.5}} satisfies that
	$m^{\v}(t,r)=0$ a.e. on $\{(t,r)\,:\,\rho^{\v}(t,r)=0\}$. Furthermore, there exists
	a function $u^{\v}(t,r)$ such that $m^{\v}(t,r)=\rho^{\v}(t,r) u^{\v}(t,r)$ a.e., and $\,u^\v(t,r)=0$ a.e.
	on $\{(t,r)\,:\,\rho^{\v}(t,r)=0\}$. Moreover, as $b\rightarrow\infty$ $($up to a subsequence$)$,
	\begin{align*}
		& m^{\v,b} \longrightarrow m^{\v}=\rho^{\v} u^{\v} \qquad  \mbox{strongly in $L^2(0,T; L^q_{\rm loc})$ for $q\in[1,\infty)$},\\
		& \frac{m^{\v,b}}{\sqrt{\rho^{\v,b}}}\longrightarrow \sqrt{\rho^{\v}} u^{\v}=:\frac{m^\v}{\sqrt{\rho^\v}}
		\qquad  \mbox{strongly in $L^2(0,T; L^2_{\rm loc})$}.
	\end{align*}
\end{lemma}

\begin{proof}
We divide the proof into three steps.

\smallskip
1. We first claim that $m^{\v}(t,r)=0\,\,$ {\it a.e.} on $\{(t,r) \,:\, \rho^{\v}(t,r)=0\}$.

To prove this claim,  for any given $T>0$ and $0<d<D<\infty$, we define
	$$
	V:=\big\{(t,r)\in[0,T]\times [d,D] \,:\, \mbox{$\rho^{\v}(t,r)=0$ and $m^{\v}(t,r)\neq 0$}\big\}\setminus \mathcal{N},
	$$
	where $\mathcal{N}$ is the set where the subsequence (still denoted) $(\rho^{\v,b}, m^{\v,b})$ does not converge to $(\rho^\v, m^\v)$
	so that the Lebesgue measure
	of $\mathcal{N}$ must be zero: $|\mathcal{N}|=0$, since $(\rho^{\v,b}, m^{\v,b})$ converges to $(\rho^\v, m^\v)$ {\it a.e.} as $b\to \infty$.

	If $|V|=0$, then we have done. If $|V|>0$, then it is clear that
	\begin{align}\label{4.14-1}
		\liminf_{b\to\infty} \frac{|m^{\v,b}(t,r)|^2}{\rho^{\v,b}(t,r)}
		=\lim_{b\to\infty} \frac{|m^{\v,b}(t,r)|^2}{\rho^{\v,b}(t,r)}
		=\infty \qquad\mbox{for}\,\, (t,r)\in V.
	\end{align}
	On the other hand, notice that  $\sqrt{\rho^{\v,b}}u^{\v,b}r^{\frac{n-1}{2}} =\frac{m^{\v,b}}{\sqrt{\rho^{\v,b}}}r^{\frac{n-1}{2}}$
	is uniformly bounded in $L^{\infty}(0,T;L^2)$. Then Fatou's lemma implies that
	\begin{align}\label{4.14-2}
     	&\int_0^T\int_d^D \liminf\limits_{b\to \infty} \frac{|m^{\v,b}(t,r)|^2}{\rho^{\v,b}(t,r)}\,r^{n-1}\dd r\dd t\\
     	&\leq \liminf\limits_{b\rightarrow\infty}\int_0^T\int_d^D \frac{|m^{\v,b}(t,r)|^2}{\rho^{\v,b}(t,r)}\,r^{n-1}\dd r\dd t
     	\leq C(T,E_0,M)<\infty.\nonumber
	\end{align}
	Combining \eqref{4.14-1} with \eqref{4.14-2} yields
	\begin{align}\label{4.14-3}
    \infty&= \int_0^T\int_d^D \liminf\limits_{b\to \infty}  \frac{|m^{\v,b}(t,r)|^2}{\rho^{\v,b}(t,r)}\,
    {\mathbf{1}}_{V}(t,r)\,r^{n-1}\dd r\dd t\\
    &\leq \int_0^T\int_d^D \liminf\limits_{b\to \infty} \frac{|m^{\v,b}(t,r)|^2}{\rho^{\v,b}(t,r)}\,r^{n-1}\dd r\dd t
      \leq C(T,E_0,M)<\infty, \nonumber	
	\end{align}
	which is impossible,
	where ${\mathbf{1}}_{V}(t,r)$ is the indicator function of set $V$.
	Therefore, it must be that $|V|=0$, which leads to the claim.

\smallskip
2. Now we can define  velocity $u^{\v}(t,r)$ as
\begin{align*}
	u^{\v}(t,r):=
	\begin{cases}
		\dis\frac{m^{\v}(t,r)}{\rho^{\v}(t,r)} \quad {\it a.e.} \,\,\, \mbox{on}\,\, \{(t,r) \,:\, \rho^{\v}(t,r)\neq 0\},\\[2mm]
		\dis 0,\qquad\quad\,\,\,\, {\it a.e.} \,\,\, \mbox{on}\,\, \{(t,r) \,:\, \rho^{\v}(t,r)= 0\},
	\end{cases}
\end{align*}
and define
\begin{align*}
	\dis \f{m^\v(t,r)}{\sqrt{\rho^{\v}(t,r)}}:=
	0 \qquad {\it a.e.} \,\, \mbox{on}\,\, \{(t,r) \,:\, \rho^{\v}(t,r)= 0\}.
\end{align*}
Then it is clear that
\begin{align}\label{4.14-4}
	m^\v(t,r)=\rho^{\v}(t,r) u^{\v}(t,r) \,\,\,\,\, a.e.,\quad\,\,\,
	\f{m^\v(t,r)}{\sqrt{\rho^{\v}(t,r)}}&=\sqrt{\rho^{\v}(t,r)} u^{\v}(t,r) \,\,\,\,\, a.e.
\end{align}
It follows from \eqref{4.14-2}, \eqref{4.14-4}, and Lemmas \ref{lem5.2} and \ref{lem5.5} that
\begin{align}\label{4.17}
	&\int_0^T\int_d^D\rho^\v |u^\v|^2\, r^{n-1}\dd r\dd t\\
	&=\int_0^T\int_d^D\rho^\v |u^\v|^2\,\mathbf{1}_{\{\rho^\v(t,r)>0\}} \, r^{n-1}\dd r\dd t\nonumber\\
	&=\int_0^T\int_d^D \frac{|m^\v|^2}{\rho^\v}\,\mathbf{1}_{\{\rho^\v(t,r)>0\}} \,r^{n-1}\dd r\dd t\nonumber\\
	&=\int_0^T\int_d^D \liminf\limits_{b\to \infty} \frac{|m^{\v,b}(t,r)|^2}{\rho^{\v,b}(t,r)}\,
       \mathbf{1}_{\{\rho^\v(t,r)>0\}} \,r^{n-1}\dd r\dd t\nonumber\\
	&\leq \int_0^T\int_d^D \liminf\limits_{b\to \infty} \frac{|m^{\v,b}(t,r)|^2}{\rho^{\v,b}(t,r)}\,r^{n-1}\dd r\dd t\nonumber\\[1mm]
	&\leq C(T,E_0,M)<\infty,\nonumber
\end{align}
where $\mathbf{1}_{\{\rho^\v(t,r)>0\}}$ is the indicator function of set $\{\rho^\v(t,r)>0\}$.
Similarly, it follows from Lemma \ref{lem2.6}  and Fatou's lemma that
\begin{align}\label{5.28}
	\int_0^T\int_{d}^{D} \rho^{\v} |u^{\v}|^3\,\dd r\dd t
	&\leq \lim_{b\rightarrow\infty}\int_0^T\int_{d}^{D} \rho^{\v,b} |u^{\v,b}|^3\,\dd r\dd t\\[2mm]
	&\leq C(d,D,M,E_0,T)<\infty. \nonumber
\end{align}

\smallskip
3. Next, since $(\rho^{\v,b}, m^{\v,b})$ converges {\it a.e.}, it is direct to know that
sequence $\sqrt{\rho^{\v,b}} u^{\v,b}=\frac{m^{\v,b}}{\sqrt{\rho^{\v,b}}}$  converges {\it a.e.}
to $\sqrt{\rho^{\v}} u^{\v}=\frac{m^{\v}}{\sqrt{\rho^{\v}}}$ on  $\{(t,r)\,:\,\rho^\v(t,r)\neq0\}$.
Moreover, for any given positive constant $k\geq1$, we have
\begin{align}\label{5.29}
	\sqrt{\rho^{\v,b}} u^{\v,b} \mathbf{1}_{\{|u^{\v,b}|\leq k\}}
	\longrightarrow \sqrt{\rho^{\v}} u^{\v} \mathbf{1}_{\{|u^{\v}|\leq k\}} \quad  \mbox{\it a.e.}.
\end{align}
It is direct to know that
\begin{align*}
	&\int_0^T\int_d^D \Big(\big|\sqrt{\rho^{\v,b}} u^{\v,b}\big|^{\f{12}{5}}
          + |\sqrt{\rho^{\v}} u^{\v}|^{\f{12}{5}}\Big)\,\dd r\dd t\nonumber\\
	&\leq C \int_0^T\int_{d}^{D}
	\big(\rho^{\v,b} |u^{\v,b}|^3 + (\rho^{\v,b})^{\gamma+1} + \rho^{\v} |u^{\v}|^3 + (\rho^\v)^{\gamma+1}\big)\,\dd r\dd t\nonumber\\[1mm]
	&\leq C(d,D,M,E_0,T),
\end{align*}
which, together with \eqref{5.29}, yields that
\begin{align}\label{5.31}
	\int_0^T\int_d^D \Big|\sqrt{\rho^{\v,b}} u^{\v,b} \mathbf{1}_{\{|u^{\v,b}|\leq k\}}
	- \sqrt{\rho^{\v}} u^{\v} \mathbf{1}_{\{|u^{\v}|\leq k\}}\Big|^2\,\dd r\dd t
	\longrightarrow 0 \,\,\,\,   \mbox{as $b\rightarrow \infty$}.
\end{align}
For $k\geq 1$, using \eqref{5.28} and Lemma \ref{lem2.6}, we have
\begin{align}\label{5.32}
&\int_0^T\int_d^D \Big( \big|\sqrt{\rho^{\v,b}} u^{\v,b} \mathbf{1}_{\{|u^{\v,b}|\geq k\}}\big|^2
+\big|\sqrt{\rho^{\v}} u^{\v} \mathbf{1}_{\{|u^{\v}|\geq k\}}\big|^2 \Big)\,\dd r\dd t\\
&\leq \frac{1}{k}\int_0^T\int_d^D \big(\rho^{\v,b} |u^{\v,b} |^3+\rho^{\v} |u^{\v} |^3\big)\,\dd r\dd t\nonumber\\
&\leq \frac{C(d,D,M,E_0,T)}{k}.\nonumber
\end{align}
Notice that
\begin{align}\label{5.30}
	&\int_0^T\int_d^D \Big|\sqrt{\rho^{\v,b}} u^{\v,b}-\sqrt{\rho^\v} u^{\v}\Big|^2\,\dd r\dd t\\
	&\leq\int_0^T\int_d^D \Big|\sqrt{\rho^{\v,b}} u^{\v,b} \mathbf{1}_{\{|u^{\v,b}|\leq k\}}
	- \sqrt{\rho^{\v}} u^{\v}\mathbf{1}_{\{|u^{\v}|\leq k\}}\Big|^2\,\dd r\dd t\nonumber\\
	&\quad + 2\int_0^T\int_d^D \Big|\sqrt{\rho^{\v,b}} u^{\v,b} \mathbf{1}_{\{|u^{\v,b}|\geq k\}}\Big|^2\,\dd r\dd t\nonumber\\
	&\quad +2\int_0^T\int_d^D\Big|\sqrt{\rho^{\v}} u^{\v} \mathbf{1}_{\{|u^{\v}|\geq k\}}\Big|^2\,\dd r\dd t.\nonumber
\end{align}
Substituting  \eqref{5.31}--\eqref{5.32} into \eqref{5.30}, we obtain
\begin{align}\nonumber
	\lim_{b\rightarrow\infty} \int_0^T\int_d^D \Big|\sqrt{\rho^{\v,b}} u^{\v,b}-\sqrt{\rho^\v} u^{\v}\Big|^2\,\dd r\dd t
	\leq \frac{C(d,D,M,E_0,T)}{k}\quad\,\, \mbox{for all $k\geq1$}.
\end{align}
Thus, by taking $k\rightarrow \infty$, we have proved that
$\frac{m^{\v,b}}{\sqrt{\rho^{\v,b}}}\longrightarrow \sqrt{\rho^{\v}} u^{\v}\equiv\frac{m^\v}{\sqrt{\rho^\v}}$
strongly in $L^2(0,T; L^2_{\rm loc})$. Therefore, the proof of Lemma \ref{lem5.6} is complete.
\end{proof}

Let  $(\rho^\v, m^\v)$ be the limit obtained above.
First, using \eqref{2.8}, \eqref{3.3-1},
Lemmas \ref{lem2.1}, \ref{lem2.2}, \ref{lem2.4}, \ref{lem2.6}, \ref{lem5.2}, and \ref{lem5.6},
Corollary 3.2, Fatou's lemma, and the lower semicontinuity, we have

\begin{proposition}\label{prop5.2}
Under assumptions \eqref{2.7}--\eqref{2.60-1},
for any fixed $\v$ and  $T>0$,  the limit functions $(\rho^\v, m^\v)=(\rho^\v, \rho^\v u^\v)$ satisfy
\begin{align}\label{5.7-1}
\begin{split}
&\hspace{0.5cm}\rho^{\v}(t,r)\geq0 \,\,\,\,\, a.e.,
\end{split}\\[2mm]
&\hspace{0.5cm} u^\v(t,r)=0, \,\,\big(\frac{m^\v}{\sqrt{\rho^\v}}\big)(t,r)=\sqrt{\rho^\v}(t,r)u^\v(t,r)=0\label{5.7} \\
&\hspace{6.0cm} a.e.\,\, \mbox{\, on $\{(t,r)\,:\, \rho^\v(t,r)=0\}$},\nonumber\\
\begin{split}
\hspace{0.5cm}\int_{0}^{\infty} \rho^\v(t,r)\,r^{n-1}\dd r\leq \frac{M}{\omega_n}\qquad\,\, \mbox{for all $t\geq0$},\label{3.36} 	
\end{split}\\[2mm]
&\int_0^\infty  \big(\frac12\rho^{\v} |u^{\v}|^2+(\rho^{\v})^{\gamma}\big)(t,r)\,r^{n-1}\dd r
 +\int_0^{\infty}r^{-n+1}\Big(\int_0^r \rho^\v(t,z) \,z^{n-1} \dd z\Big)^2\,\dd r\label{3.36b}\\
&\quad+\int_0^{\infty} \Big(\int_0^r \rho^\v z^{n-1} \dd z\Big)\rho^\v(t,r)\,r\dd r
+\v\int_0^t\int_0^\infty  (\rho^{\v} |u^{\v}|^2)(s,r)\, r^{n-3} \dd r \dd s \nonumber\\[1.5mm]
&\quad\leq C(M,E_0)\qquad \mbox{for all $t \geq0$},\nonumber\\[2mm]
&\hspace{0.2cm}\v^2 \int_0^\infty \big|(\sqrt{\rho^{\v}(t,r)})_r\big|^2\,r^{n-1}\dd r
+\v\int_0^T\int_{0}^\infty \big|\big((\rho^{\v}(s,r))^{\frac{\gamma}{2}}\big)_r\big|^2\,r^{n-1}\dd r\dd s\label{5.8-1}\\[1mm]
&\hspace{0.2cm}\quad\leq C(M,E_0,T) \qquad\,\,\mbox{for $t\in[0,T]$}, \nonumber\\[2mm]
\begin{split}
\int_0^T\int_d^D  \big(\rho^{\v}|u^{\v}|^3+(\rho^{\v})^{\gamma+\theta}+(\rho^{\v})^{\g+1}\big)(t,r)\, r^{n-1}\dd r \dd t
\leq C(d,D,M,E_0,T),\label{5.10-1}
\end{split}
\end{align}
where $[d,D]\Subset (0,\infty)$.
\end{proposition}

The following lemma is devoted to the convergence of the potential functions $\Phi^{\v,b}$.
\begin{lemma}\label{lem5.1}
For fixed $\v>0$,	there exists a function $\Phi^\v(t,\mathbf{x})=\Phi^{\v}(t,r)$ such that,
as $b\rightarrow\infty$ {\rm(}up to a subsequence{\rm)},
\begin{align}
	&\Phi^{\v,b} \rightharpoonup \Phi^\v \quad
	\mbox{weak-$\ast$ in $L^\infty(0,T; H^1_{\rm loc}(\mathbb{R}^n))$ and weakly in $L^2(0,T; H^1_{\rm loc}(\mathbb{R}^n))$},\label{5.2}\\
	&\Phi_r^{\v,b}(t,r)r^{n-1} \longrightarrow \Phi_r^\v(t,r)r^{n-1} =\kappa\int_0^r  \rho^{\v}(t,z)\, z^{n-1} \dd z\quad
	\mbox{in $C_{\rm loc}([0,T]\times[0,\infty))$},\label{5.3}	
\end{align}
and
\begin{align}\label{5.7-0}
	\|\Phi^\v(t)\|_{L^{\frac{2n}{n-2}}(\mathbb{R}^n)}+\|\nabla \Phi^\v(t)\|_{L^2(\mathbb{R}^n)}\leq C(M,E_0)\qquad\mbox{for $t\geq0$}.
\end{align}
Moreover, if $\g>\frac{2n}{n+2}$,
\begin{align}\label{5.4-1}
	\int_0^\infty  |(\Phi_r^{\v,b}-\Phi_r^{\v})(t,r)|^2\, r^{n-1} \dd r\longrightarrow 0\quad \mbox{as $b\to \infty\, $ {\rm (}\mbox{up to a subsequence}{\rm )}}.
\end{align}
\end{lemma}

\begin{proof}
The proof of \eqref{5.2} and \eqref{5.7-0} are direct by using Corollary \ref{cor3.1}.

We now prove \eqref{5.3}.  For any $D>0$ and $(t,r)\in[0,T]\times[0,D]$,
taking $b$ sufficiently large, then it follows from \eqref{2.24-1} that
\begin{align}\label{5.3-1}
	&\Big| \Phi_r^{\v,b}(t,r)r^{n-1} - \kappa\int_0^r  \rho^{\v}(t,z)\,z^{n-1}\dd z\Big|\\
	& = \left|\kappa\int_0^D  \big(\rho^{\v,b}-\rho^{\v}\big)(t,z)\,z^{n-1}\dd z \right|\nonumber\\
	&\leq C\left|\int_\sigma^D \big(\rho^{\v,b}-\rho^{\v}\big)(t,z)\,z^{n-1}\dd z \right|
	+C\left|\int_0^\sigma \big(\rho^{\v,b}-\rho^{\v}\big)(t,z)\,z^{n-1}\dd z \right|.\nonumber
\end{align}
Using \eqref{5.12}, we see that, for any fixed $\sigma>0$,
\begin{align}\label{5.3-2}
\lim_{b\rightarrow\infty} \sup_{t\in[0,T]}\left|\int_\sigma^D \big(\rho^{\v,b}-\rho^{\v}\big)(t,z)\,z^{n-1}\dd z \right|=0.
\end{align}
It follows from \eqref{5.7-1} and Lemma \ref{lem2.1} that
\begin{align*}
&\left|\int_0^\sigma \big(\rho^{\v,b}-\rho^{\v}\big)(t,z)\,z^{n-1}\dd z \right|\nonumber\\
&\leq C\Big(\int_0^\sigma \big((\rho^{\v,b})^{\gamma}+(\rho^{\v})^{\gamma}\big)\,z^{n-1} \dd r \Big)^{\frac1\gamma} \,
\Big(\int_0^\sigma z^{n-1} \dd r\Big)^{1-\frac1\gamma}\nonumber\\
&\leq C(M,E_0) \sigma^{n(1-\frac1\gamma)}\longrightarrow 0\qquad  \mbox{as $\sigma \rightarrow0$},
\end{align*}
which, together with \eqref{5.3-1}--\eqref{5.3-2}, yields that
\begin{equation}\nonumber
\sup_{[0,T]\times[0,D]}\Big| \Phi_r^{\v,b}(t,r)r^{n-1} - \kappa\int_0^r  \rho^{\v}(t,z)\,z^{n-1} \dd z\Big|\longrightarrow 0
\qquad \mbox{as $b\rightarrow\infty$},
\end{equation}
which leads to \eqref{5.3}.

For \eqref{5.4-1}, we first notice that
\begin{align}\nonumber
\frac{1}{r^{n-1}} \Big|\int_0^{r} (\rho^{\v,b}-\rho^\v)(t,z)\, z^{n-1} \dd z\Big|^2\leq C(M) r^{-n+1}\qquad \mbox{for $r>0$},
\end{align}
which yields that
\begin{align}\label{4.33}
\int_k^\infty  \frac{1}{r^{n-1}} \Big|\int_0^{r} (\rho^{\v,b}-\rho^\v)(t,z)\, z^{n-1} \dd z\Big|^2 \dd r\leq C(M) k^{-n+2}.
\end{align}
Using the H\"{o}lder inequality, we have
\begin{align}\label{4.34}
&\frac{1}{r^{n-1}} \Big|\int_0^{r} (\rho^{\v,b}-\rho^\v)(t,z)\, z^{n-1} \dd z\Big|^2\\
& \leq Cr^{-n+1} \left(\int_0^r \big((\rho^{\v,b})^\g+(\rho^{\v})^\g\big) z^{n-1} \dd z\right)^{\f2\g}\, r^{2n(1-\f1\g)}\nonumber\\
&\leq C(E_0,M) r^{n+1-\frac{2n}{\g}}.\nonumber
\end{align}
Since $n+1-\frac{2n}{\g}>-1$ for $\g>\frac{2n}{n+2}$, it follows from \eqref{5.12}, \eqref{4.34},
and Lebesgue's dominated convergence theorem that, for any given $k>0$,
\begin{align}\nonumber
\int_0^k \frac{1}{r^{n-1}} \Big|\int_0^{r} (\rho^{\v,b}-\rho^\v)(t,z)\,z^{n-1} \dd z\Big|^2 \dd r\longrightarrow 0\qquad \mbox{as $b\to \infty$},
\end{align}
which, together with \eqref{4.33},  yields  that
\begin{align*}
&\lim_{b\to\infty}\int_0^\infty  |(\Phi_r^{\v,b}-\Phi_r^{\v})(t,r)|^2\, r^{n-1} \dd r\\
& = \lim_{b\to\infty}\int_0^\infty \frac{1}{r^{n-1}} \Big|\int_0^{r} (\rho^{\v,b}-\rho^\v)(t,z)\, z^{n-1} \dd z\Big|^2 \dd r\nonumber\\
&\leq C(M) k^{-n+2}+\lim_{b\to\infty}\int_0^k \frac{1}{r^{n-1}} \Big|\int_0^{r} (\rho^{\v,b}-\rho^\v)(t,z)\, z^{n-1} \dd z\Big|^2 \dd r\\
& \leq C(M) k^{-n+2}.
\end{align*}
Then \eqref{5.4-1} follows by taking $k\to \infty$.
\end{proof}

\begin{remark}
The convergence result \eqref{5.4-1} is essential to prove the energy inequality for the case that $\kappa=1$ $($gaseous stars$)$.
Moreover, from \eqref{5.3}, it is direct to know that $\Phi^\v$ satisfies the Poisson equation in the classical sense except the origin:
\begin{equation}\nonumber
	\Delta \Phi^\v(t,\mathbf{x})=\kappa \rho^\v(t,\mathbf{x})\qquad \mbox{for}\  (t,\mathbf{x})\in [0,\infty)\times \mathbb{R}^n\backslash \{\mathbf{0}\}.
\end{equation}
\end{remark}

\begin{lemma}\label{lem4.8-1}
Let  $\g>1$ for $\kappa=-1$ $($plasmas$)$ and $\g>\frac{2n}{n+2}$ for $\kappa=1$ $($gaseous stars$)$. Then
\begin{align}\label{4.35}
	&\int_0^{\infty}\Big(\frac{1}{2}\Big|\frac{m^\v}{\sqrt{\rho^\v}}\Big|^2 + \rho^\v e(\rho^\v)\Big)(t,r)\,r^{n-1} \dd r
	-\frac{\kappa}{2} \int_0^\infty  |\Phi^\v(t,r)|^2\, r^{n-1}\dd r\\
	&\leq \int_0^{\infty}\Big(\frac{1}{2} \Big|\frac{m_0^\v}{\sqrt{\rho_0^\v}}\Big|^2 +\rho_0^\v e(\rho_0^\v)\Big)(r)\,r^{n-1} \dd r
	-\frac{\kappa}{2} \int_0^\infty  |\Phi_0^\v(r)|^2\, r^{n-1}\dd r.\nonumber
\end{align}
\end{lemma}

\begin{proof} For $\kappa=-1$ $($plasmas$)$,  \eqref{4.35} follows directly from  \eqref{2.13-1}, \eqref{2.23}, and Fatou's lemma.

For $\kappa=1$ $($gaseous stars$)$,   \eqref{4.35}  follows from \eqref{2.13-1}, \eqref{2.23}, \eqref{5.4-1}, and Fatou's lemma,
where the strong convergence of the gravitational potentials \eqref{5.4-1} plays a key role.
\end{proof}

Denote
\begin{equation*}
(\rho^\v,\M^\v, \Phi^{\v})(t,\mathbf{x}):=(\rho^\v(t,r),\  m^\v(t,r)\frac{\mathbf{x}}{r}, \  \Phi^\v(t,r)).
\end{equation*}
We show that $(\rho^\v,\M^\v,\Phi^\v)$ is a global weak solution of the Cauchy problem for CNSPEs \eqref{1.1}
in $\mathbb{R}^n$ in the sense of Definition \ref{definition-NSP}.

\begin{lemma}\label{lem5.8}
Let $0\leq t_1<t_2\leq T$, and let $\zeta(t,\mathbf{x})\in C^1([0,T]\times\mathbb{R}^n)$ be any smooth function with compact support.
Then
\begin{align}\label{5.39}
	&\int_{\mathbb{R}^n} \rho^{\v}(t_2,\mathbf{x}) \zeta(t_2,\mathbf{x})\, \dd\mathbf{x}\\
	& =\int_{\mathbb{R}^n} \rho^{\v}(t_1,\mathbf{x}) \zeta(t_1,\mathbf{x})\, \dd\mathbf{x}
	+\int_{t_1}^{t_2} \int_{\mathbb{R}^n} \big(\rho^{\v} \zeta_t + \M^{\v}\cdot\nabla\zeta\big)\,\dd\mathbf{x}\dd t.\nonumber
\end{align}
Moreover, the total mass is conserved:
\begin{equation}\label{3.27}
	\int_{\mathbb{R}^n} \rho^\v(t,\mathbf{x})\, \dd\mathbf{x}=\int_{\mathbb{R}^n} \rho_0^\v(\mathbf{x})\, \dd\mathbf{x}=M\qquad\, \mbox{for $t\geq0$}.
\end{equation}
\end{lemma}

\begin{proof}
Using \eqref{2.60}, we can choose sufficiently large $b\gg 1$ so that
${\rm supp}\,\zeta (t,\cdot) \subset B_{b/2}(\mathbf{0})$ for $t\in[0,T]$.
Then it follows from $\eqref{5.5-2}_1$ and a direct calculation that
\begin{align}\label{5.40}
0&=\int_{t_1}^{t_2} \int_{\mathbb{R}^n\backslash B_a(\mathbf{0})}
\big((\rho^{\v,b})_t +\mbox{div}\, \M^{\v,b}\big) \zeta(t,\mathbf{x})\, \dd\mathbf{x}\dd t\\
&=\int_{\mathbb{R}^n\backslash B_a(\mathbf{0})} \rho^{\v,b} \zeta\, \dd\mathbf{x}\bigg|_{t_1}^{t_2}
-\int_{t_1}^{t_2} \int_{\mathbb{R}^n\backslash B_a(\mathbf{0})} \big(\rho^{\v,b} \zeta_t + \M^{\v,b} \cdot\nabla\zeta\big)
\,\dd\mathbf{x}\dd t\nonumber\\
&=\int_{\mathbb{R}^n} \rho^{\v,b} \zeta\,\dd\mathbf{x}\bigg|_{t_1}^{t_2}
-\int_{t_1}^{t_2} \int_{\mathbb{R}^n} \big(\rho^{\v,b} \zeta_t + \M^{\v,b} \cdot\nabla\zeta\big)\, \dd\mathbf{x}\dd t,\nonumber
\end{align}
where we have used the fact that $(\rho^{\v,b},m^{\v,b})$ is extended by zero in $[0,T]\times[0,a)$.

Notice that, for $i=1,2$,
\begin{align}\label{5.41}
&\Big| \int_{\mathbb{R}^n} \big(\rho^{\v,b}-\rho^{\v}\big)(t_i,\mathbf{x}) \zeta(t_i,\mathbf{x})\, \dd\mathbf{x} \Big|\\
&\leq \Big| \int_{\mathbb{R}^n\backslash B_\s(\mathbf{0})} \big(\rho^{\v,b}-\rho^{\v}\big)(t_i,\mathbf{x}) \zeta(t_i,\mathbf{x})\,
\dd\mathbf{x} \Big|\nonumber\\
&\quad +\Big|\int_{B_\s(\mathbf{0})}\big(\rho^{\v,b}-\rho^{\v}\big)(t_i,\mathbf{x}) \zeta(t_i,\mathbf{x})\, \dd\mathbf{x}\Big|.\nonumber
\end{align}
We denote
\begin{align}\label{3.39-1}
\phi(t,r):=\int_{\partial B_1(\mathbf{0})}\zeta(t,r\bom)\, \dd\bom\in C^1_{0}([0,T]\times[0,\infty)),
\end{align}
which, together with \eqref{5.12},  yields that, for any fixed $\sigma>0$,
\begin{align}\label{5.42}
&\lim_{b\rightarrow\infty}\Big|\int_{\mathbb{R}^n\backslash B_\s(\mathbf{0})} \big(\rho^{\v,b}-\rho^{\v}\big)(t_i,\mathbf{x}) \zeta(t_i,\mathbf{x})\,
\dd\mathbf{x} \Big|\\
&=\lim_{b\rightarrow\infty}\Big|\int_{\sigma}^{\infty} \omega_n\big(\rho^{\v,b}-\rho^{\v}\big)(t_i,r) \phi(t_i,r)\, r^{n-1}\dd r\Big|=0.\nonumber
\end{align}
Using Lemma \ref{lem2.1} and \eqref{5.7-1}, we have
\begin{align}\label{5.42-1}
&\Big|\int_{B_\s(\mathbf{0})} \big(\rho^{\v,b}-\rho^{\v}\big)(t_i,\mathbf{x}) \zeta(t_i,\mathbf{x})\, \dd\mathbf{x} \Big|\\
&\leq C\|\zeta\|_{L^\infty}\Big\{\int_0^\sigma \big((\rho^{\v,b})^{\gamma}+(\rho^{\v})^{\gamma}\big)\,r^{n-1}\dd r \Big\}^{\frac1\gamma}
\, \Big\{\int_0^\sigma r^{n-1} \dd r\Big\}^{1-\frac1\gamma}\nonumber\\
&\leq C(M,E_0) \|\zeta\|_{L^\infty} \sigma^{n(1-\frac1\gamma)}\longrightarrow0\qquad \mbox{as $\sigma \rightarrow0$},\nonumber
\end{align}
which, together with \eqref{5.41} and \eqref{5.42}, leads to
\begin{align}\label{5.43}
\lim_{b\rightarrow\infty}\int_{\mathbb{R}^n} \rho^{\v,b}(t_i,\mathbf{x}) \zeta(t_i,\mathbf{x})\,\dd\mathbf{x}
=\int_{\mathbb{R}^n} \rho^{\v}(t_i,\mathbf{x}) \zeta(t_i,\mathbf{x})\,\dd\mathbf{x}\qquad\mbox{for $i=1,2$. }
\end{align}

From \eqref{3.39-1}, it is direct to show that
\begin{align}\label{3.43-1}
\phi_r(t,r)=\int_{\partial B_1(\mathbf{0})}\bom\cdot \nabla\zeta(t,r\bom)\,\dd\bom,
\end{align}
which, together with \eqref{5.12} and Lemma \ref{lem5.6}, implies that
\begin{align}\label{5.44}
&\lim_{b\rightarrow\infty}\int_{t_1}^{t_2} \int_{\mathbb{R}^n\backslash B_\s(\mathbf{0})}
\big(\rho^{\v,b} \zeta_t + \M^{\v,b} \cdot\nabla\zeta\big)\,\dd\mathbf{x}\dd t\\
&=\lim_{b\rightarrow\infty}\int_{t_1}^{t_2}\int_\sigma^\infty \big(\rho^{\v,b} \phi_t + m^{\v,b} \phi_r\big)\,\omega_n r^{n-1} \dd r\dd t\nonumber\\
&=\int_{t_1}^{t_2}\int_\sigma^\infty \big(\rho^{\v} \phi_t + m^{\v} \phi_r\big)\,\omega_n r^{n-1} \dd r\dd t\nonumber\\
&=\int_{t_1}^{t_2} \int_{\mathbb{R}^n\backslash B_\s(\mathbf{0})} \big(\rho^{\v} \zeta_t + \M^{\v} \cdot\nabla\zeta\big)\, \dd\mathbf{x}\dd t.\nonumber
\end{align}
Similar to those as in \eqref{5.42-1}, we have
\begin{align*}
&\Big|\int_{t_1}^{t_2} \int_{B_\s(\mathbf{0})}  \big(\rho^{\v,b}-\rho^{\v}\big) \zeta_t\, \dd\mathbf{x}\dd t\Big|
\leq C(E_0,T) \|\zeta_t\|_{L^\infty}\,\sigma^{n(1-\frac1\gamma)},\\[2mm]
&\Big|\int_{t_1}^{t_2} \int_{B_\s(\mathbf{0})} \big(\M^{\v,b}-\M^\v\big) \cdot\nabla \zeta\, \dd\mathbf{x}\dd t\Big|\\
&\leq C\|\nabla\zeta\|_{L^\infty}\Big\{\int_{t_1}^{t_2}\int_0^\sigma \big(\rho^{\v,b}+\rho^{\v}\big)(t,r)\,r^{n-1}\dd r\dd t\Big\}^{\frac12}\\
&\quad\times \Big\{\int_{t_1}^{t_2}
\int_0^\sigma \big(\rho^{\v,b}|u^{\v,b}|^2+\rho^{\v}|u^{\v}|^2\big)(t,r)\,r^{n-1}\dd r\dd t\Big\}^{\frac12}\\
&\leq C(M,E_0,T) \|\nabla\zeta\|_{L^\infty} \sigma^{\frac{n}{2}(1-\frac1\gamma)},
\end{align*}
which, together with   \eqref{5.44}, yields that
\begin{align}\label{5.47}
\lim_{b\rightarrow\infty}\int_{t_1}^{t_2} \int_{\mathbb{R}^n} \big(\rho^{\v,b} \zeta_t + \M^{\v,b} \cdot\nabla\zeta\big)\, \dd\mathbf{x}\dd t
=\int_{t_1}^{t_2} \int_{\mathbb{R}^n} \big(\rho^{\v} \zeta_t + \M^{\v} \cdot\nabla\zeta\big) \dd\mathbf{x}\dd t.
\end{align}
Combining \eqref{5.40} with \eqref{5.43} and \eqref{5.47}, we conclude \eqref{5.39}.

Finally, we prove the  conservation of mass \eqref{3.27}.
We take smooth test functions $\zeta(t,\mathbf{x})=\phi_k(r)$ in \eqref{5.39} with
\begin{align}
&\phi_k(r)=
\begin{cases}
	1 \qquad &\mbox{for $r\in[0,k]$}, \\
	\mbox{smooth} \quad &\mbox{for $r\in[k,k+1]$},\\
	0 \qquad &\mbox{for $r\in[k+1,\infty)$},
\end{cases} \label{3.37}\\[2mm]
&|\phi_k'(r)|\leq C \qquad \mbox{for all $r\in [0,\infty)$}, \label{3.38}
\end{align}
where $C>0$ is a constant independent of $k$.
Now it follows from  \eqref{5.39} and \eqref{3.37}--\eqref{3.38} that
\begin{equation}\label{3.39}
\int_{\mathbb{R}^n} \rho^\v(t,\mathbf{x}) \phi_k\, \dd\mathbf{x}
=\int_{\mathbb{R}^n}\rho_0^\v(\mathbf{x}) \phi_k\, \dd\mathbf{x}
+\int_0^t\int_{k\leq |\mathbf{x}|\leq k+1} \M^\v\cdot \nabla\phi_k\,\dd\mathbf{x}\dd s.
\end{equation}
From \eqref{1.6-1},  \eqref{3.36}, and Lebesgue's dominated convergence theorem, we have
\begin{align}\label{3.40}
\lim_{k\rightarrow\infty} \int_{\mathbb{R}^n} (\rho^\v(t,\mathbf{x}),\,\rho_0^\v(\mathbf{x})) \phi_k\,\dd\mathbf{x}
=\int_{\mathbb{R}^n} (\rho^\v(t,\mathbf{x}),\,\rho_0^\v(\mathbf{x})) \,\dd\mathbf{x}.
\end{align}
Since
\begin{align*}
&\left|\int_0^T\int_{k\leq |\mathbf{x}|\leq k+1} \M^\v\cdot \nabla\phi_k\,\dd\mathbf{x}\dd t\right|\\
&=\left|\int_0^T\int_k^{k+1} \rho^\v u^\v \phi_{k}'(r)\,\omega_n\, r^{n-1} \dd r\dd t\right|\nonumber\\
&\leq  C\left\{\int_0^T\int_k^{k+1} \rho^\v \,r^{n-1}\dd r\dd t \right\}^{\frac12}\,
\left\{\int_0^T\int_k^{k+1} \rho^\v |u^\v|^2\, r^{n-1} \dd r\dd t \right\}^{\frac12}\longrightarrow 0\, \,\,\,\, \mbox{as $k\rightarrow\infty$},
\end{align*}
together with \eqref{3.39}--\eqref{3.40}, we conclude \eqref{3.27}.
\end{proof}

\begin{lemma}\label{lem5.9}
Let $\bp(t,\mathbf{x})\in \left(C^2_0([0,T]\times \mathbb{R}^n)\right)^n$ be any smooth function with compact support
so that $\bp(T,\mathbf{x})={\bf 0}$. Then
\begin{align}
	&\int_{\mathbb{R}_+^{n+1}} \Big\{\M^{\v} \cdot\partial_t\bp +\frac{\M^{\v}}{\sqrt{\rho^{\v}}} \cdot \big(\frac{\M^{\v}}{\sqrt{\rho^{\v}}}\cdot \nabla\big)\bp
	+p(\rho^{\v}) \,\mbox{\rm div}\,\bp -\rho^\v \nabla\Phi^\v\cdot \bp\Big\}\,\dd\mathbf{x} \dd t\label{5.48}\\
	&\quad +\int_{\mathbb{R}^n} \M_0^\v\cdot \bp(0,\mathbf{x})\,\dd\mathbf{x}
	\nonumber\\
	&=-\v\int_{\mathbb{R}_+^{n+1}}
	\Big\{\frac{1}{2}\M^{\v}\cdot \big(\Delta \bp+\nabla\mbox{\rm div}\,\bp \big)
	+ \frac{\M^{\v}}{\sqrt{\rho^{\v}}} \cdot \big(\nabla\sqrt{\rho^{\v}}\cdot \nabla\big)\bp\nonumber\\
	&\hspace{5.5cm}\,\,\,\,+ \nabla\sqrt{\rho^{\v}}  \cdot \big(\frac{\M^{\v}}{\sqrt{\rho^{\v}}}\cdot \nabla\big)\bp\Big\}\,\dd\mathbf{x}\dd t\nonumber\\
	&=\sqrt{\v}\int_{\mathbb{R}_+^{n+1}}
	\sqrt{\rho^{\v}} \Big\{V^{\v}  \frac{\mathbf{x}\otimes\mathbf{x}}{r^2}
	+\frac{\sqrt{\v}}{r}\frac{m^\v}{\sqrt{\rho^\v}}\big(I_{n\times n}-\frac{\mathbf{x}\otimes\mathbf{x}}{r^2}\big)\Big\}: \nabla\bp\, \dd\mathbf{x}\dd t,\nonumber
\end{align}
where $V^{\v}(t,r)\in L^2(0,T; L^2(\mathbb{R}^n))$ is a function such that
$$
\displaystyle\int_0^T\int_{\mathbb{R}^n} |V^{\v}(t,\mathbf{x})|^2\, \dd\mathbf{x}\dd t\leq C(E_0,M)
\quad \mbox{for some $C(E_0,M)>0$ independent of $T>0$}.
$$
\end{lemma}

\begin{proof}
For any given $\sigma\in (0,1]$, let $\chi_\sigma(r)\in C^\infty(\mathbb{R})$ be a cut-off function satisfying
\begin{align}\label{3.49-1}
\chi_\sigma(r)=0\,\,\, \mbox{for $r\leq \sigma$};\,\,\,\,\,
\chi_\sigma(r)=1\,\,\, \mbox{for $r\geq 2\sigma$}; \,\,\,\,\,
|\chi_\sigma'(r)|+\sigma|\chi_\sigma''(r)|\leq \frac{C}{\sigma}\,\,\, \mbox{for $r\in\R$}.
\end{align}
Denote $\Psi_\sigma(t,\textbf{x})=\bp(t,\textbf{x}) \chi_\sigma(|\textbf{x}|)$.
Taking $b\gg1$ large enough so that  $a=b^{-1}\leq \sigma$,
then it follows from $\eqref{5.5-2}_2$ and integration by parts that
\begin{align*}
&\int_{\mathbb{R}_+^{n+1}} \Big\{\M^{\v,b} \cdot\partial_t\Psi_\sigma +\frac{\M^{\v,b}}{\sqrt{\rho^{\v,b}}}\cdot \big(\frac{\M^{\v,b}}{\sqrt{\rho^{\v,b}}}\cdot \nabla\big)\Psi_{\sigma}
+ p(\rho^{\v,b})\,\mbox{\rm div}\,\Psi_{\sigma} \Big\}\,  \dd\mathbf{x}\dd t
\nonumber\\
&\quad +\int_{\mathbb{R}^n} \M_0^{\v,b}\cdot \Psi_\s(0,\mathbf{x})\,\dd\textbf{x}\nonumber\\
&=\int_{\mathbb{R}_+^{n+1}} \rho^{\v,b} \nabla\Phi^{\v,b}\cdot \Psi_\sigma\,\dd\mathbf{x}\dd t+R^{\v,b},
\end{align*}
where
\begin{align*}
R^{\v,b}&=-\v\int_{\mathbb{R}_+^{n+1}} \Big\{\frac12\M^{\v,b} \cdot \big(\Delta \Psi_{\sigma}+\nabla\mbox{div}\,\Psi_{\sigma} \big)
+\frac{\M^{\v,b}}{\sqrt{\rho^{\v,b}}} \cdot \big(\nabla\sqrt{\rho^{\v,b}}\cdot \nabla\big)\Psi_{\sigma} \nonumber\\
&\qquad\qquad\quad\,\,\,+\nabla\sqrt{\rho^{\v,b}}  \cdot \big(\frac{\M^{\v,b}}{\sqrt{\rho^{\v,b}}}\cdot \nabla\big)\Psi_{\sigma}\Big\}\, \dd\mathbf{x} \dd t\\
&=\v\int_{\mathbb{R}_+^{n+1}} \sqrt{\rho^{\v,b}}\Big(\sqrt{\rho^{\v,b}} D(\frac{\M^{\v,b}}{\rho^{\v,b}})\Big): \nabla\Psi_{\sigma}\, \dd\mathbf{x} \dd t.
\end{align*}

For the term involving the potentials, using \eqref{5.12} and \eqref{5.2}, we have
\begin{align}\label{4.52-1}
\lim_{b\rightarrow \infty}\int_{\mathbb{R}_+^{n+1}} \rho^{\v,b} \nabla\Phi^{\v,b}\cdot  \Psi_\sigma\,\dd\mathbf{x}\dd t
=\int_{\mathbb{R}_+^{n+1}} \rho^{\v} \nabla\Phi^{\v}\cdot  \Psi_\sigma\,\dd\mathbf{x}\dd t.
\end{align}

For the convergence of the viscous term, it follows from \eqref{5.5-1} and a direct calculation that
\begin{align*}
\partial_i \Big(\frac{\mathcal{M}^{\v,b}_j}{\rho^{\v,b}}\Big)=u^{\v,b}_r \frac{x_ix_j}{r^2}+\frac{u^{\v,b}}{r} \big(\d_{ij}-\frac{x_i x_j}{r^2}\big).
\end{align*}
Thus, using Lemma \ref{lem2.1},  there exists a function $V^\v(t,r)$ so that
\begin{align}\label{5.52-1}
\sqrt{\v}\sqrt{\rho^{\v,b}} D\Big(\frac{\mathcal{M}^{\v,b}_j}{\rho^{\v,b}}\Big)\, {\longrightharpoonup}\,
V^{\v} \frac{{\mathbf{x}}\otimes{\mathbf{x}}}{r^2}
+\frac{\sqrt{\v}}{r}\sqrt{\rho^\v}u^{\v}\Big(I_{n\times n}-\frac{{\mathbf{x}}\otimes{\mathbf{x}}}{r^2}\Big)
\end{align}
in  $ L^2(0,T; (L^2(B_{\s^{-1}}(\mathbf{0})\backslash B_\sigma(\mathbf{0})))^{n\times n}) $ as  $b\rightarrow\infty$ for any given $\sigma>0$, and
\begin{align}\label{5.52-2}
\int_0^T\int_{\mathbb{R}^n} \left| V^\v\right|^2 \dd\mathbf{x} \dd t \leq C (E_0,M).
\end{align}

Denote
\begin{align}
\phi_{1\sigma}(t,r):=\int_{\partial B_1(\mathbf{0})}\big\{\bom\cdot(\Delta \Psi_{\sigma})(t,r\bom)
+\bom\cdot(\nabla\mbox{div}\,\Psi_{\sigma})(t,r\bom) \big\}\, \dd\bom.\nonumber
\end{align}
Then $\phi_{1\sigma}\in C_0([0,T]\times (0,\infty))$. Hence,  using Lemma \ref{lem5.6}, we find that, as $b\rightarrow\infty$,
\begin{align}\label{3.55-1}
&\int_{\mathbb{R}_+^{n+1}} \M^{\v,b} \cdot \big\{\Delta \Psi_{\sigma}+\nabla\mbox{div}\,\Psi_{\sigma} \big\}\, \dd\mathbf{x}\dd t
=\int_{\mathbb{R}_+^{2}} m^{\v,b} \phi_{1\sigma}\,\omega_n r^{n-1} \dd r \dd t\\
&\longrightarrow \int_{\mathbb{R}_+^{2}}  m^{\v} \phi_{1\sigma}\,\omega_n r^{n-1} \dd r \dd t
=\int_{\mathbb{R}_+^{n+1}} \M^{\v} \cdot \big\{\Delta \Psi_{\sigma}+\nabla\mbox{div}\,\Psi_{\sigma} \big\}\, \dd\mathbf{x}\dd t.\nonumber
\end{align}
Similarly, using Lemmas \ref{lem5.2} and \ref{lem5.6},  we see that, as $b\rightarrow\infty$,
\begin{align}\label{3.56-1}
& \int_{\mathbb{R}^{n+1}_+} \Big\{\frac{\M^{\v,b}}{\sqrt{\rho^{\v,b}}} \cdot \big(\nabla\sqrt{\rho^{\v,b}}\cdot \nabla\big)\Psi_{\sigma}
+ \nabla\sqrt{\rho^{\v,b}}  \cdot \big(\frac{\M^{\v,b}}{\sqrt{\rho^{\v,b}}}\cdot \nabla\big)\Psi_{\sigma}\Big\}\, \dd\mathbf{x} \dd t\\
&\longrightarrow \int_{\mathbb{R}^{n+1}_+} \Big\{ \frac{\M^{\v}}{\sqrt{\rho^{\v}}}\cdot\big(\nabla\sqrt{\rho^{\v}}\cdot \nabla\big)\Psi_{\sigma}
+ \nabla\sqrt{\rho^{\v}}  \cdot \big(\frac{\M^{\v}}{\sqrt{\rho^{\v}}}\cdot \nabla\big)\Psi_{\sigma} \Big\}\,\dd\mathbf{x}  \dd t.\nonumber
\end{align}
Combining \eqref{5.52-1}  and \eqref{3.55-1}--\eqref{3.56-1}, we obtain that, as  $b\rightarrow\infty$,
\begin{align}\label{3.57-2}
R^{\v,b}\longrightarrow &-\v\int_{\mathbb{R}_+^{n+1}}
\Big\{\frac{1}{2}\M^{\v}\cdot \big(\Delta \Psi_{\sigma}+\nabla\mbox{div}\,\Psi_{\sigma} \big)\\
&\hspace{1.8cm} +\Big(\frac{\M^{\v}}{\sqrt{\rho^{\v}}}\cdot \big(\nabla\sqrt{\rho^{\v}}\cdot \nabla\big)
  + \nabla\sqrt{\rho^{\v}}\cdot \big(\frac{\M^{\v}}{\sqrt{\rho^{\v}}}\cdot \nabla\big)\Big)\Psi_{\sigma}\Big\}\,\dd\mathbf{x}\dd t
\nonumber\\
&=\sqrt{\v}\int_{\mathbb{R}_+^{n+1}}\sqrt{\rho^{\v}} \Big\{V^{\v}  \frac{\mathbf{x}\otimes\mathbf{x}}{r^2}
+\frac{\sqrt{\v}}{r}\sqrt{\rho^\v}u^{\v}\Big(I_{n\times n}-\frac{\mathbf{x}\otimes\mathbf{x}}{r^2}\Big)\Big\}: \nabla\Psi_{\sigma}\, \dd\mathbf{x} \dd t.\nonumber
\end{align}

Also, by similar arguments as in \eqref{3.55-1}, using Lemma \ref{lem5.2},   Corollary \ref{cor5.3},
and Lemma \ref{lem5.6}, we have
\begin{align*}
&\int_{\mathbb{R}_+^{n+1}} \Big\{\M^{\v,b} \cdot\partial_t\Psi_\sigma
+\frac{\M^{\v,b}}{\sqrt{\rho^{\v,b}}}\cdot \big(\frac{\M^{\v,b}}{\sqrt{\rho^{\v,b}}}\cdot \nabla\big)\Psi_{\sigma}
+ p(\rho^{\v,b})\,\mbox{\rm div}\,\Psi_{\sigma} \Big\}\,\dd\mathbf{x}\dd t\\
&\qquad\quad+\int_{\mathbb{R}^n} \M_0^{\v,b}\cdot \Psi_\s(0,\mathbf{x})\,\dd\mathbf{x}\\
&\longrightarrow \int_{\mathbb{R}_+^{n+1}} \Big\{\M^{\v} \cdot\partial_t\Psi_\sigma
+\frac{\M^{\v}}{\sqrt{\rho^{\v}}}\cdot \big(\frac{\M^{\v}}{\sqrt{\rho^{\v}}}\cdot \nabla\big)\Psi_{\sigma}
+ p(\rho^{\v})\,\mbox{\rm div}\,\Psi_{\sigma} \Big\}\,\dd\mathbf{x}\dd t\\
&\qquad\quad+\int_{\mathbb{R}^n} \M_0^{\v}\cdot \Psi_\s(0,\mathbf{x})\,\dd\mathbf{x}
\end{align*}
as $b\rightarrow \infty$, which, together with \eqref{4.52-1} and  \eqref{3.57-2}, yields that
\begin{align}\label{3.59-1}
&\int_{\mathbb{R}_+^{n+1}} \Big\{\M^{\v} \cdot\partial_t\Psi_\sigma
+\frac{\M^{\v}}{\sqrt{\rho^{\v}}}\cdot \big(\frac{\M^{\v}}{\sqrt{\rho^{\v}}}\cdot \nabla\big)\Psi_{\sigma}
+ p(\rho^{\v})\,\mbox{\rm div}\,\Psi_{\sigma}-\rho^{\v} \nabla\Phi^{\v}\cdot  \Psi_\sigma \Big\}\,\dd\mathbf{x}\dd t\\
&\quad +\int_{\mathbb{R}^n} \M_0^{\v}\cdot \Psi_\s(0,\mathbf{x})\,\dd\mathbf{x} \nonumber\\
&=-\v\int_{\mathbb{R}_+^{n+1}}
\Big\{\frac{1}{2}\M^{\v}\cdot \big(\Delta \Psi_{\sigma}+\nabla\mbox{div}\,\Psi_{\sigma} \big)\nonumber\\
&\qquad\qquad\quad\,\,\,\,\, +\Big(\frac{\M^{\v}}{\sqrt{\rho^{\v}}}\cdot \big(\nabla\sqrt{\rho^{\v}}\cdot \nabla\big)
+ \nabla\sqrt{\rho^{\v}}\cdot \big(\frac{\M^{\v}}{\sqrt{\rho^{\v}}}\cdot \nabla\big)\Big)\Psi_{\sigma}\Big\}\,\dd\mathbf{x}\dd t\nonumber\\
&=\sqrt{\v}\int_{\mathbb{R}_+^{n+1}}\sqrt{\rho^{\v}} \Big\{V^{\v}  \frac{\mathbf{x}\otimes\mathbf{x}}{r^2}
+\frac{\sqrt{\v}}{r}\sqrt{\rho^\v}u^{\v}\Big(I_{n\times n}-\frac{\mathbf{x}\otimes\mathbf{x}}{r^2}\Big)\Big\}: \nabla\Psi_{\sigma}\, \dd\mathbf{x} \dd t.\nonumber
\end{align}

Next, we consider the limit: $\sigma\rightarrow0$ in \eqref{3.59-1}. First, we define
\begin{align}\label{3.63-2}
	\varphi(t,r):=\int_{\partial B_1(\mathbf{0})} \bom\cdot\bp(t,r\bom)\,\dd\bom
  &=\frac{1}{r^{n-1}}\int_{\partial B_r(\mathbf{0})} \bom\cdot\bp(t,\mathbf{y})\,\dd S_{\mathbf{y}}\\
	&=\frac{1}{r^{n-1}}\int_{B_r(\mathbf{0})} \mbox{div}\,\bp(t,\mathbf{y})\,\dd\mathbf{y},\nonumber
\end{align}
which implies that
\begin{align}\label{3.63-1}
|\varphi(t,r)|\leq C(\|\bp\|_{C^1}) r;
\end{align}
also see \cite{Jiang-Zhang-2001, Schrecker}.
For the term involving the potential, we notice from \eqref{5.3} and  \eqref{3.63-1} that
\begin{align*}
\Big|\rho^{\v}\,(\Phi_r^\v r^{n-1})\,\varphi\Big|
\leq C(\|\bp\|_{C^1}) \rho^{\v}(t,r)r \int_0^r  \rho^{\v}(t,z)\,z^{n-1}\dd z
\end{align*}
for $(t,r)\in [0,\infty)\times[0,\infty)$, which, together with Lebesgue's dominated convergence theorem
and \eqref{3.36b}, yields that
\begin{align}\label{3.63-5}
\lim_{\sigma\rightarrow0}  \int_{\mathbb{R}^{n+1}_+} \rho^{\v} \nabla\Phi^{\v}\cdot  \Psi_\sigma\,\dd\mathbf{x}\dd t
&=\lim_{\sigma\rightarrow0}  \int_{\R_+^2} \rho^{\v} \Phi_r^\v \,\varphi \chi_{\sigma}(r)\,\omega_n r^{n-1}\dd r\dd t\\
&= \int_{\R_+^2} \rho^{\v} \Phi_r^\v\,\varphi\, \omega_n r^{n-1}\dd r\dd t\nonumber\\
&=\int_{\mathbb{R}^{n+1}_+} \rho^{\v} \nabla\Phi^{\v}\cdot  \bp\, \dd\mathbf{x} \dd t.\nonumber	
\end{align}

Using \eqref{3.63-2}, Lebesgue's dominated convergence theorem, and Proposition \ref{prop5.2}, we have
\begin{align}\label{3.61}
	&\lim_{\s\rightarrow0}\Big\{\int_{\mathbb{R}^{n+1}_+} \M^{\v} \cdot\partial_t\Psi_\sigma\,\dd\mathbf{x} \dd t
	+ \int_{\mathbb{R}^n} \M^\v_0\cdot \Psi_\s(0,\mathbf{x})\,\dd\textbf{x}\Big\}\\
	&=\lim_{\s\rightarrow0}\Big\{ \int_{\R_+^2}  m^{\v} \, \partial_t\varphi \, \chi_\sigma(r)\,\omega_n r^{n-1}  \dd r \dd t
	+\int_{0}^\infty m_0^\v\, \varphi(0,r) \chi_\sigma(r)\,\omega_n  r^{n-1} \dd r\Big\}\nonumber\\
	&= \int_{\R_+^2}  m^{\v} \, \partial_t\varphi \,\omega_n r^{n-1}  \dd r \dd t
	+\int_{0}^\infty m_0^\v\, \varphi(0,r)\, \omega_n r^{n-1} \dd r\nonumber\\
	&=\int_{\mathbb{R}^{n+1}_+} \M^{\v} \cdot\partial_t\bp\, \dd\mathbf{x}\dd t + \int_{\mathbb{R}^n} \M^\v_0\cdot \bp(0,\mathbf{x})\,\dd\mathbf{x}.\nonumber
\end{align}

Employing \eqref{3.63-1} and Proposition \ref{prop5.2}, we have
\begin{align}\label{3.64}
&\Big|\int_{\mathbb{R}^{n+1}_+}  \big(\frac{|m^\v|^2}{\rho^\v} +p(\rho^\v)\big)\,\bp\cdot \frac{\mathbf{x}}{r}\, \chi'_\sigma(r)\, \dd\mathbf{x} \dd t\Big|\\
&\leq C\int_0^T\int_{\sigma}^{2\sigma} \big(\frac{|m^\v|^2}{\rho^\v} +p(\rho^\v)\big)\big|\varphi(t,r) \chi'_\sigma(r)\big|\, r^{n-1}\dd r \dd t\nonumber\\
&\leq C\int_0^T\int_{\sigma}^{2\sigma} \big(\frac{|m^\v|^2}{\rho^\v} +p(\rho^\v)\big)\, r^{n-1}\dd r \dd t \longrightarrow 0\qquad  \mbox{as $\sigma\rightarrow0$},\nonumber
\end{align}
\begin{align}\label{3.66}
&\Big|\v\int_{\mathbb{R}^{n+1}_+} \frac{m^{\v}}{\sqrt{\rho^{\v}}} (\sqrt{\rho^{\v}})_r \bp\cdot \frac{\mathbf{x}}{r} \chi'_\sigma(r)\,\dd\mathbf{x} \dd t\Big|\\
&\leq C\v\int_0^T\int_{\sigma}^{2\sigma} \big|\frac{m^{\v}}{\sqrt{\rho^{\v}}}(\sqrt{\rho^{\v}})_r \varphi(t,r) \chi'_\sigma(r)\big|\, r^{n-1}\dd r\dd t\nonumber\\
&\leq C\int_0^T\int_{\sigma}^{2\sigma} \big(\frac{|m^{\v}|^2}{\rho^{\v}}+ \v^2|(\sqrt{\rho^{\v}})_r|^2 \big)\, r^{n-1}\dd r\dd t\longrightarrow 0
\qquad \mbox{as $\sigma\rightarrow0$},\nonumber
\end{align}
and
\begin{align}\label{3.67}
&\Big|\sqrt{\v} \int_{\mathbb{R}^{n+1}_+}\chi'_\sigma(r)\sqrt{\rho^{\v}} \Big\{V^{\v}  \frac{\mathbf{x}\otimes\mathbf{x}}{r^2}
+\frac{\sqrt{\v}}{r}\sqrt{\rho^\v}u^{\v}\Big(I_{n\times n}-\frac{\mathbf{x}\otimes\mathbf{x}}{r^2}\Big)\Big\}: \Big(\bp\otimes \frac{\mathbf{x}}{r}\Big)\,
\dd\mathbf{x}\, \dd t\Big|\\
&\quad=\Big|\sqrt{\v}\int_0^T\int_{\mathbb{R}^n} \chi'_\sigma(r)\sqrt{\rho^{\v}}V^{\v}   \bp\cdot  \frac{\mathbf{x}}{r}\,\dd\mathbf{x} \dd t\Big|\nonumber\\
&\leq C\Big|\sqrt{\v}\int_0^T\int_{\sigma}^{2\sigma} \sqrt{\rho^{\v}}V^{\v}\, r^{n-1}\dd r \dd t\Big|\longrightarrow 0
\qquad \mbox{as $\sigma\rightarrow0$}.\nonumber
\end{align}
Using \eqref{3.64}--\eqref{3.67},  Lebesgue's dominated convergence theorem, and Proposition \ref{prop5.2}, we obtain
\begin{align}
	&\lim_{\s\rightarrow0} \int_{\mathbb{R}^{n+1}_+}
	\Big\{\frac{\M^\v}{\sqrt{\rho^{\v}}} \cdot \big(\frac{\M^\v}{\sqrt{\rho^{\v}}}\cdot \nabla\big)\Psi_{\sigma}
   +p(\rho^{\v})\,\mbox{div}\,\Psi_{\sigma}\Big\}
	\,\dd\mathbf{x} \dd t\label{3.68} \\
	&\quad =\int_{\mathbb{R}^{n+1}_+}\Big\{\frac{\M^\v}{\sqrt{\rho^{\v}}}\cdot \big(\frac{\M^\v}{\sqrt{\rho^{\v}}} \cdot \nabla\big)\bp
	+p(\rho^{\v})\,\mbox{div}\,\bp\Big\}\, \dd\mathbf{x} \dd t,\nonumber\\[2mm]
	&\lim_{\s\rightarrow0}\v \int_{\mathbb{R}^{n+1}_+}
	\Big\{\frac{\M^\v}{\sqrt{\rho^{\v}}} \cdot \big(\nabla\sqrt{\rho^{\v}}\cdot \nabla\big)\Psi_{\sigma}
	+ \nabla\sqrt{\rho^{\v}}  \cdot \big(\frac{\M^\v}{\sqrt{\rho^{\v}}}\cdot \nabla\big)\Psi_{\sigma}\Big\}\, \dd\mathbf{x} \dd t\label{3.69}\\
	&\quad =\v \int_{\mathbb{R}^{n+1}_+}\Big\{\frac{\M^\v}{\sqrt{\rho^{\v}}} \cdot \big(\nabla\sqrt{\rho^{\v}}\cdot \nabla\big)\bp
	+ \nabla\sqrt{\rho^{\v}}  \cdot \big(\frac{\M^\v}{\sqrt{\rho^{\v}}}\cdot \nabla\big)\bp\Big\}\, \dd\mathbf{x} \dd t,\nonumber\\[2mm]
	&\lim_{\s\rightarrow0}\sqrt{\v} \int_{\mathbb{R}^{n+1}_+}\sqrt{\rho^{\v}} \Big\{V^{\v}  \frac{\mathbf{x}\otimes\mathbf{x}}{r^2}
	+\frac{\sqrt{\v}}{r}\sqrt{\rho^\v}u^{\v}\big(I_{n\times n}-\frac{\mathbf{x}\otimes\mathbf{x}}{r^2}\big)\Big\}: \nabla\Psi_{\sigma}\, \dd\mathbf{x} \dd t\label{3.70}\\
	&\quad =\sqrt{\v}\int_{\mathbb{R}^{n+1}_+} \sqrt{\rho^{\v}} \Big\{V^{\v}  \frac{\mathbf{x}\otimes\mathbf{x}}{r^2}
	+\frac{\sqrt{\v}}{r}\sqrt{\rho^\v}u^{\v}\big(I_{n\times n}-\frac{\mathbf{x}\otimes\mathbf{x}}{r^2}\big)\Big\}: \nabla\bp\,  \dd\mathbf{x} \dd t.\nonumber
\end{align}

We notice that
\begin{align}\label{3.70-1}
	&\Delta (\Psi_{\sigma})_i=\chi_\sigma(r) \Delta \psi_i +2\nabla\psi_i \cdot \nabla \chi_\sigma(r)+\psi_i \Delta \chi_\sigma(r),\\
	&\partial_i(\mbox{div}\,\Psi_{\sigma})= \chi_\sigma(r) \partial_i (\mbox{div}\,\bp)+\mbox{div}\,\bp \partial_i \chi_\sigma(r)
	+\partial_i \bp\cdot \nabla \chi_\sigma(r)\nonumber\\
	&\qquad\qquad\quad\,+\frac{x_i}{r}\chi''_\sigma(r) \bp\cdot \frac{\mathbf{x}}{r}+\chi'_{\sigma}(r) \big( \bp\cdot \frac{\nabla x_i}{r}-\bp \cdot \frac{\mathbf{x}}{r} \frac{x_i}{r^2}\big).\nonumber
\end{align}
It follows from \eqref{3.63-1} and Proposition \ref{prop5.2} that
\begin{align}\label{3.71}
&\left|\sum_{i=1}^n\v\int_{\mathbb{R}^{n+1}_+} m^{\v}\,\frac{x_i}{r} \,\Big\{ 2\nabla\psi_i \cdot \nabla \chi_\sigma
+\psi_i \Delta \chi_\sigma+\mbox{div}\,\bp\, \partial_i \chi_\sigma(r)
+\partial_i \bp\cdot \nabla \chi_\sigma(r) \right.\\
&\qquad\qquad\qquad\quad\,\,\,\,   \left.+\frac{x_i}{r}\chi''_\sigma(r) \bp\cdot \frac{\mathbf{x}}{r}
+\chi'_{\sigma}(r) \big( \bp\cdot \frac{\nabla x_i}{r}-\bp \cdot \frac{\mathbf{x}}{r} \frac{x_i}{r^2}\big) \Big\}\, \dd\mathbf{x} \dd t\right|\nonumber\\
&\leq C(\|\bp\|_{C^1})\int_0^T\int_{\sigma}^{2\sigma} \v| m^{\v}| \Big(|\chi'_\sigma(r)|+\frac1{r}\varphi(r)|\chi'_\sigma(r)|+\varphi(r)|\chi''_\sigma(r)|
\Big)\,r^{n-1}\dd r \dd t\nonumber\\
&\leq C(\|\bp\|_{C^1})\int_0^T\int_{\sigma}^{2\sigma} \v|m^{\v}|\, r^{n-2} \dd r \dd t\nonumber\\
&\leq C(\|\bp\|_{C^1})\left\{\int_0^T\int_{\sigma}^{2\sigma} \rho^\v\, r^{n-1} \dd r \dd t \right\}^{\frac12} \,
\left\{\v\int_0^T\int_{\sigma}^{2\sigma} \frac{|m^{\v}|^2}{\rho^\v}\, r^{n-3} \dd r \dd t \right\}^{\frac12}\longrightarrow 0\nonumber
\end{align}
as $\sigma\rightarrow0$.
Thus, using \eqref{3.70-1}--\eqref{3.71}, Lebesgue's dominated convergence theorem, and Proposition \ref{prop5.2},
we have
\begin{align}\label{3.72}
&\lim_{\sigma\rightarrow0}\int_{\mathbb{R}^{n+1}_+} \M^{\v} \cdot \big\{\Delta \Psi_{\sigma}+\nabla\mbox{div}\,\Psi_{\sigma} \big\}\,
\dd\mathbf{x}\dd t\\
&=\int_{\mathbb{R}^{n+1}_+} \M^{\v} \cdot \big\{\Delta \bp+\nabla\mbox{div}\,\bp\big\}\,\dd\mathbf{x}\dd t.\nonumber
\end{align}
Substituting \eqref{3.63-5}--\eqref{3.61}, \eqref{3.68}--\eqref{3.70}, and \eqref{3.72}
into \eqref{3.59-1} leads to \eqref{5.48}.
\end{proof}

\begin{lemma}\label{lem5.10}
Let $\xi(\mathbf{x})\in C_0^1(\mathbb{R}^n)$ be any smooth function with compact support.
Then
\begin{align} \label{3.74}
	\int_{\mathbb{R}^n} \nabla\Phi^\v(t,\mathbf{x}) \cdot \nabla \xi(\mathbf{x})\, \dd\mathbf{x}
	=-\kappa\int_{\mathbb{R}^n} \rho^{\v}(t,\mathbf{x}) \xi(\mathbf{x})\, \dd\mathbf{x}\qquad \mbox{for $t\geq0$}.
\end{align}
\end{lemma}

The proof is direct by using $\eqref{5.5-2}_3$, \eqref{5.3}, and \eqref{5.43}, so we omit the details here.

\subsection{$H^{-1}_{\rm loc}$-Compactness}

To use the compensated compactness framework in \cite{Perepelitsa},
we need the $H_{\rm loc}^{-1}$-compactness of entropy dissipation measures.

\begin{lemma}[$H_{\rm loc}^{-1}$--compactness]\label{lem7.1}
Let $(\eta, q)$ be a weak entropy pair defined in \eqref{weakentropy}
for any smooth compact supported function $\psi(s)$ on $\mathbb{R}$.
Then, for $\v\in (0, \v_0]$,
\begin{align}\label{7.1}
	\partial_t\eta(\rho^\v,m^\v)+\partial_rq(\rho^\v,m^\v) \qquad   \mbox{is compact in $ H^{-1}_{\rm loc}(\mathbb{R}^2_+)$}.
\end{align}
\end{lemma}

\begin{proof} To obtain \eqref{7.1},
we have to make the argument in the weak sense, since $(\rho^\v, \mathcal{M}^{\v}, \Phi^\v)$ is a weak solution
of CNSPEs \eqref{1.1}.
In fact, we first have to study the equation for $\partial_t\eta(\rho^\v,m^\v)+\partial_rq(\rho^\v,m^\v) $ in
the distributional sense, which is more complicated than that in \cite{Perepelitsa,Chen6}.
We divide the proof into five steps.

1. Since
\begin{align*}
\eta(\rho,\rho u)&=\rho \int_{-1}^1 \psi(u+\rho^{\theta}s)[1-s^2]_+^{\fb}\,\dd s,\\
q(\rho,\rho u)&=\rho \int_{-1}^1 (u+\theta\rho^{\theta}s)\psi(u+\rho^{\theta}s)[1-s^2]_+^{\fb}\,\dd s,
\end{align*}
then it follows from  \cite[Lemma 2.1]{Perepelitsa} that
\begin{align}
&|\eta(\rho,\rho u)|+|q(\rho,\rho u)|\leq C_{\psi} \rho \qquad \mbox{for $\gamma\in(1,3]$},\label{7.2}\\
&|\eta(\rho,\rho u)|\leq C_{\psi} \rho, \quad   |q(\rho,\rho u)|\leq C_{\psi} \big(\rho+\rho^{1+\theta}\big)
\quad   \mbox{for $\gamma\in (3,\infty)$},\label{7.2-1}\\
&|\partial_\rho\eta (\rho,\rho u)|\leq C_{\psi} \big(1+\rho^{\theta}\big),\quad   |\partial_m \eta(\rho,\rho u)|\leq C_{\psi}.\label{7.6}
\end{align}
On the other hand, if we regard $\partial_m \eta(\rho,m)$ as a function of $(\rho, u)$, then
\begin{align}\label{7.3}
|\partial_{m\rho}\eta|\leq C_{\psi} \rho^{\theta-1}, \qquad |\partial_{mu}\eta|\leq C_{\psi}.
\end{align}

2. Denote $(\eta^{\v,b}, q^{\v,b}):=(\eta, q)(\rho^{\v,b}, m^{\v,b})$
and $(\eta^{\v}, q^{\v}):=(\eta, q)(\rho^{\v}, m^{\v})$ for simplicity. Multiplying $\eqref{2.1}_1$ by $\eta_\rho(\rho^{\v,b}, m^{\v,b})$,  $\eqref{2.1}_2$ by $\eta_m(\rho^{\v,b}, m^{\v,b})$, and add them together to obtain
\begin{align}\label{7.4}
&\partial_t \eta^{\v,b}+\partial_r q^{\v,b}\\
& =-\frac{n-1}{r} m^{\v,b} \big(\eta_\rho^{\v,b}+u^{\v,b} \eta_m^{\v,b}\big)
-\kappa \eta_m^{\v,b} \frac{\rho^{\v,b}}{r^{n-1}} \int_0^r \rho^{\v,b}(t,z)\,z^{n-1} \dd z\nonumber\\
&\quad + \v \eta^{\v,b}_m\Big\{ \big(\rho^{\v,b} (u_r^{\v,b}+\frac{n-1}{r} u^{\v,b})\big)_r
-\frac{n-1}{r} u^{\v,b}\rho^{\v,b}_r \Big\},\nonumber
\end{align}
where $\rho^{\v,b}$ is understood to be zero in domain $[0,T]\times [0,a)$
so that $\int_a^r \rho^{\v,b}(t,z)\,z^{n-1} \dd z$ can be rewritten
as $\int_0^r \rho^{\v,b}(t,z)\,z^{n-1}\dd z$ in the potential term.

Let $\phi(t,r)\in C_0^\infty(\mathbb{R}^2_+)$ and $b\gg1$ so that ${\rm supp}\,\phi(t,\cdot)\Subset(a,b(t))$.
Then multiplying \eqref{7.4}  by $\phi$ and integrating by parts yield that
\begin{align}\label{7.7}
&\int_{\mathbb{R}^2_+} \big(\partial_t \eta^{\v,b}+\partial_r q^{\v,b}\big) \phi\, \dd r\dd t\\
&=-\int_{\mathbb{R}^2_+} \frac{n-1}{r} m^{\v,b} \big(\eta_\rho^{\v,b}+u^{\v,b} \eta_m^{\v,b}\big) \phi\, \dd r\dd t\nonumber\\
&\quad -\v \int_{\mathbb{R}^2_+}\rho^{\v,b} (\eta_m^{\v,b})_r \, \Big(u_r^{\v,b}+\frac{n-1}{r} u^{\v,b}\Big) \phi\, \dd r\dd t\nonumber\\
&\quad -\v \int_{\mathbb{R}^2_+}  \rho^{\v,b} \eta_m^{\v,b} \, \Big(u_r^{\v,b}+\frac{n-1}{r} u^{\v,b}\Big) \phi_r\,\dd r\dd t\nonumber\\
&\quad-\v \int_{\mathbb{R}^2_+} \eta_m^{\v,b} \frac{n-1}{r}  u^{\v,b}\rho^{\v,b}_r  \phi\, \dd r\dd t\nonumber\\
&\quad -\kappa \int_{\mathbb{R}^2_+} \eta_m^{\v,b}  \frac{\rho^{\v,b}}{r^{n-1}} \Big(\int_0^r \rho^{\v,b}(t,z)\, z^{n-1} \dd z\Big)\phi\,\dd r\dd t\nonumber\\
&:= \sum_{i=1}^5I_{i}^{\v,b}.\nonumber
\end{align}

3. It is direct to see that
\begin{align}\label{7.8}
\eta^{\v,b}\longrightarrow \eta^\v\qquad  \mbox{{\it a.e.} in  $\{(t,r)\,:\, \rho^{\v}(t,r)\neq0 \}\,$ as $b\rightarrow\infty$}.
\end{align}
In $\{(t,r)\,:\, \rho^{\v}(t,r)=0 \}$, it follows from \eqref{7.2}--\eqref{7.2-1} that
\begin{align}\label{7.9}
|\eta^{\v,b}|\leq C_{\psi} \rho^{\v,b}\longrightarrow0=\eta^{\v}\qquad\, \mbox{as $b\rightarrow\infty$}.
\end{align}
Combining  \eqref{7.8}--\eqref{7.9} together, we have
\begin{align}\label{7.10}
\eta^{\v,b}\longrightarrow \eta^\v \qquad  \mbox{{\it a.e.} as $b\rightarrow\infty$}.
\end{align}
Similarly, we also have
\begin{align}\label{7.11}
q^{\v,b}\longrightarrow q^{\v} \qquad  \mbox{{\it a.e.}  as $b\rightarrow\infty$}.
\end{align}
Let $K\Subset(0,\infty)$ be any compact subset.  For $\gamma\in(1,3]$,
it follows from \eqref{2.65} and \eqref{7.2} that
\begin{align}\label{7.12}
\int_0^T\int_{K} \big(|\eta^{\v,b}|+|q^{\v,b}|\big)^{\gamma+1} \dd r\dd t
&\leq C_{\psi}\int_0^T\int_{K} |\rho^{\v,b}|^{\gamma+1}\, \dd r\dd t\\
&\leq C_{\psi}(K,M,E_0,T). \nonumber
\end{align}
For $\gamma\in(3,\infty)$, it follows from \eqref{7.2-1} and \eqref{2.95} that
\begin{align}\label{7.13}
	\int_0^T\int_{K} \big(|\eta^{\v,b}|+|q^{\v,b}|\big)^{\frac{\gamma+\theta}{1+\theta}} \dd r\dd t
	&\leq C_{\psi}\int_0^T\int_{K} \big(|\rho^{\v,b}|^{\frac{\gamma+\theta}{1+\theta}}+|\rho^{\v,b}|^{\gamma+\theta}\big)\,\dd r\dd t\\
	&\leq C_{\psi}(K,M,E_0,T). \nonumber
\end{align}
We take  $p_1=\gamma+1>2$ when $\gamma\in(1,3]$, and $p_1=\frac{\gamma+\theta}{1+\theta}>2$ when $\gamma\in (3,\infty)$.
Then it follows from \eqref{7.12}--\eqref{7.13}  that
\begin{align}\nonumber
(\eta^{\v,b}, \  q^{\v,b})  \qquad   \mbox{is uniformly bounded in } L^{p_1}_{\rm loc}(\mathbb{R}^2_+),
\end{align}
which, together with \eqref{7.10} and \eqref{7.11}, yields that, up to a subsequence,
\begin{align}\nonumber
(\eta^{\v,b}, \  q^{\v,b}) \to (\eta^{\v}, \  q^{\v})
\qquad \mbox{in $L^2_{\rm loc}(\mathbb{R}^2_+)$  as $b\rightarrow\infty$}.
\end{align}
Thus, for any $\phi\in C^1_0(\mathbb{R}^2_+)$, we see that,
as $b\rightarrow \infty$ (up to a subsequence),
\begin{align}\label{7.16}
	\int_{\mathbb{R}^2_+} \big(\partial_t \eta^{\v,b}+\partial_r q^{\v,b}\big) \phi\, \dd r\dd t
	&=-\int_{\mathbb{R}^2_+} \big(\eta^{\v,b}\partial_t \phi+ q^{\v,b} \partial_r \phi\big)\, \dd r\dd t
	\\
	&\longrightarrow-\int_{\mathbb{R}^2_+} \big(\eta^{\v}\partial_t \phi+ q^{\v} \partial_r \phi\big)\, \dd r\dd t.\nonumber
\end{align}
Furthermore,
$(\eta^{\v}, q^{\v})$  is  uniformly bounded in $L^{p_1}_{\rm loc}(\mathbb{R}^2_+)$ for some $p_1>2$,
which implies that
\begin{align}\label{7.17-1}
\partial_t\eta^{\v}+\partial_rq^{\v} \qquad  \mbox{is uniformly bounded in $\v>0$ in $W^{-1,p_1}_{\rm loc}(\mathbb{R}_+^2)$}.
\end{align}

4. Now we estimate the terms of \eqref{7.7}-RHS.
For $I_1^{\v,b}$, a direct calculation shows that
$$
|\eta_\rho+u\eta_m|\leq C_{\psi} (1+\rho^{\theta})
$$
which,
together with Lemma \ref{lem5.6} and  similar arguments in \eqref{7.8}--\eqref{7.10},  yields that
\begin{align}\label{7.18}
\frac{n-1}{r} m^{\v,b} \big(\eta_\rho^{\v,b}+u^{\v,b} \eta_m^{\v,b}\big)
\longrightarrow \frac{n-1}{r} m^{\v} \big(\eta_\rho^\v+u^{\v} \eta_m^\v\big)\quad\mbox{{\it a.e.} as  $b\rightarrow\infty$}.
\end{align}
Then it follows from
\eqref{5.10-1}  that
\begin{align}\label{7.19}
&\int_0^T\int_{K} \Big|\frac{n-1}{r} m^{\v,b}\big(\eta_\rho^{\v,b}+u^{\v,b} \eta_m^{\v,b}\big)\Big|^{\frac76}\dd r\dd t\\
&\leq C\int_0^T\int_{K}  \big(\rho^{\v,b}|u^{\v,b}|^2 +\rho^{\v,b} +(\rho^{\v,b})^{\gamma}\big)^{\frac76}\,  \dd r\dd t\nonumber\\
&\leq
\begin{cases}
	\displaystyle C\Big(1+\int_0^T\int_{K} \rho^{\v,b} |u^{\v,b}|^3\, \dd r\dd t \Big)^{\frac79}
	\Big(\int_0^T\int_{K}\big(1+|\rho^{\v,b}|^{\gamma+1}\big)\dd r\dd t \Big)^{\frac29}\,\,\, \mbox{for $\gamma\in(1,3]$},\\[3mm]
	\displaystyle C\Big(1+\int_0^T\int_{K}\rho^{\v,b}|u^{\v,b}|^3\, \dd r\dd t \Big)^{\frac79}
	\Big(\int_0^T\int_{K}\big(1+|\rho^{\v,b}|^{\gamma+\theta}\big)\dd r\dd t \Big)^{\frac29}\,\,\, \mbox{for $\gamma\in(3,\infty)$}
\end{cases}
\nonumber\\[2mm]
&\leq C(K,M, E_0, T).\nonumber
\end{align}
Using \eqref{7.18}--\eqref{7.19},  up to a subsequence, we have
\begin{align}
&I^{\v,b}_1 \longrightarrow -\int_{\mathbb{R}^2_+} \frac{n-1}{r} m^{\v} \big(\eta_\rho^{\v}+u^{\v} \eta_m^{\v}\big) \phi\, \dd r\dd t
\qquad  \mbox{as $b\rightarrow\infty$},\label{7.20}\\
&\int_0^T\int_{K}\Big| \frac{n-1}{r} m^{\v} \big(\eta_\rho^{\v}+u^{\v} \eta_m^{\v}\big) \Big|^{\frac76}\,\dd r\dd t \leq C(K,M,E_0,T).\label{7.20-1}
\end{align}

For $I_2^{\v,b}, I_4^{\v,b}$, and $I_5^{\v,b}$, it follows from Lemmas \ref{lem2.1} and \ref{lem2.2}, and \eqref{7.6}--\eqref{7.3} that
\begin{align*}
&\int_0^T\int_{K}  \Big|\v \rho^{\v,b} (\eta_m^{\v,b})_r \, (u_r^{\v,b}+\frac{n-1}{r} u^{\v,b}) \Big|\,\dd r\dd t\\
&\quad\leq C_{\psi}(K) \int_0^T\int_{K}\Big(\v \rho^{\v,b}|u_r^{\v,b}|^2+\v(\rho^{\v,b})^{\gamma-2}|\rho_r^{\v,b}|^2+\rho^{\v,b}|u^{\v,b}|^2\Big)
\,\dd r\dd t\nonumber\\[1mm]
&\quad\leq C_{\psi}(K,M,E_0,T), \\[2mm]
&\int_0^T\int_{K}  \Big|\v \eta_m^{\v,b}\frac{n-1}{r} \rho^{\v,b}_r u^{\v,b}\Big|\, \dd r\dd t\\
&\quad\leq C_{\psi}(K) \Big(\v^2\int_0^T\int_{K} |(\sqrt{\rho^{\v,b}})_r|^2\, \dd r\dd t\Big)^{\frac12}\Big(\int_0^T\int_{K} \rho^{\v,b}|u^{\v,b}|^2\, \dd r\dd t\Big)^{\frac12}\\[1mm]
&\quad\leq C_{\psi}(K,M,E_0,T),\\[2mm]
&\int_0^T\int_{K}  \Big|\kappa \eta^{\v,b}_m \frac{\rho^{\v,b} }{r^{n-1}} \int_0^r \rho^{\v,b}(t,z)\,z^{n-1} \dd z \Big| \dd r\dd t
\leq C_{\psi}(K,M) \int_0^T\int_{K} \rho^{\v,b}\,\dd r\dd t\nonumber\\[1mm]
&\hspace{6.6cm}\leq C_{\psi}(K,M,E_0,T).
\end{align*}
Thus, there exist  local bounded Radon measures $\mu_1^\v, \mu_2^\v$, and $ \mu_3^\v$ on $\mathbb{R}^2_+$ so that,
as $b\to \infty$ (up to a subsequence),
\begin{align}
&-\v \rho^{\v,b} (\eta_m^{\v,b})_r \, (u_r^{\v,b}+\frac{n-1}{r} u^{\v,b}) \longrightharpoonup \mu_1^{\v},
\nonumber\\[2mm]
&-\v \eta_m^{\v,b} \frac{n-1}{r} \rho^{\v,b}_r u^{\v,b} \longrightharpoonup \mu_2^{\v},\nonumber\\[2mm]
&-\kappa \eta_m^{\v,b}  \frac{\rho^{\v,b} }{r^{n-1}} \int_0^r \rho^{\v,b}(t,z)\, z^{n-1} \dd z \longrightharpoonup \mu_3^{\v}.\nonumber
\end{align}
In addition,
\begin{equation}\label{7.21-2}
\mu_i^\v((0,T)\times \mathcal{O})\leq C_{\psi}(K,T,E_0)\qquad\,\,   \mbox{for $i=1,2,3$},
\end{equation}
for each open subset $\mathcal{O}\subset K$.
Then, up to a subsequence, we have
\begin{align}\label{7.23}
I_{2}^{\v,b}+I_{4}^{\v,b}+I_{5}^{\v,b}\longrightarrow \langle \mu_1^{\v}+\mu_2^{\v}+\mu_3^{\v}, \phi\rangle
\qquad \mbox{as $b\rightarrow\infty$}.
\end{align}

For $I_3^{\v,b}$, we notice from Lemma \ref{lem2.1} that
\begin{align}
& \int_0^T\int_{K} \Big|\sqrt{\v}\rho^{\v,b} \eta_m^{\v,b} \,(u_r^{\v,b}+\frac{n-1}{r} u^{\v,b}) \Big|^{\frac43}\dd r\dd t \nonumber\\
&\leq C_{\psi}(K)\int_0^T\int_{K} \Big|\sqrt{\v}\rho^{\v,b}  (|u_r^{\v,b}|+|u^{\v,b}|) \Big|^{\frac43}\dd r\dd t\nonumber\\
&\leq C_{\psi}(K) \Big(\v\int_0^T\int_{K}\big(\rho^{\v,b}|u_r^{\v,b}|^2+\rho^{\v,b}|u^{\v,b}|^2\big)\,\dd r\dd t\Big)^{\frac23} \Big(\int_0^T\int_{K} |\rho^{\v,b}|^2\dd r\dd t\Big)^{\frac13}\nonumber\\[1mm]
&\leq C_{\psi}(K,M,E_0,T).\nonumber
\end{align}
Then there exists a function $f^\v$ such that, as $b\to \infty$ (up to a subsequence),
\begin{align}
&-\sqrt{\v}\rho^{\v,b} \eta_m^{\v,b} \, \big(u_r^{\v,b}+\frac{n-1}{r} u^{\v,b}\big)\, {\longrightharpoonup}\, f^\v \qquad\,
\mbox{weakly in $L^{\frac43}_{\rm loc}(\mathbb{R}^2_+)$},\label{7.25}\\
&\int_0^T\int_{K} |f^{\v}|^{\frac43}\,\dd r\dd t\leq C_{\psi}(K,M,E_0,T).\label{7.26}
\end{align}
Thus, it follows from \eqref{7.25} that
\begin{align}\label{7.27}
I_3^{\v,b}\longrightarrow \sqrt{\v}\int_0^T\int_{K} f^{\v} \phi_r \,\dd r\dd t
\qquad \mbox{as $b\to \infty$ (up to a subsequence)}.
\end{align}

\smallskip
5. Taking $b\rightarrow\infty$ (up to a subsequence) on both sides of \eqref{7.7}, then it follows
from \eqref{7.16}, \eqref{7.20}, \eqref{7.23}, and \eqref{7.27} that
\begin{align*}
\partial_t \eta^{\v}+\partial_r q^{\v}
=- \frac{n-1}{r} \rho^{\v} u^{\v}\big(\eta_\rho^{\v}+u^{\v} \eta_m^{\v}\big)
+\mu_1^{\v}+\mu_2^{\v}+\mu_3^{\v}-\sqrt{\v} f^{\v}_r
\end{align*}
in the sense of distributions. Noting \eqref{7.20-1}--\eqref{7.21-2}, we know that
\begin{align}\label{7.31}
- \frac{n-1}{r} \rho^{\v} u^{\v}\big(\eta_\rho^{\v}+u^{\v} \eta_m^{\v}\big) +\mu_1^{\v}+\mu_2^{\v}+\mu_3^{\v}
\,\,\,  \mbox{is a local bounded Radon measure},
\end{align}
and the bound is uniform in $\v>0$.
From \eqref{7.26}, we know that
\begin{align}\label{7.32}
\sqrt{\v} f^{\v}_r \longrightarrow 0 \qquad  \mbox{in $W^{-1,\frac43}_{\rm loc}(\mathbb{R}^2_+)$ as $\v\rightarrow0+$}.
\end{align}
Then it follows from \eqref{7.31}--\eqref{7.32} that
\begin{align}\label{7.33}
\partial_t \eta^{\v}+\partial_r q^{\v} \,\,\, \mbox{are confined in a compact subset of $W^{-1, p_2}_{\rm loc}(\mathbb{R}^2_+)$}
\end{align}
for some $p_2\in(1,2)$.

The interpolation compactness theorem ({\it cf}. \cite{Chen2,Chen5}) indicates that,
for $p_2>1$, $p_1\in(p_2,\infty]$, and $p_0\in[p_2, p_1)$,
\begin{align*}
&(\mbox{compact set of} \  W^{-1,p_2}_{\rm loc}(\mathbb{R}_+^2)) \cap (\mbox{bounded  set of} \  W^{-1,p_1}_{\rm loc}(\mathbb{R}_+^2))\\
&\subset (\mbox{compact set of} \  W^{-1,p_0}_{\rm loc}(\mathbb{R}_+^2)),\nonumber
\end{align*}
which is a generalization of Murat's lemma in \cite{F. Murat,L. Tartar}.
Combining this
theorem for $1<p_2<2$, $p_1>2$, and $p_0=2$ with the facts in \eqref{7.17-1} and \eqref{7.33}, we conclude \eqref{7.1}.
\end{proof}

\begin{remark}
	Since  $(\rho^\v, m^{\v})$ are the weak solutions of CNSPEs, it is not convenient to use the weak formulation to
	prove the $H_{\rm loc}^{-1}$-compactness directly.
	Therefore, in this section above, we first study the equation satisfied by $\partial_t\eta^{\v,b}+\partial_rq^{\v,b}$,
	then take limit $b\rightarrow \infty$ to obtain the equation satisfied by $\partial_t\eta^{\v}+\partial_rq^{\v}$ in the distributional sense,
	and finally use the equation to establish the $H_{\rm loc}^{-1}$--compactness.
\end{remark}

Combining Proposition \ref{prop5.2} with   Lemmas \ref{lem5.1} and \ref{lem4.8-1}--\ref{lem7.1}, we have

\begin{theorem}\label{thm7.2}
Let $(\rho_0^\v, m_0^\v)$ be the initial data satisfying \eqref{1.13-2}--\eqref{1.13-4}.
Then, for each $\v>0$, there exists a global spherically symmetric weak
solution
$$
(\rho^\v, \M^\v, \Phi^\v)(t,\mathbf{x}):=(\rho^\v(t,r), m^\v(t,r)\frac{\mathbf{x}}{r}, \Phi^\v(t,r))
$$
of CNSPEs \eqref{1.1} in the sense of Definition {\rm \ref{definition-NSP}}.
Moreover, $m^\v(t,r)=\rho^\v(t,r) u^\v(t,r)$, with
$u^\v(t,r):=\frac{m^\v(t,r)}{\rho^\v(t,r)}\,$\, a.e. on  $\{(t,r)\,:\,\rho^\v(t,r)\ne 0\}$
and $u^\v(t,r):=0\,$ \, a.e. on $\{(t,r)\,:\, \rho^\v(t,r)=0\,\, \mbox{or $\,r=0$}\}$,
satisfies the following properties{\rm :}
\begin{align}
\begin{split}
  &\hspace{0.5cm}	\rho^{\v}(t,r)\geq0 \quad a.e.,
\end{split}\\[2mm]
  &\hspace{0.5cm}u^\v(t,r)=0, \,\,\big(\frac{m^\v}{\sqrt{\rho^\v}}\big)(t,r)=\sqrt{\rho^\v}(t,r)u^\v(t,r)=0\label{8.1}\\
  &\hspace{6cm} a.e.\,\, \mbox{on $\{(t,r)\,:\, \rho^\v(t,r)=0\}$}, \nonumber\\[2mm]
\begin{split}
  &\hspace{0.5cm}\Phi_r^\v(t,r)=\frac{\kappa}{r^{n-1}} \int_0^r  \rho^\v(t,z)\, z^{n-1}\dd z \qquad\,\, \mbox{for}\,\, (t,r)\in \R_+^2,\label{8.6}\\
\end{split}\\
\begin{split}
  &\hspace{0.5cm}\int_{0}^{\infty} \rho^\v(t,r)\, r^{n-1} \dd r= \frac{M}{\omega_n} \qquad\,\, \mbox{for all $t\geq0$},\label{8.2}
\end{split}\\[2mm]
&\int_0^\infty  \big(\frac12\rho^{\v} |u^{\v}|^2+(\rho^{\v})^{\gamma}\big)(t,r)\,r^{n-1} \dd r
+\v\int_0^t\int_0^\infty  (\rho^{\v} |u^{\v}|^2)(s,r)\, r^{n-3} \dd r \dd s \\
&\quad +\int_0^{\infty} |\Phi_r^\v(t,r)|^2\,r^{n-1}\dd r
+\int_0^{\infty}  \Big(\int_0^r \rho^\v(t,z)\, z^{n-1} \dd z\Big)\rho^\v(t,r)\,r \dd r\nonumber\\[1mm]
&\quad\leq C(M,E_0)\qquad\,\, \mbox{for all $t \geq0$},\nonumber\\[2mm]
	&\v^2 \int_0^\infty r^{n-1}\big|\big(\sqrt{\rho^{\v}(t,r)}\big)_r\big|^2\dd r
	+\v\int_0^T\int_{0}^\infty r^{n-1}\big|\big((\rho^{\v}(s,r))^{\frac{\gamma}{2}}\big)_r\big|^2\dd r\dd s \label{8.3}\\[1mm]
	&\quad \leq C(M,E_0,T) \qquad\,\,\mbox{for $t\in[0,T]$},\nonumber\\[2mm]
	&\int_0^T\int_d^D  \big(\rho^{\v}|u^{\v}|^3+(\rho^{\v})^{\gamma+\theta}+(\rho^{\v})^{\g+1}\big)(t,r)\, r^{n-1} \dd r\dd t\label{8.5}\\[1mm]
    &\quad \leq C(d,D,M,E_0,T),\nonumber		
\end{align}
for any fixed $T>0$ and any compact subset $[d,D]\Subset (0,\infty)$. Moreover, the following energy inequality holds{\rm :}
\begin{align}\label{8.7}
&\int_{\R^n}\Big(\frac{1}{2}\Big|\frac{\M^\v}{\sqrt{\rho^\v}}\Big|^2
+ \rho^\v e(\rho^\v)-\frac{\kappa}{2}  |\nabla\Phi^\v|^2\Big)(t,\mathbf{x})\,\dd\mathbf{x}\\
&\leq\int_{\R^n}\Big(\frac{1}{2}\Big|\frac{\M_0^\v}{\sqrt{\rho_0^\v}}\Big|^2
+ \rho^\v_0 e(\rho_0^\v)-\frac{\kappa}{2}   |\nabla\Phi_0^\v|^2\Big)(\mathbf{x})\, \dd \mathbf{x}\qquad\mbox{for $t\geq0$}.
  \nonumber	
\end{align}
Furthermore, let $(\eta^\psi, q^\psi)$ be an entropy pair defined in \eqref{weakentropy}
for a smooth function $\psi(s)$ with compact support on $\mathbb{R}$.
Then, for $\v\in (0,\v_0]$,
\begin{align*}
	\partial_t\eta^\psi(\rho^\v,m^\v)+\partial_rq^\psi(\rho^\v,m^\v) \qquad  \mbox{is compact in $H^{-1}_{\rm loc}(\mathbb{R}^2_+)$}.\nonumber
\end{align*}
\end{theorem}

\section{Proof of the Main Theorems}

In this section, we give a complete proof of Main Theorem II: Theorem \ref{thm1.2},
which leads to Main Theorem I: Theorem \ref{thm:2.1}, as indicated in Remark \ref{remark:2.5a}.
We divide the proof into several steps.

1. The uniform estimates and compactness properties obtained in Theorem \ref{thm7.2} imply that
the weak solutions
$$
\dis (\rho^{\v}, \M^{\v}, \Phi^\v)=(\rho^{\v}, m^{\v}\frac{\mathbf{x}}{r}, \Phi^\v)
$$
of CNSPEs \eqref{1.1}
satisfy the compensated compactness framework in
\cite{Perepelitsa}.
Then the compactness theorem in \cite{Perepelitsa} for the whole range $\gamma>1$
implies that there exists a vector function $(\rho,m)(t,r)$ such that
\begin{align}\label{7.51-1}
(\rho^\v,m^{\v})\longrightarrow (\rho,m) \quad  \mbox{{\it a.e.} $\,(t,r)\in \mathbb{R}^2_+\,$ as $\v\rightarrow0+$ (up to a subsequence)}.
\end{align}

By similar arguments as in the proof of Lemma \ref{lem5.6}, we find that
$m(t,r)=0$ {\it a.e.} on $\{(t,r) \,:\, \rho(t,r)=0\}$.
We can define the  limit velocity $u(t,r)$ by setting $u(t,r):=\frac{m(t,r)}{\rho(t,r)}$ {\it a.e.} on $\{(t,r)\,:\,\rho(t,r)\neq0\}$
and $u(t,r):=0$ {\it a.e.} on $\{(t,r)\,: \rho(t,r)=0\,\,\, \mbox{or $r=0$}\}$.
Then we have
\begin{align}
m(t,r)=\rho(t,r) u(t,r).\nonumber
\end{align}
We can also define $(\frac{m}{\sqrt{\rho}})(t,r):=\sqrt{\rho(t,r)} u(t,r)$, which is $0$ {\it a.e.} on $\{(t,r) \ :\ \rho(t,r)=0\}$.
Moreover, we obtain that, as $\v\rightarrow0+$,
\begin{equation*}
\frac{m^{\v}}{\sqrt{\rho^{\v}}}\equiv\sqrt{\rho^{\v}} u^{\v}\longrightarrow \frac{m}{\sqrt{\rho}}\equiv\sqrt{\rho} u \quad
\mbox{strongly in $L^2([0,T]\times[d,D], r^{n-1}\dd r \dd t)$}
\end{equation*}
for any given $T$ and $[d,D]\Subset (0,\infty)$.

Notice that $|m|^{\frac{3(\gamma+1)}{\gamma+3}}\leq C\big(\rho|u|^3+\rho^{\gamma+1} \big)$
which,  together with
\eqref{8.5}, yields that
\begin{align}\label{7.51}
(\rho^\v,m^{\v})\longrightarrow (\rho,m) \qquad \mbox{in $L^{p_1}_{\rm loc}(\mathbb{R}^2_+)\times L^{p_2}_{\rm loc}(\mathbb{R}^2_+)$}
\end{align}
for $p_1\in[1,\gamma+1)$ and $p_2\in[1,\frac{3(\gamma+1)}{\gamma+3})$,
where $L^{p_j}_{\rm loc}(\mathbb{R}^2_+)$ represents $L^{p_j}([0,T]\times K)$ for any $T>0$ and $K\Subset (0,\infty)$, $j=1,2$.

From the same estimates, we also obtain the convergence of the mechanical energy as $\v\rightarrow0+:$
\begin{align}\nonumber
\eta^{\ast}(\rho^\v,m^\v)\longrightarrow \eta^{\ast}(\rho,m)\qquad \mbox{in $L^1_{\rm loc}(\mathbb{R}_+^2)$}.
\end{align}
Since $\eta^{\ast}(\rho,m)$ is a convex function, by passing limit in \eqref{5.7-1} and \eqref{3.36}, we have
\begin{align}\label{7.82}
\int_{t_1}^{t_2} \int_0^\infty \big\{\eta^{\ast}(\rho,m)(t,r) +\rho(t,r)\big\}\,r^{n-1}\dd r\dd t\leq C(M,E_0)(t_2-t_1),
\end{align}
which indicates that
\begin{align}\label{7.83}
	\sup_{0\le t\le T}\int_0^\infty \big\{\eta^{\ast}(\rho,m)(t,r)+\rho(t,r)\big\}\,r^{n-1} \dd r\leq C(M,E_0).
\end{align}
That is, $\rho(t,r)\in L^\infty([0, T]; L^\gamma(\mathbb{R}; r^{n-1}\dd r))$, since $\eta^{\ast}(\rho,m)$ contains a term: $\rho^\gamma$.
This indicates that $\rho(t,\mathbf{x})$ is a function in  $L^\infty([0, T]; L^\gamma(\mathbb{R}^n))$ for $\gamma>1$ (rather than a measure in space-time),
so that no delta measure ({\it i.e.}, concentration) is formed in the density $\rho$ in the time interval $[0,T]$, especially
at the origin $r=0$.

2. A direct calculation shows that
\begin{align}
&\lim_{\v\rightarrow0+}\int_0^T\int_0^D\left|\int_0^r  \rho^{\v}(t,z)\, z^{n-1}\dd z-\int_0^r \rho(t,z)\,z^{n-1}\dd z\right|\dd r \dd t\nonumber\\
&\leq D\lim_{\v\rightarrow0+} \int_0^T\int_\sigma^D  |\rho^{\v}(t,z) - \rho(t,z)|\, z^{n-1}\dd z \dd t
+C(D,T,M,E_0)\sigma^{n(1-\frac1\gamma)}\nonumber\\
&\leq C(D,T,M,E_0)\sigma^{n(1-\frac1\gamma)},\nonumber
\end{align}
which,  as $\v\rightarrow 0+$ (up to a subsequence),  yields that
\begin{align}\label{5.1}
\Phi_r^{\v}(t,r)r^{n-1} & =\kappa\int_0^r  \rho^{\v}(t,z)\, z^{n-1}\dd z\\
&\longrightarrow \kappa\int_0^r  \rho(t,z)\,z^{n-1}\dd z
\quad\mbox{\it a.e. $(t,r)\in\mathbb{R}_+^2$}. \nonumber
\end{align}
Then, using Fatou's lemma, \eqref{5.7-1}, \eqref{7.51-1}, \eqref{5.1}, and similar arguments as in \eqref{7.82}--\eqref{7.83},
we have
\begin{align}\label{7.100}
\int_0^{\infty} \Big(\int_0^r\rho(t,z)\, z^{n-1}\dd z\Big)\rho(t,r)\, r\dd r\leq C(M,E_0)
\quad \mbox{for {\it a.e.} $t\geq0$}.
\end{align}

Estimate \eqref{5.7-0} implies that there exists a function $\Phi(t,\mathbf{x})=\Phi(t,r)$ such that, as $\v\rightarrow0+$ (up to a subsequence),
\begin{align}
&\Phi^{\v} \longrightharpoonup \Phi \,\,\,
\mbox{{weak-}$\ast$ in $L^\infty(0,T; H^1_{\rm loc}(\mathbb{R}^n))$ and weakly in $L^2(0,T; H^1_{\rm loc}(\mathbb{R}^n))$},\label{7.102}\\[2mm]
&\|\Phi(t)\|_{L^{\frac{2n}{n-2}}(\mathbb{R}^n)}+\|\nabla \Phi(t)\|_{L^2(\mathbb{R}^n)}\leq C(M,E_0)
\qquad\mbox{{\it a.e.} $t\geq0$}.\label{7.103}
\end{align}
It follows from the uniqueness of  limit and \eqref{5.1} that
\begin{align}\label{7.101}
\Phi_r(t,r)r^{n-1} =\kappa\int_0^r \rho(t,z)\, z^{n-1}\dd z\qquad\mbox{{\it a.e.} $(t,r)\in\mathbb{R}_+^2$}.
\end{align}

To prove the energy inequality for the case that $\kappa=1$ (gaseous stars), we need the strong convergence of the potential functions.
Notice that
\begin{align}\nonumber
\frac{1}{r^{n-1}} \Big|\int_0^{r} (\rho^{\v}-\rho)(t,z)\, z^{n-1} \dd z\Big|^2\leq C(M) r^{-n+1}\qquad \mbox{for $r>0$ and {\it a.e.} $t\geq0$},
\end{align}
which yields that
\begin{align}\label{8.8}
\int_k^\infty  \frac{1}{r^{n-1}} \Big|\int_0^{r} (\rho^{\v}-\rho)(t,z)\,z^{n-1}\dd z\Big|^2 \dd r
\leq C(M) k^{-n+2} \qquad \mbox{{\it a.e.} $t\geq0$}.
\end{align}
Using the H\"{o}lder inequality,
\begin{align}\label{8.9}
&\frac{1}{r^{n-1}} \Big|\int_0^{r} (\rho^{\v}-\rho)(t,z)\, z^{n-1} \dd z\Big|^2\\
& \leq Cr^{-n+1} \left(\int_0^r \big( (\rho^{\v})^\g+\rho^\g\big)\,z^{n-1}\dd z\right)^{\f2\g}\, r^{2n(1-\f1\g)}\nonumber\\
&\leq C(E_0,M)\,r^{n+1-\frac{2n}{\g}}\qquad  \mbox{{\it a.e.} $t\geq0$}.\nonumber
\end{align}
Since $n+1-\frac{2n}{\g}>-1$ for $\g>\frac{2n}{n+2}$, then it follows from \eqref{7.51-1}, \eqref{8.9},
and Lebesgue's dominated convergence theorem that, for any given $k>0$,
\begin{align}\nonumber
\int_0^T\int_0^k \frac{1}{r^{n-1}} \Big|\int_0^{r} (\rho^{\v}-\rho)(t,z)\, z^{n-1} \dd z\Big|^2 \dd r\dd t\to 0\qquad \mbox{as $\v\to 0+$},
\end{align}
which, together with \eqref{8.8},  yields  that
\begin{align*}
&\lim_{\v\to 0+} \int_{0}^{T}\int_0^\infty  |(\Phi_r^{\v}-\Phi_r)(t,r)|^2\,r^{n-1} \dd r
\dd t\\
&= \lim_{\v\to 0+} \int_{0}^{T}\int_0^\infty \frac{1}{r^{n-1}} \Big|\int_0^{r} (\rho^{\v}-\rho)(t,z)\, z^{n-1} \dd z\Big|^2 \dd r \dd t\nonumber\\
&\leq C(M) Tk^{-n+2}+\lim_{\v\to 0+} \int_{0}^{T}\int_0^k \frac{1}{r^{n-1}} \Big|\int_0^{r} (\rho^{\v}-\rho)(t,z)\,z^{n-1} \dd z\Big|^2 \dd r\dd t\\
&\leq C(M,T)\,k^{-n+2}.
\end{align*}
Then \eqref{5.4-1} follows by taking $k\to \infty$ to obtain
\begin{equation}\label{8.10}
\lim_{\v\to 0+} \int_{0}^{T}\int_0^\infty  |(\Phi_r^{\v}-\Phi_r)(t,r)|^2\,r^{n-1}\dd r
\dd t=0 \qquad \mbox{if $\g>\frac{2n}{n+2}$}.
\end{equation}

\smallskip
3. Now we define
\begin{align*}
(\rho,\M,\Phi)(t,\mathbf{x})=(\rho(t,r), m(t,r) \frac{\mathbf{x}}{r},\Phi(t,r)).
\end{align*}
Using \eqref{8.7},  Fatou's lemma, and \eqref{8.10}, we have
\begin{align*}
&\int_{t_1}^{t_2}\int_{\R^n}
\Big(\frac{1}{2}\Big|\frac{\M}{\sqrt{\rho}}\Big|^2 + \rho e(\rho)-\frac{\kappa}{2}  |\nabla\Phi|^2\Big)(t,\mathbf{x})\,\dd\mathbf{x} \dd t\nonumber\\
&\leq (t_2-t_1)\int_{\R^n}\Big(\frac{1}{2}\Big|\frac{\M_0}{\sqrt{\rho_0}}\Big|^2
+\rho_0 e(\rho_0)-\frac{\kappa}{2}|\nabla\Phi_0|^2\Big)(\mathbf{x})\,\dd\mathbf{x},
\end{align*}
which implies that, for {\it a.e.} $t\geq0$,
\begin{align*}
&\int_{\R^n}\Big(\frac{1}{2}\Big|\frac{\M}{\sqrt{\rho}}\Big|^2 + \rho e(\rho)-\frac{\kappa}{2}|\nabla\Phi|^2\Big)(t,\mathbf{x})\,\dd\mathbf{x}\nonumber\\
&\leq \int_{\R^n}\Big(\frac{1}{2}\Big|\frac{\M_0}{\sqrt{\rho_0}}\Big|^2 + \rho_0 e(\rho_0)
-\frac{\kappa}{2}|\nabla\Phi_0|^2\Big)(\mathbf{x})\,\dd \mathbf{x}.
\end{align*}
This leads to \eqref{1.20-1}--\eqref{1.20-2}.

\smallskip
4. In the following, we prove that $(\rho,\M, \Phi)(t,\mathbf{x})$ is  indeed a global weak solution of
problem \eqref{1.1-1}--\eqref{1.1-3} in $\mathbb{R}^n$.

Let $\zeta(t,\mathbf{x})\in C_0^1(\mathbb{R}\times \mathbb{R}^n)$ be a smooth function with compact support.
Then it follows from \eqref{5.39} that
\begin{align}\label{7.60}
\int_{\mathbb{R}^{n+1}_+}\big(\rho^{\v} \zeta_t + \M^{\v}\cdot\nabla \zeta\big)\, \dd\mathbf{x}\dd t
+\int_{\mathbb{R}^n}\rho_0^\v(\mathbf{x}) \zeta(0,\mathbf{x})\, \dd\mathbf{x}=0.
\end{align}
Without loss of generality, we assume that ${\rm supp}\,\zeta \subset [-T,T]\times B_D(\mathbf{0})$ for some $T,D>0$.
Let $\phi(t,r)$ be the corresponding function defined in \eqref{3.39-1}.
Using \eqref{3.43-1}, \eqref{7.51}, and similar arguments as in the proof of Lemma \ref{lem5.8},
we see that, for any fixed $\sigma>0$,
\begin{align}\label{3.116}
	&\lim_{\v\rightarrow0+}\int_{0}^{\infty} \int_{\mathbb{R}^n\backslash B_\s(\mathbf{0})}
	\big(\rho^{\v} \zeta_t + \M^{\v} \cdot\nabla\zeta\big)\,\dd\mathbf{x}\dd t\\
	&=\lim_{\v\rightarrow0+}\int_{0}^{\infty}\int_\sigma^\infty \big(\rho^{\v} \phi_t + m^{\v} \phi_r\big)
	\,\omega_n r^{n-1}\dd r\dd t\nonumber\\
	&=\int_{0}^{\infty}\int_\sigma^\infty  \big(\rho \phi_t + m \phi_r\big)\,\omega_n r^{n-1} \dd r \dd t\nonumber\\
	&=\int_{0}^{\infty} \int_{\mathbb{R}^n\backslash B_\s(\mathbf{0})} \big(\rho\zeta_t + \M \cdot\nabla\zeta\big)\, \dd\mathbf{x}\dd t.\nonumber
\end{align}
Using \eqref{7.83} and similar arguments as in \eqref{5.42-1}, we have
\begin{align}\label{3.117}
&\Big|\int_{0}^{\infty} \int_{B_\s(\mathbf{0})}  \big(\rho^{\v}-\rho\big) \zeta_t \;\dd\mathbf{x}\dd t\Big|\\
&\leq C(\|\zeta\|_{C^1},T)\Big\{\int_0^T\int_0^\sigma\big((\rho^{\v})^{\gamma}+ \rho^{\gamma}\big)\,r^{n-1} \dd r\dd t \Big\}^{\frac1\gamma}
\,\Big\{\int_0^\sigma r^{n-1} \dd r \Big\}^{1-\frac1\gamma}\nonumber\\
&\leq C(E_0,M,\|\zeta\|_{C^1},T)\, \sigma^{n(1-\frac1\gamma)}\longrightarrow 0\qquad \mbox{as $\sigma\rightarrow0$},\nonumber
\end{align}
and
\begin{align}\label{3.117-1}
	&\Big|\int_{0}^{T} \int_{B_\s(\mathbf{0})}  \big(\M^{\v}-\M\big) \cdot\nabla \zeta\,\dd\mathbf{x} \dd t\Big|\\
	&\leq C\Big\{\int_{0}^{T}\int_0^\sigma \Big(\frac{|m^\v|^2}{\rho^\v}+\frac{m^2}{\rho}\Big)(t,r)|\phi_r|\,r^{n-1}\dd r \dd t\Big\}^{\frac12}\nonumber\\
	&\qquad\times \Big\{\int_{0}^{T}\int_0^\sigma \big(\rho^{\v}+\rho\big)(t,r)\, |\phi_r| r^{n-1}\dd r \dd t\Big\}^{\frac12}\nonumber\\
	&\leq C(E_0,M,\|\zeta\|_{C^1},T)\sigma^{\frac{n}{2}(1-\frac1\gamma)}\longrightarrow 0\qquad \mbox{as $\sigma\rightarrow0$},\nonumber
\end{align}
which, together with \eqref{3.116}--\eqref{3.117-1}, yields that
\begin{align}\label{3.118}
\lim_{\v\rightarrow0+}\int_{0}^{\infty} \int_{\mathbb{R}^n} \big(\rho^{\v} \zeta_t + \M^{\v} \cdot\nabla\zeta\big)\,\dd\mathbf{x} \dd t
=\int_{0}^{\infty} \int_{\mathbb{R}^n} \big(\rho \zeta_t + \M \cdot\nabla\zeta\big)\,\dd\mathbf{x} \dd t.
\end{align}
Taking $\v\rightarrow0+$ in \eqref{7.60} and using \eqref{3.118}, we conclude that $(\rho,\M)$ satisfies \eqref{1.18-1}.

Next, we consider the momentum equation.
Let $\bp=(\psi_1,\cdots,\psi_n) \in \big(C_0^2(\mathbb{R}\times\mathbb{R}^n)\big)^n$ be a smooth function with compact support,
and let $\chi_\sigma(r)\in C^\infty(\mathbb{R})$ be a cut-off function satisfying \eqref{3.49-1}.
Without loss of generality, we assume that ${\rm supp}\,\bp \subset [-T,T]\times B_D(\mathbf{0})$ for some $T,D>0$.
Denote $\Psi_{\sigma}=\bp \chi_\sigma$. Then we have
\begin{align}\label{7.62}
&\Big| \v \int_{\R_+^{n+1}}\Big\{\frac{1}{2}\M^{\v}\cdot \big(\Delta \Psi_{\sigma}+\nabla\mbox{\rm div}\,\Psi_{\sigma}\big)
+ \frac{\M^{\v}}{\sqrt{\rho^{\v}}} \cdot \big(\nabla\sqrt{\rho^{\v}}\cdot \nabla\big)\Psi_{\sigma}\\
&\hspace{6cm}+ \nabla\sqrt{\rho^{\v}}  \cdot \big(\frac{\M^{\v}}{\sqrt{\rho^{\v}}}\cdot \nabla\big)\Psi_{\sigma}\Big\}\dd\mathbf{x} \dd t \Big|\nonumber\\
&=\Big|\sqrt{\v}\int_{\mathbb{R}^{n+1}_+} \sqrt{\rho^{\v}} \Big\{V^{\v}  \frac{\mathbf{x}\otimes\mathbf{x}}{r^2}
+\frac{\sqrt{\v}}{r}\sqrt{\rho^\v}u^{\v}\Big(I_{n\times n}-\frac{\mathbf{x}\otimes\mathbf{x}}{r^2}\Big)\Big\}: \nabla\Psi_{\sigma}\,\dd\mathbf{x} \dd t\Big|\nonumber\\
&\leq C\Big\{\int_{0}^\infty\int_{\mathbb{R}^n}|V^\v|^2\,\dd\mathbf{x} \dd t+\v\int_{\mathbb{R}_+^2}  \frac{|m^\v|^2}{\rho^\v} r^{n-3}\,\dd r\dd t\Big\}^{\frac12}
\Big\{\v\int_0^\infty \int_{\mathbb{R}^n} \rho^\v |\nabla\Psi_{\sigma}|^2\, \dd\mathbf{x} \dd t\Big\}^{\frac12}\nonumber\\
&\leq C(\sigma,E_0, D,T)\sqrt{\v}\longrightarrow 0\qquad \mbox{as $\v\rightarrow 0+$}.\nonumber
\end{align}
For the potential term, it follows from \eqref{7.51} and \eqref{7.102} that
\begin{align}\label{7.104}
\lim_{\v\rightarrow0+} \int_{\mathbb{R}^{n+1}_+} \rho^\v \nabla\Phi^\v\cdot \Psi_\sigma\,\dd\mathbf{x} \dd t
=\int_{\mathbb{R}^{n+1}_+} \rho \nabla\Phi\cdot \Psi_\sigma\, \dd\mathbf{x} \dd t.
\end{align}
Using \eqref{7.62}--\eqref{7.104}  and passing limit $\v\rightarrow0+$ (up to a subsequence) in \eqref{3.59-1} yield that
\begin{align}\label{7.63}
	&\int_{\mathbb{R}^{n+1}_+}\Big\{\M\cdot\partial_t\Psi_\sigma +\frac{\M}{\sqrt{\rho}} \cdot \big(\frac{\M}{\sqrt{\rho}} \cdot \nabla\big)\Psi_{\sigma}
	+p(\rho)\, \mbox{div}\,\Psi_{\sigma} - \rho \nabla\Phi \cdot \Psi_\sigma\Big\}\, \dd\mathbf{x} \dd t\\
	& +\int_{\mathbb{R}^n} \M_0\cdot \Psi_\s(0,\mathbf{x})\,\dd\mathbf{x}  =0.\nonumber
\end{align}

In the following, we take limit $\sigma\rightarrow0$ in \eqref{7.63}.
Notice that, for any $T>0$ and $D>0$,
\begin{equation}\label{3.122}
\int_0^T\int_0^D \big(\frac{m^2}{\rho}+ p(\rho)\big)(t,r)\,r^{n-1} \dd r\dd t\leq C(E_0,M,D,T),
\end{equation}
which, together with similar arguments as in \eqref{3.61}, yields that
\begin{align}\label{3.123}
	&\lim_{\sigma\rightarrow0+} \bigg(\int_{\mathbb{R}^{n+1}_+} \M \cdot\partial_t\Psi_\sigma\, \dd\mathbf{x} \dd t
	+\int_{\mathbb{R}^n} \M_0\cdot \Psi_\s(0,\mathbf{x})\,\dd\mathbf{x} \bigg)\\
	&=\int_{\mathbb{R}^{n+1}_+}  \M \cdot\partial_t\bp\, \dd\mathbf{x} \dd t
	+\int_{\mathbb{R}^n} \M_0\cdot \bp(0,\mathbf{x})\,\dd\mathbf{x}.\nonumber
\end{align}
Using  \eqref{3.63-2}--\eqref{3.63-1} and \eqref{3.122}, we have
\begin{align}
&\Big|\int_{\mathbb{R}^{n+1}_+} \big(\frac{m^2}{\rho} +p(\rho) \big) \, \bp\cdot \frac{\mathbf{x}}{r}\, \chi'_\sigma(r)\,
\dd\mathbf{x} \dd t\Big|\nonumber\\
&\leq \int_0^\infty\int_{\sigma}^{2\sigma} \big(\frac{m^2}{\rho}+p(\rho)\big)\,|\varphi(t,r) \chi'_\sigma(r)|\, r^{n-1}\dd r \dd t\nonumber\\
&\leq C\int_0^T\int_{\sigma}^{2\sigma}  \big(\frac{m^2}{\rho}+p(\rho) \big)\,r^{n-1} \dd r \dd t\longrightarrow 0\qquad \mbox{as $\sigma\rightarrow0+$},
\nonumber
\end{align}
which, together with \eqref{3.122} and Lebesgue's dominated convergence theorem, yields that
\begin{align}\label{3.125}
&\lim_{\sigma\rightarrow0+} \int_{\mathbb{R}^{n+1}_+} \bigg\{\frac{\M}{\sqrt{\rho}} \cdot \big(\frac{\M}{\sqrt{\rho}} \cdot \nabla\big)\Psi_{\sigma}
+ p(\rho)\,\mbox{div}\,\Psi_{\sigma} \bigg\}\,\dd\mathbf{x} \dd t\\
&= \int_{\mathbb{R}^{n+1}_+} \bigg\{ \frac{\M}{\sqrt{\rho}} \cdot \big(\frac{\M}{\sqrt{\rho}} \cdot \nabla\big)\bp
+ p(\rho)\, \mbox{div}\,\bp\bigg\}\, \dd\mathbf{x} \dd t.\nonumber
\end{align}

For the term involving the potential, using \eqref{3.63-2}--\eqref{3.63-1}, we see that
\begin{align}\nonumber
\Big|\rho(t,r) \big(\int_0^r  \rho(t,z)\,z^{n-1}\dd z\big)\, \varphi(t,r) \chi_{\sigma}(r)\Big|
\leq C(\|\bp\|_{C^1})\, \rho(t,r) r \int_0^r  \rho(t,z)\, z^{n-1}\dd z,
\end{align}
which, together with \eqref{7.100}, \eqref{7.101}, and Lebesgue's dominated convergence theorem, yields that
\begin{align}\label{7.105}
&\lim_{\sigma\rightarrow0+}\int_{\mathbb{R}^{n+1}_+} \rho \nabla\Phi \cdot \Psi_\sigma\, \dd\mathbf{x}\dd t\\
&=\lim_{\sigma\rightarrow0+}\int_{\R_+^2} \rho\Phi_r \, \varphi(t,r) \chi_{\sigma}(r)\, \omega_n r^{n-1}\dd r \dd t\nonumber\\
&=\kappa \lim_{\sigma\rightarrow0+} \int_{\R_+^2}  \rho(t,r)
\big(\int_0^r  \rho(t,z)\,\omega_n z^{n-1}\dd z\big)\varphi(t,r) \chi_{\sigma}(r)\, \dd r\dd t\nonumber\\
&=\kappa \int_{\R_+^2}  \rho(t,r)   \big(\int_0^r  \rho(t,z)\,\omega_n z^{n-1}\dd z\big) \varphi(t,r)\,\dd r \dd t\nonumber\\
&=\int_{\mathbb{R}^{n+1}_+} \rho \nabla\Phi \cdot \bp\, \dd\mathbf{x} \dd t.\nonumber
\end{align}
Substituting \eqref{3.123}--\eqref{7.105}
into \eqref{7.63}, we conclude that $(\rho,\M,\Phi)$ satisfies \eqref{1.19-1}. By the Lebesgue theorem,
we can weaken the assumption that $\bp\in (C_0^2)^n$ as $\bp\in (C_0^1)^n$.

Finally, we consider the Poisson equation.
Let  $\xi(\mathbf{x})\in C_0^1(\mathbb{R}^n)$ be any smooth function with compact support.
For any $t_2> t_1\geq 0$,  we use \eqref{3.74}, \eqref{7.102}, and similar arguments as in \eqref{3.118},
and then pass limit $\v\rightarrow0+$ (up to a subsequence) to obtain
\begin{align} \label{7.106}
-\int_{t_1}^{t_2}\int_{\mathbb{R}^n} \nabla\Phi(s,\mathbf{x}) \cdot \nabla \xi(\mathbf{x})\,\dd\mathbf{x} \dd s
=\kappa\int_{t_1}^{t_2}\int_{\mathbb{R}^n} \rho(s,\mathbf{x}) \xi(\mathbf{x})\,\dd\mathbf{x} \dd s.
\end{align}
Applying the Lebesgue point theorem, we obtain that, for {\it a.e.} $t\geq0$,
\begin{align}
&\lim_{t_2,t_1\rightarrow t}\frac{1}{t_2-t_1}\int_{t_1}^{t_2}
\int_{\mathbb{R}^n} \nabla\Phi(s,\mathbf{x}) \cdot \nabla \xi(\mathbf{x})\,\dd\mathbf{x} \dd s
=\int_{\mathbb{R}^n} \nabla\Phi(t,\mathbf{x}) \cdot \nabla \xi(\mathbf{x})\, \dd\mathbf{x},\label{7.107}\\[2mm]
&\lim_{t_2,t_1\rightarrow t}\frac{1}{t_2-t_1}\int_{t_1}^{t_2}\int_{\mathbb{R}^n} \rho(s,\mathbf{x}) \xi(\mathbf{x})\, \dd\mathbf{x} \dd s
=\int_{\mathbb{R}^n} \rho(t,\mathbf{x}) \xi(\mathbf{x})\, \dd\mathbf{x}.\label{7.108}
\end{align}
Combining \eqref{7.106}--\eqref{7.108} together, we conclude that $(\rho,\M,\Phi)$ satisfies \eqref{1.19-2}.
$\hfill\Box$

\appendix
\section{Sobolev's Inequality and Construction of the Approximate Initial Data Sequences}\setcounter{equation}{0}

In this appendix, we first state Sobolev's inequality used in \S3--\S4; see {\it e.g.} \cite[\S 8.3]{Lieb-Loss}.

\begin{lemma}[Sobolev's Inequality]\label{Sobolev}
	For $n\geq 3$, let $\nabla f\in L^2(\mathbb{R}^n)$ and $\lim_{|\mathbf{x}|\rightarrow\infty}f(\mathbf{x})=0$.
	Then
	\begin{equation}\label{A.1}
		\|f\|^2_{L^{\frac{2n}{n-2}}}\leq A_n \|\nabla f\|_{L^2}^2,
	\end{equation}
	where $A_n$ is the best constant which is given by
	\begin{equation}\label{A.2}
		A_n=\frac{4}{n(n-2)} \omega_{n+1}^{-\frac{2}{n}}
	\end{equation}
	with $\omega_{n+1}=\frac{2\pi^{\frac{n+1}{2}}}{\Gamma(\frac{n+1}{2})}$ as the surface area of unit sphere in $\mathbb{R}^{n+1}$.
\end{lemma}

We now construct the approximate initial data sequences $(\rho^\v_0,m_0^\v)$ and  $(\rho^{\v,b}_0, u^{\v,b}_0)$
with desired estimates, regularity, and boundary compatibility.

To keep the $L^p$--properties of mollification, it is more convenient to smooth out the initial data
in the original coordinates in $\mathbb{R}^n$;
so we do not distinguish functions $(\rho_0,m_0)(r)$ from $(\rho_0,m_0)(\mathbf{x})=(\rho_0,m_0)(|\mathbf{x}|)$ for simplicity below.

For the initial data $(\rho_0,m_0)$, we assume that
\begin{align}\label{A.3}
	\begin{split}
		&\int_{\mathbb{R}^n} \Big(\rho_0+ \rho_0^\gamma+ \rho_0^{\frac{2n}{n+2}}
		+\Big|\frac{m_0}{\sqrt{\rho_0}}\Big|^2 \Big)\,\dd\mathbf{x}<\infty\quad \mbox{for $\kappa=-1$ (plasmas) with $\g>1$},\\
		&\int_{\mathbb{R}^n} \Big(\rho_0+ \rho_0^\gamma +\Big|\frac{m_0}{\sqrt{\rho_0}}\Big|^2 \Big)\,\dd\mathbf{x}<\infty\qquad\qquad
		\mbox{for $\kappa=1$ (gaseous stars) with $\g>\frac{2n}{n+2}$},
	\end{split}
\end{align}
and denote
\begin{equation*}\label{A.3-1}
	\Delta\Phi_0=\kappa \rho_0,
\end{equation*}
which is well-defined, under assumption \eqref{A.3}.

From now on, we denote $C>0$ is a universal constant independent of $\v,\delta$, and $b$.

Let  $J(\mathbf{x})$ be the standard mollification function and
$\dis J_{\delta}(\textbf{x}):=\frac{1}{\delta^n}J(\frac{\mathbf{x}}{\delta})$ for $\delta\in(0,1)$.
For later use, we take $\delta=\v^{\frac12}$ and define  $\tilde{\rho}_{0}^{\v}(\mathbf{x})$  as
\begin{align}\label{A.4}
	\tilde{\rho}_{0}^{\v}(\mathbf{x})
	:=\Big(\int_{\mathbb{R}^n} \sqrt{\rho_{0}(\mathbf{x-y})} J_{\sqrt{\v}}(\mathbf{y})\,\dd\mathbf{y}+\v e^{-|\mathbf{x}|^2}\Big)^2.
\end{align}
Then $\tilde{\rho}_{0}^{\v}(\mathbf{x})$ is still a spherically symmetric function,
{\it i.e.}, $\tilde{\rho}_{0}^{\v}(\mathbf{x})=\tilde{\rho}_{0}^{\v}(|\mathbf{x}|)$.
It is also direct to know that $\tilde{\rho}_{0}^\v(\mathbf{x})\geq  \v^2 e^{-2|\mathbf{x}|^2}>0$.

\begin{lemma}\label{lemA.2}
	Let $q\in \{1, \gamma\}$ for $\kappa=1$ $($gaseous stars$)$, and $q\in \{1,\gamma, \frac{2n}{n+2}\}$ for $\kappa=-1$ $($plasmas$)$. Then
	\begin{align}
		&\|\tilde{\rho}_0^\v\|_{L^q}\leq \|\rho_0\|_{L^q}+C \v \qquad\mbox{as $\v\in(0,1]$},\label{A.5}\\
		&\lim_{\v\rightarrow0+} \big\{\|\tilde{\rho}_0^\v-\rho_0\|_{L^q}+\|\sqrt{\tilde{\rho}_0^\v}-\sqrt{\rho_0}\|_{L^{2q}}\big\}=0,\label{A.6}\\
		&\v^2\int_{\R^n}\Big|\nabla_{\mathbf{x}}\sqrt{\tilde{\rho}^\v_{0}(\mathbf{x})}\Big|^2\dd\mathbf{x}
		\leq C\v \big(\|\rho_0\|_{L^1}+1\big)\longrightarrow0 \qquad\mbox{as $\v\rightarrow0+$}.\label{A.7}
	\end{align}
\end{lemma}

\begin{proof}
 It is direct to see that \eqref{A.5}--\eqref{A.6} follow
from the standard property of mollifier operator.
For \eqref{A.7}, we notice that
\begin{align}
	\v^2 \int_{\R^n}\Big|\nabla_{\mathbf{x}}\sqrt{\tilde{\rho}^\v_{0}(\mathbf{x})}\Big|^2\,\dd\mathbf{x}
	&\leq C\v^2\int_{\R^n}\left|\int_{\mathbb{R}^n} \sqrt{\rho_{0}(\mathbf{x-y})}
	\nabla_{\mathbf{y}}J_{\sqrt{\v}}(\mathbf{y}) \dd\mathbf{y}\right|^2 \dd\mathbf{x}+C\v^4 \nonumber\\[1mm]
	& \leq C \|\rho_0\|_{L^1}\v+C\v^4\longrightarrow0\qquad\,\,\,\,\mbox{as $\v\rightarrow0+$}.\nonumber
\end{align}
\end{proof}

In general, since $\int_{\R^n} \tilde{\rho}_{0}^{\v}(\mathbf{x})\,\dd\mathbf{x}
\neq \int_{\R^n} \rho_0(\mathbf{x})\,\dd\mathbf{x} =M$,  we define
\begin{align}\label{A.7-1}
	\rho_0^{\v}(\mathbf{x}):=\frac{M }{\int_{\R^n} \tilde{\rho}_{0}^{\v}(\mathbf{x})\,\dd\mathbf{x}}
	\tilde{\rho}_{0}^{\v}(\mathbf{x}).
\end{align}
Combining Lemma \ref{lemA.2} and \eqref{A.7-1}, we have

\begin{lemma}
	Let $q=1$ and $\gamma>1$ for $\kappa=1$ $($gaseous stars$)$, and $q=1,\gamma$, and $\frac{2n}{n+2}$ for $\kappa=-1$ $($plasmas$)$.
	Then
	\begin{align}
		&\int_{\R^n}\rho_0^\v(\mathbf{x})\,\dd\mathbf{x}=M, \qquad  \|\rho_0^\v\|_{L^q}\leq C\big(\|\rho_0\|_{L^p}+1\big)
		\qquad\mbox{for all $\v\in(0,1]$},\label{A.5-1}\\
		&\lim_{\v\rightarrow0+} \big\{\|\rho_0^\v-\rho_0\|_{L^q}+\|\sqrt{\rho_0^\v}-\sqrt{\rho_0}\|_{L^{2q}}\big\}=0,\label{A.6-1}\\
		&\v^2\int_{\R^n}\Big|\nabla_{\mathbf{x}}\sqrt{\rho^\v_{0}(\mathbf{x})}\Big|^2 \dd\mathbf{x}
		\leq C\v \big(\|\rho_0\|_{L^1}+1\big)\longrightarrow0\qquad\mbox{as $\v\rightarrow0+$}.\label{A.7-2}
	\end{align}
\end{lemma}

We define
\begin{align}\label{A.8}
	\Delta\Phi_0^\v=\kappa \rho_{0}^\v.
\end{align}
Then a direct calculation yields that
\begin{lemma}\label{lemA.3}
	$\Phi_0^\v$ satisfy
	\begin{align}\label{A.9}
		\|\nabla_{\mathbf{x}}(\Phi_0^\v-\Phi_0)\|_{L^2}\leq C\|\rho_0^\v-\rho_0\|_{L^{\frac{2n}{n+2}}}\longrightarrow0
		\qquad\mbox{as $\v\rightarrow0+$}.
	\end{align}
\end{lemma}

From \eqref{A.4}, we know that $\rho_{0}^{\v}(\mathbf{x})$ is a good approximation.
However, we don't know yet whether
\begin{equation}\label{A.18}
	\rho_{0}^{\v}(b)\cong b^{-n+\alpha}\qquad \mbox{with $\alpha:=\min\{\frac12,(1-\frac1\gamma)n\}$}
\end{equation}
is satisfied. In fact,  \eqref{A.18} (which agrees with condition \eqref{2.60-1}) is required
in the proof of Lemmas \ref{lem2.2}--\ref{lem2.3}.
To solve this problem, we denote $S=S(z)\in C^\infty(\mathbb{R})$ to be a cut-off function satisfying
\begin{align}\label{A.19-1}
	\begin{cases}
		S(z)=0\quad   \mbox{if $z\in(-\infty,0]$},\\
		S(z)=1\quad  \mbox{if $z\in[1,\infty)$},\\
		S(z) \quad  \mbox{is monotonic increasing in $[0,1]$}.
	\end{cases}
\end{align}
Now we define $\tilde{\rho}_{0}^{\v,b}(r)$ by
\begin{align}\label{A.19}
	\tilde{\rho}_{0}^{\v,b}(\mathbf{x})
	&:=\Big\{\sqrt{\rho_{0}^{\v}(\mathbf{x})}\,\big\{1-S(2 (|\mathbf{x}|-(b-1)))\big\}
	+b^{-\frac{n-\alpha}{2}}\, S(2 (|\mathbf{x}|-(b-1)))\Big\}^2.
\end{align}
It is direct to check that $\tilde{\rho}_{0}^{\v,b}(b)=b^{-(n-\alpha)}$,
which clearly satisfies \eqref{A.18}.

\begin{lemma}\label{lemA.6}
	The smooth functions $\tilde{\rho}_{0}^{\v,b}(\mathbf{x})$ defined in \eqref{A.19} satisfy \eqref{2.60-1}  and
	\begin{align}
		&\int_{|\mathbf{x}|\leq b} \Big(|\tilde{\rho}_0^{\v,b}(\mathbf{x})-\rho_0^\v(\mathbf{x})|^q
		+\Big| \sqrt{\tilde{\rho}_{0}^{\v,b}(\mathbf{x})}
		-\sqrt{\rho_{0}^{\v}(\mathbf{x})}\Big|^{2q}\Big)\,\dd \mathbf{x}\rightarrow0 \,\,\,\mbox{as $b\to \infty$},\label{A.20}\\[1mm]
		&\v^2\int_{|\mathbf{x}|\leq b}\Big|\nabla_{\mathbf{x}}\sqrt{\tilde{\rho}_0^{\v,b}(\mathbf{x})}\Big|^2 \dd\mathbf{x}
		\leq C\big(1+\|\rho_0\|_{L^1}\big)\,\v, \label{A.20-1}
	\end{align}
	where $q\in \{1,\gamma\}$ for $\kappa=1$ $($gaseous stars$)$, and $q\in \{1,\gamma, \frac{2n}{n+2}\}$ for $\kappa=-1$ $($plasmas$)$.
\end{lemma}

\begin{proof}
 Using \eqref{A.19}, a direct calculation shows that
\begin{align}\label{A.21}
	\int_{|\mathbf{x}|\leq b} |\tilde{\rho}_0^{\v,b}(\mathbf{x})|^q\,\dd\mathbf{x}
	&\leq C\int_{|\mathbf{x}|\leq b} |\rho_0^{\v}(\mathbf{x})|^q\,\dd\mathbf{x}+C b^{-n+\alpha}\int_{b-1}^{b} r^{n-1} \dd r\\[1mm]
	&\leq C\big(\|\rho_{0}\|_{L^q}+1+b^{-1+\alpha}\big)\nonumber\\[1mm]
    &\leq C\big(\|\rho_{0}\|_{L^q}+1\big).\nonumber
\end{align}
Using \eqref{A.5} and \eqref{A.21}, we have
\begin{align}\label{A.22}
	&\int_{|\mathbf{x}|\leq b} \Big| \tilde{\rho}_{0}^{\v,b}(\mathbf{x})
	-\rho_{0}^{\v}(\mathbf{x})\Big|^{q} \dd \mathbf{x}
	+\int_{|\mathbf{x}|\leq b} \Big| \sqrt{\tilde{\rho}_{0}^{\v,b}(\mathbf{x})}
	-\sqrt{\rho_{0}^{\v}(\mathbf{x})}\Big|^{2q}  \dd \mathbf{x}\\
	&\leq \Big(\int_{|\mathbf{x}|\leq b} \Big| \sqrt{\tilde{\rho}_{0}^{\v,b}(\mathbf{x})}
	-\sqrt{\rho_{0}^{\v}(\mathbf{x})}\Big|^{2q} \dd\mathbf{x}\Big)^{\f12}
	\Big(\int_{|\mathbf{x}|\leq b} \Big| \sqrt{\tilde{\rho}_{0}^{\v,b}(\mathbf{x})}
	+\sqrt{\rho_{0}^{\v}(\mathbf{x})}\Big|^{2q}  \dd\mathbf{x}\Big)^{\f12}\nonumber\\
	&\leq  C\big(\|\rho_{0}\|_{L^q}+1\big)
	\Big(\int_{b-1\leq |\mathbf{x}|\leq b} |\rho_{0}^{\v}(\mathbf{x})|^{q}\,\dd\mathbf{x}
	+b^{-n+\alpha}\int_{b-1}^{b}r^{n-1} \dd r\Big)^{\f12}\nonumber\\
	&\leq C\big(\|\rho_{0}\|_{L^q}+1\big)
	\Big(\int_{b-2\leq|\mathbf{x}|\leq b+1} \big( |\rho_{0}(\mathbf{x})|^{q} +e^{-r^2}\big)\,r^{n-1} \dd r
	+b^{-1+\alpha} \Big)^{\f12} \longrightarrow 0\nonumber
\end{align}
as $b\rightarrow\infty$, which implies \eqref{A.20}.

For \eqref{A.20-1}, a direct calculation shows that
\begin{align}
	\nabla\sqrt{\tilde{\rho}_{0}^{\v,b}(\mathbf{x})}
	&=\nabla\sqrt{\rho_{0}^{\v}(\mathbf{x})}\,\big(1-S(2 (|\mathbf{x}|-(b-1)))\big)\nonumber\\
	&\quad-2\big(\sqrt{\rho_{0}^{\v}(\mathbf{x})}-\sqrt{b^{-n+\alpha}}\big)\, S'(2(|\mathbf{x}|-(b-1)))\frac{\mathbf{x}}{|\mathbf{x}|},
	\nonumber
\end{align}
which, together with  \eqref{A.7}, yields that
\begin{align}
	&\v^2\int_{0}^{b} \Big|\nabla\sqrt{\tilde{\rho}_0^{\v,b}(\mathbf{x})}\Big|^2\,\dd\mathbf{x}\nonumber\\
	&\leq C\v^2\int_{0}^{b} \Big|\nabla\sqrt{\rho_0^{\v}(\mathbf{x})}\Big|^2\, \dd \mathbf{x}
	+C\v^2\int_{b-1\leq |\mathbf{x}|\leq b} \big(\rho_0^{\v}(\mathbf{x})+b^{-n+\alpha}\big)\, \dd \mathbf{x}\nonumber\\
	&\leq C\big(\|\rho_0\|_{L^1}+1\big)\v^2.\nonumber
\end{align}
\end{proof}

In general, since $\dis \int_{b^{-1}\leq |\mathbf{x}|\leq b} \tilde{\rho}_{0}^{\v,b}(\mathbf{x})\,\dd\mathbf{x}\neq M$, we define
\begin{align}\label{A.23-1}
	\rho_0^{\v,b}(\mathbf{x}):=\frac{M }{\int_{b^{-1}\leq |\mathbf{x}|\leq b} \tilde{\rho}_{0}^{\v,b}(\mathbf{x})\,\dd\mathbf{x}}
	\tilde{\rho}_{0}^{\v,b}(\mathbf{x}).
\end{align}
Combining \eqref{A.23-1} and  Lemma \ref{lemA.6}, we have

\begin{lemma}\label{lemA.6-1}
	The smooth function $\rho_{0}^{\v,b}(\mathbf{x})$ defined in \eqref{A.23-1} satisfies \eqref{2.60-1}
	and
	\begin{align}
		&\int_{b^{-1}\leq |\mathbf{x}|\leq b} \rho_{0}^{\v,b}(\mathbf{x})\,\dd\mathbf{x}= M \qquad\mbox{for all $\v\in (0,1]$ and $b>1$},\\
		&\int_{|\mathbf{x}|\leq b} \Big(\big|\rho_0^{\v,b}(\mathbf{x})-\rho_0^\v(\mathbf{x})\big|^q
		+\Big| \sqrt{\rho_{0}^{\v,b}(\mathbf{x})}
		-\sqrt{\rho_{0}^{\v}(\mathbf{x})}\Big|^{2q}\Big)\,\dd \mathbf{x}\rightarrow0 \,\,\mbox{as $b\to \infty$},\label{A.20-2}\\
		&\v\int_{|\mathbf{x}|\leq b}\Big|\nabla_{\mathbf{x}}\sqrt{\rho_0^{\v,b}(\mathbf{x})}\Big|^2 \dd\mathbf{x}
		\leq C\big(\|\rho_0\|_{L^1}+1\big),\label{A.20-3}
	\end{align}
	where $q\in\{1,\gamma\}$ for $\kappa=1$ $($gaseous stars$)$, and $q\in \{1,\gamma, \frac{2n}{n+2}\}$ for $\kappa=-1$ $($plasmas$)$.
\end{lemma}

We define
\begin{align}\label{A.50}
	\Delta\Phi_0^{\v,b}=\kappa \rho_{0}^{\v,b}\,\mathbf{1}_{\{b^{-1}\le |\mathbf{x}|\leq b\}}(\mathbf{x}),
\end{align}
where ${\mathbf{1}}_{\{b^{-1}\le |\mathbf{x}|\leq b\}}(\mathbf{x})$ is the indicator function of set $\{b^{-1}\le |\mathbf{x}|\leq b\}$.
Using \eqref{A.20-2} and  by a direct calculation yield
\begin{lemma}\label{lemA.9}
	$\Phi_0^{\v,b}$ satisfy
	\begin{align}\label{A.51a}
		\|\nabla_{\mathbf{x}}(\Phi_0^{\v,b}-\Phi_0^\v)\|_{L^2}
		\leq C\|\rho_0^{\v,b}\,\mathbf{1}_{\{b^{-1}\le |\mathbf{x}|\leq b\}}(\mathbf{x})-\rho_0^\v\|_{L^{\frac{2n}{n+2}}}\rightarrow0
		\,\,\mbox{as $b\rightarrow \infty$}.
	\end{align}
\end{lemma}

Next, we construct the approximate initial data for the velocity.
We denote $\mathbf{1}_{[4b^{-1},b-2]}$ to be the indicator
function of $\{\mathbf{x}\in\mathbb{R}^n\ :\ 4b^{-1}\leq|\mathbf{x}|\leq b-2\}$
and define $u_0^\v(\mathbf{x})$ and $\tilde{u}_{0}^{\v,b}(\mathbf{x})$ as
\begin{align}
	&u_{0}^{\v}(\mathbf{x})
	:=\frac{1}{\sqrt{\rho_0^{\v}(\mathbf{x})}} (\frac{m_0}{\sqrt{\rho_0}})(\mathbf{x}),\label{A.31-1}\\
	&\tilde{u}_{0}^{\v,b}(\textbf{x})
	:=\frac{1}{\sqrt{\rho_0^{\v,b}(\mathbf{x})}}\int_{\mathbb{R}^n}\Big(\frac{m_0 \mathbf{1}_{[4b^{-1},b-2]}}{\sqrt{\rho_0}}\Big)
	(\mathbf{x}-\mathbf{y})\,J_{b^{-1}}(\textbf{y})\,\dd\textbf{y},\label{A.31}
\end{align}
where $\rho_0^{\v}$ is the function defined in \eqref{A.7-1}.
Clearly, $\tilde{u}_{0}^{\v,b}(\mathbf{x})$ is a spherically symmetric function,
{\it i.e.}, $\tilde{u}_{0}^{\v,b}(\mathbf{x})=\tilde{u}_{0}^{\v,b}(r)$.

\begin{lemma}\label{lemA.5}
	$u^\v_0(\mathbf{x})$ defined in \eqref{A.31-1}  satisfies
	\begin{align}
		&\int_{\mathbb{R}^n}\rho_0^\v(\mathbf{x}) |u_0^\v(\mathbf{x})|^2\dd\mathbf{x}
		\equiv\int_{\mathbb{R}^n}\frac{|m_0(\mathbf{x})|^2}{\rho_0(\mathbf{x})}\dd\mathbf{x}
		\qquad\,\,\,\mbox{for any $\v\in (0,1]$}, \label{8.19}\\[1.5mm]
		&\lim_{\v\rightarrow0+} \|\rho^\v_0 u^\v_0-m_0\|_{L^1(\mathbb{R}^n)}=0.\label{8.19-1}
	\end{align}
	Moreover, $\tilde{u}_0^{\v,b}(\mathbf{x})$ defined in \eqref{A.31} is in $C_0^\infty(\mathbb{R}^{n})$
	and satisfies
	\begin{align}
		&{\rm supp}\, \tilde{u}_0^{\v,b}\subset \big\{\mathbf{x} \in \mathbb{R}^n~ : ~ 2b^{-1}\leq |\mathbf{x}|\leq b-1\big\},\label{8.20o}\\
		&\lim_{b\rightarrow\infty} \int_{\mathbb{R}^n}\rho_0^{\v,b}(\mathbf{x}) |\tilde{u}_0^{\v,b}(\mathbf{x})|^2\,\dd\mathbf{x}
		= \int_{\mathbb{R}^n}\rho_0^\v(\mathbf{x}) |u_0^\v(\mathbf{x})|^2\dd\mathbf{x},\label{8.20}\\[1mm]
		&\lim_{b\rightarrow\infty} \|\rho^{\v,b}_0 \tilde{u}^{\v,b}_0-\rho^\v_0 u^{\v}_0\|_{L^1(\mathbb{R}^n)}=0.\label{8.20-1}
	\end{align}
\end{lemma}

\begin{proof}
\eqref{8.19} follows directly from \eqref{A.31-1}.
Using \eqref{A.6-1} and \eqref{A.31-1}, we have
\begin{align}
	&\int_{\R^n}
	\big|(\rho^\v_0 u^\v_0-m_0)(\mathbf{x})\big|\, \dd\mathbf{x}\nonumber\\
	&=\int_{\R^n}
	\big|\big(\sqrt{\rho^\v_0}-\sqrt{\rho_0}\big)(\mathbf{x})\big(\frac{m_0}{\sqrt{\rho_0}}\big)(\mathbf{x})\big|\,\dd\mathbf{x}\nonumber\\
	&\leq \Big(\int_{\mathbb{R}^n}\frac{|m_0(\mathbf{x})|^2}{\rho_0(\mathbf{x})}\dd\mathbf{x}\Big)^{\frac12}
	\Big( \int_{\R^n}
	\big|\big(\sqrt{\rho^\v_0}-\sqrt{\rho_0}\big)(\mathbf{x})\big|^2\dd\mathbf{x}\Big)^{\frac12}
	\longrightarrow0 \quad\,\, \mbox{as $\v\rightarrow0+$},\nonumber
\end{align}
which leads to \eqref{8.19-1}.

From \eqref{A.31}, it is clear that $\tilde{u}_0^{\v,b}(\mathbf{x})\in C_0^\infty(\mathbb{R}^{n})$
and  ${\rm supp}\, \tilde{u}_0^{\v,b}\subset \{\mathbf{x} \in \mathbb{R}^n~ : ~ 2b^{-1}\leq |\mathbf{x}|\leq b-1\}$.
For any given small constant $\epsilon>0$, there exists small $\sigma=\sigma(\epsilon)>0$ and large $N=N(\epsilon)\gg1$ such that
\begin{align}\label{A.32}
	\int_{B_{2\sigma}(\mathbf{0})\cup\{|\mathbf{x}|\geq N\}} \frac{|m_0(\mathbf{x})|^2}{\rho_0(\mathbf{x})}\dd\mathbf{x}
	\leq \epsilon.
\end{align}
Taking $b>0$ large enough so that $\sigma\geq 6 b^{-1}$, then it follows from  \eqref{A.31} that
\begin{align}\label{A.33}
	\int_{\sigma\leq |\mathbf{x}|\leq N+2} \Big|\Big(\sqrt{\rho_0^{\v,b}} \tilde{u}_0^{\v,b}
	- \frac{m_0}{\sqrt{\rho_0}}\Big)(\mathbf{x})\Big|^2\dd\mathbf{x}\longrightarrow0\qquad\,\, \mbox{as $b\rightarrow\infty$}.
\end{align}
Since $\sigma\geq 6b^{-1}$, we use \eqref{A.32} to obtain
\begin{align}\label{A.34}
&\int_{B_{\sigma}(\mathbf{0})\cup\{|\mathbf{x}|\ge N+1\}}
\big|\sqrt{\rho_0^{\v,b}(\mathbf{x})} \tilde{u}_0^{\v,b}(\mathbf{x})\big|^2\dd\mathbf{x}\\
&\leq \int_{B_{\sigma}(\mathbf{0})\cup\{|\mathbf{x}|\ge N+1\}}
\Big|\int_{\mathbb{R}^{n}} \big(\frac{m_0}{\sqrt{\rho_0}}\mathbf{1}_{[4b^{-1},b-2]}\big)(\mathbf{x-y})
J_{b^{-1}}(\mathbf{y})\, \dd\mathbf{y}\Big|^2\dd\mathbf{x}
\nonumber\\
&  \leq \int_{B_{2\sigma}(\mathbf{0})\cup\{|\mathbf{x}|\ge N\}} \frac{|m_0(\mathbf{x})|^2}{\rho_0(\mathbf{x})}\dd\mathbf{x}\leq \epsilon.\nonumber
\end{align}
It follows from \eqref{A.31-1}--\eqref{A.31} and  \eqref{A.33}--\eqref{A.34} that
\begin{align} \label{A.27}
	&\int_{\mathbb{R}^n} \big|\big(\sqrt{\rho_0^{\v,b}} \tilde{u}_0^{\v,b}
	- \sqrt{\rho_0^{\v}} u_0^{\v} \big)(\mathbf{x})\big|^2\dd\mathbf{x}\\
	&=\int_{\mathbb{R}^n} \big|\big(\sqrt{\rho_0^{\v,b}} \tilde{u}_0^{\v,b}
	- \frac{m_0}{\sqrt{\rho_0}}\big)(\mathbf{x})\big|^2\dd\mathbf{x}\nonumber\\
	&\leq \int_{\sigma\leq |\mathbf{x}|\leq N+2} \big|\big(\sqrt{\rho_0^{\v,b}} \tilde{u}_0^{\v,b}
	- \frac{m_0}{\sqrt{\rho_0}}\big)(\mathbf{x})\big|^2\dd\mathbf{x}\nonumber\\
	&\quad +C\int_{B_{2\sigma}(\mathbf{0})\cup\{|\mathbf{x}|\geq N\}}\frac{|m_0(\mathbf{x})|^2}{\rho_0(\mathbf{x})}\dd\mathbf{x}\nonumber\\
	&\rightarrow 0\nonumber
\end{align}
as $b\rightarrow\infty$, which yields \eqref{8.20}.

Using \eqref{A.20-2}, \eqref{8.20o}, and \eqref{A.27}, we see that, as $b\to \infty$,
\begin{align}
	&\int_{\mathbb{R}^n} \big|(\rho^{\v,b}_0 \tilde{u}^{\v,b}_0-\rho^\v_0 u^{\v}_0)(\mathbf{x})\big| \dd\mathbf{x}\nonumber\\
	&\leq \Big( \int_{\mathbb{R}^n}\rho^\v_0(\mathbf{x})\,\dd\mathbf{x}\Big)^{\frac12}
	\Big(\int_{\mathbb{R}^n} \big|\big(\sqrt{\rho_0^{\v,b}} \tilde{u}_0^{\v,b}
	- \sqrt{\rho_0^{\v}} u_0^{\v}\big)(\mathbf{x})\big|^2\dd\mathbf{x}\Big)^{\frac12}\nonumber\\
	&\quad +\int_{|\mathbf{x}|\leq b} \Big|\big(\sqrt{\rho_{0}^{\v,b}}-\sqrt{\rho_0^\v}\big)(\mathbf{x}) \, \big(\sqrt{\rho_{0}^{\v,b}} \tilde{u}^{\v,b}_0\big)(\mathbf{x})\Big|\,\dd\mathbf{x}\longrightarrow 0,
	\nonumber
\end{align}
which implies \eqref{8.20-1}.
\end{proof}

We still need to modify $\tilde{u}_{0}^{\v,b}(\mathbf{x})$ so that it satisfies the stress-free boundary condition \eqref{2.5} at $r=b$.
Let $S(\cdot)$ be the function in \eqref{A.19-1}. Then we define
\begin{equation}\label{A.37}
	u_{0}^{\v,b}(\mathbf{x}):=\tilde{u}_{0}^{\v,b}(\mathbf{x})
	-\frac{1}{\varepsilon}\,S(4(|\mathbf{x}|-(b-\frac{1}{2})))
	\,\frac{1}{|\mathbf{x}|^{n-1}}\int_{|\mathbf{x}|}^{b} \frac{p(\rho_{0}^{\v,b}(z))}{\rho_{0}^{\v,b}(z)}\,z^{n-1} \dd z.
\end{equation}
It is direct to check that  $u_{0}^{\v,b}(\textbf{x})\in C^{\infty}([b^{-1},b])$ satisfies
the following boundary conditions:
\begin{align}\label{A.40}
	u_{0}^{\v,b}(b^{-1})=0,\qquad \,\, \Big\{p(\rho_{0}^{\v,b})-\v \rho_{0}^{\v,b} \big(u^{\v,b}_{0r}+\frac{n-1}{r} u^{\v,b}_0 \big)\Big\}\Big|_{r=b}=0.
\end{align}
A direct calculation shows from \eqref{A.19} that
\begin{align}\label{A.38}
&\int_{|\textbf{x}|\leq b}\Big|\sqrt{\rho_{0}^{\v,b}(\textbf{x})}u_{0}^{\v,b}(\textbf{x})-\sqrt{\rho_{0}^{\v,b}(\mathbf{x})}\tilde{u}_{0}^{\v,b}(\mathbf{x})\Big|^2
\dd\textbf{x}\\
&\leq\frac{a_0^2}{\varepsilon^2}\int_{b-\frac{1}{2}}^{b}\rho_{0}^{\v,b}(r) r^{n-1}
\Big|\frac{1}{r^{n-1}}\int_{r}^{b}(\rho_{0}^{\v,b}(z))^{\gamma-1}\,z^{n-1} \dd z\Big|^2\dd r\nonumber\\
&\leq C\varepsilon^{-2}b^{-2(\gamma-1)(n-\alpha)} \longrightarrow0\qquad\mbox{as $b\rightarrow \infty$}.\nonumber
\end{align}
Therefore, combining \eqref{8.20o}, \eqref{A.27}, and \eqref{A.38}, we conclude
\begin{lemma}\label{lemA.7}
	For fixed $\v>0$,
	\begin{align*}
		&\lim_{b\rightarrow\infty}\int_{|\mathbf{x}|\leq b}\big|\sqrt{\rho_{0}^{\v,b}(\mathbf{x})}u_{0}^{\v,b}(\mathbf{x})\big|^2\dd\mathbf{x}
		=\int_{\mathbb{R}^n} \big|\sqrt{\rho_0^\v(\mathbf{x})} u_0^\v(\mathbf{x})\big|^2\dd\mathbf{x},\\[2mm]
		&\sqrt{\rho_{0}^{\v,b}(\mathbf{x})}u_{0}^{\v,b}(\mathbf{x})\longrightarrow \sqrt{\rho_{0}^{\v}(\mathbf{x})}u_{0}^{\v}(\mathbf{x})
		\quad \mbox{in $L^2(\{|\mathbf{x}|\leq b\})\, $  as $b\rightarrow\infty$}.
	\end{align*}	
\end{lemma}

With $\rho^{\v}_0(\mathbf{x})$, $\rho^{\v,b}_0(\mathbf{x})$, $u_0^{\v}(\mathbf{x})$,
and $u_0^{\v,b}(\mathbf{x})$ defined respectively  in \eqref{A.7-1}, \eqref{A.23-1}, \eqref{A.31-1}, and \eqref{A.37},
we can construct the approximate initial data
$(\rho_0^{\v,b}, m_0^{\v,b})(r)=(\rho_0^{\v,b}, \rho_0^{\v,b} u_0^{\v,b})(r)$ for \eqref{2.1}--\eqref{initial}:
For $b\gg1$,  define
\begin{align}\label{A.44}
	(\rho^{\v,b}_0, u^{\v,b}_0)(r):=(\rho_{0}^{\v,b}(\mathbf{x}), u^{\v,b}_0(\mathbf{x})) \mathbf{1}_{[b^{-1},b]}(\mathbf{x}).
\end{align}
Then, collecting all the above estimates, we have the following results.

\begin{lemma}\label{propA.1}
	Let $\big(\rho^{\v,b}_0, u^{\v,b}_0\big)(r)$ be the functions defined in \eqref{A.44}
	so that $(\rho^{\v,b}_0,  u^{\v,b}_0)$ is in $C^\infty([b^{-1},b])$ and satisfies the boundary condition \eqref{A.40}.
	Let $q\in \{1, \gamma\}$ for $\kappa=1$ $($gaseous stars$)$, and $q\in\{1,\gamma,\frac{2n}{n+2}\}$ for $\kappa=-1$ $($plasmas$)$.
	Then
	\begin{enumerate}
		\item[\rm (i)]  For all $\v\in(0,1]$,
		\begin{equation*}
			\qquad \int_0^\infty \rho_0^\v(r)\,\omega_n r^{n-1} \dd r=M,\qquad
			\v^2\int_{0}^\infty \big|\partial_r\sqrt{\rho^\v_{0}(r)}\big|^2\, r^{n-1}\dd r
			\leq C\v (M+1).
		\end{equation*}
		Moreover, as $\v\rightarrow0+$,
		\begin{align*}
			&\qquad (E^\v_0, E_1^\v)\longrightarrow (E_0,0), \\
			&\qquad (\rho_0^{\v}, m_0^{\v})(r)\longrightarrow (\rho_0, m_0)(r) \quad \mbox{in $L^q([0,\infty); r^{n-1}\dd r)\times L^1([0,\infty); r^{n-1}\dd r)$},\\
			&\qquad \Phi_{0r}^\v\longrightarrow \Phi_{0r} \quad \mbox{in $L^2([0,\infty); r^{n-1} \dd r)$},
		\end{align*}
		where $E^\v_0, E_1^\v$, and $E_0$ are defined in \eqref{1.13-3}, \eqref{1.13-1}, and \eqref{1.20a}, respectively.
		Moreover, there exists a small constant $\v_0>0$ such that
		\begin{equation*}
			M<M_{\rm c}^\v(\gamma) \quad\mbox{for all $\v\in (0,\v_0]$ and $\g\in(\frac{2n}{n+2}, \frac{2(n-1)}{n}]$},
		\end{equation*}
		where $M_{\rm c}^\v(\gamma)$ is defined in \eqref{2.17}.
		
		\smallskip
		\item[\rm (ii)] For $\v\in(0,1]$ and $b>1$,
		\begin{equation*}
			\qquad \int_{b^{-1}}^{b} \rho_{0}^{\v,b}(r)\,\omega_n r^{n-1} \dd r= M, \quad
			\v^2\int_{0}^{b} \big|\partial_r\sqrt{\rho^{\v,b}_{0}(r)}\big|^2\, r^{n-1}\dd r
			\leq C\v (M+1).
		\end{equation*}
		Moreover, for any fixed $\v\in(0,1]$, as $b\rightarrow\infty$,
		\begin{align*}
			\begin{split}
				&\qquad (E_0^{\v,b}, E_1^{\v,b})\longrightarrow (E_0^\v, E_1^\v), \\
				&\qquad (\rho_0^{\v,b}, m_0^{\v,b})(r)\longrightarrow (\rho_0^\v, m^\v_0)(r) \quad
				\mbox{\rm in $L^q([0, b]; r^{n-1}\dd r)\times L^1([0, b]; r^{n-1}\dd r)$},\\
				&\qquad \Phi_{0r}^{\v,b}\to \Phi^\v_{0r} \quad \mbox{\rm in $L^2([0,\infty); r^{n-1} \dd r)$},
			\end{split}
		\end{align*}
		where $E_0^{\v,b}$ and $E_1^{\v,b}$ are defined in \eqref{2.7-1} and \eqref{2.7-2}, respectively.
		Furthermore, there exists a constant $\mathfrak{B}(\v)\gg 1$ such that
		\begin{equation}\label{A.51}
			M<M_{\rm c}^{\v,b}(\gamma)\qquad\mbox{for $\v\in (0,\v_0],\,\, b\geq \mathfrak{B}(\v)$, and $\g\in(\frac{2n}{n+2}, \frac{2(n-1)}{n}]$},
		\end{equation}
		where $M_{\rm c}^{\v,b}(\gamma)$ is defined in \eqref{2.36}.
	\end{enumerate}
\end{lemma}

\noindent{\bf Acknowledgments.}
The research of Gui-Qiang G. Chen was supported in part by the UK Engineering and Physical
Sciences Research Council Awards EP/L015811/1, EP/V008854, and EP/V051121/1. The research
of Lin He was supported in part by the National Natural Science Foundation of China Grant No.
12001388 and No. 12371223, and the Sichuan Younth Science and Technology Foundation Grant
No. 2021JDTD0024, the National Key R\&D Program of China Grant No. 2022YFA1007700, and
the Natural Science Foundation of Sichuan Province Grant No. 2023NSFSC1367. The research of
Yong Wang was supported in part by the National Natural Science Foundation of China Grant No.
12022114, No. 12288201, and No. 12071397, the CAS Project for Young Scientists in Basic Research
Grant No. YSBR-031, and the Youth Innovation Promotion Association of the Chinese Academy
of Science Grant No. 2019002. The research of Difan Yuan was supported in part by the National
Natural Science Foundation of China Grant No. 12001045, the China Postdoctoral Science Foundation Grant No. 2020M680428 and No. 2021T140063,
and the UK Engineering and Physical Sciences
Research Council Award EP/V051121/1

\frenchspacing
\bibliographystyle{cpam}

\begin{thebibliography}{99}
	
\bibitem{BD-2003-CRMASP}
 Bresch D.; Desjardins B. \,
 On viscous shallow-water equations (Saint-Venant model) and the quasi-geostrophic limit,
\textit{C. R. Math. Acad. Sci. Paris}, \textbf{335} (2002), 1079--1084.

\bibitem{BD-2004-CRMASP}
 Bresch D.;  Desjardins B. \,
Some diffusive capillary models of Korteweg type,
\textit{ C. R. Math. Acad. Sci. Paris, Section M\'{e}canique}, \textbf{332} (2004), 881--886.

\bibitem{BD-2007-JMPA}
 Bresch D.;  Desjardins B. \,
On the existence of global weak solutions to the Navier-Stokes equations for viscous compressible and heat conducting fluids,
\textit{J. Math. Pures Appl.}  \textbf{87} (2007), 57--90.

\bibitem{BD-2003-CPDE}
 Bresch D.; Desjardins B.;  Lin C. K. \,
On some compressible fluid models: Korteweg, lubrication, and shallow water systems,
\textit{Comm. Partial Diff. Eqs.}  \textbf{28} (2003), 843--868.

\bibitem{Chandrasekhar}
Chandrasekhar S. \,
\textit{ An Introduction to the Study of Stellar Structures},
University of Chicago Press, Chicago, 1938.

\bibitem{Chen2}
Chen G.-Q. \,
Convergence of the Lax-Friedrichs scheme
for isentropic gas dynamics (III),
\textit{ Acta Math. Sci.}  \textbf{6B} (1986), 75--120 (in English);
\textbf{8A} (1988), 243--276 (in Chinese).

\bibitem{Chen-Feldman2018}
Chen G.-Q.;  Feldman M. \,
\textit{The Mathematics of Shock Reflection-Diffraction and von Neumann's Conjectures},
Research Monograph, Annals of Mathematics Studies, {197},
Princeton University Press, Princeton, 2018.

\bibitem{Perepelitsa}
 Chen G.-Q.;  Perepelista M. \,
Vanishing viscosity limit of the Navier-Stokes equations to the Euler equations for compressible fluid flow,
\textit{Comm. Pure Appl. Math.}  \textbf{63} (2010), 1469--1504.

\bibitem{Chen6}
 Chen G.-Q.;  Perepelista M. \,
Vanishing viscosity limit solutions of the compressible Euler equations
with spherically symmetry and large initial data,
\textit{Commun. Math. Phys.}  \textbf{338} (2015), 771--800.

\bibitem{Chen-Schrecker-2018}
 Chen G.-Q.; Schrecker M. \,
Vanishing viscosity approach to the compressible Euler equations for transonic nozzle and spherically symmetric flows,
\textit{Arch. Ration. Mech. Anal.}  \textbf{229} (2018), 1239--1279.

\bibitem{Chen-Wang}
 Chen G.-Q.;  Wang  D.-H. \,
Formation of singularities in compressible Euler-Poisson fluids with heat diffusion and damping relaxation,
\textit{J. Differential Equations}  \textbf{144} (1998), 44--65.

\bibitem{Chen-Wang-2018}
 Chen G.-Q.;   Wang Y. \,
Global solutions of the compressible Euler equations with large initial data of spherical symmetry and positive far-field density,
\textit{Arch. Ration. Mech. Anal.}  \textbf{243} (2022), 1699--1771.


\bibitem{CS-2012}
Coutand D.;  Shkoller S. \,
Well-posedness in smooth function spaces for the moving-boundary three
dimensional compressible Euler equations in physical vacuum,
\textit{Arch. Ration. Mech. Anal.}  \textbf{206} (2012), 515--616.

\bibitem{Courant-Friedrichs-1962}
 Courant R.;  Friedrichs K. O. \,
\textit{Supersonic Flow and Shock Waves}, Springer-Verlag, New York, 1962.

\bibitem{Dafermos}
 Dafermos C. M. \,
\textit{Hyperbolic Conservation Laws in Continuum Physics}, Springer-Verlag, Berlin, 2016.

\bibitem{Chen5}
 Ding X.;  Chen G.-Q.;   Luo P. \,
Convergence of the Lax-Friedrichs scheme for the isentropic gas dynamics (I)--(II),
\textit{Acta Math. Sci.}  \textbf{5B} (1985), 483--500, 501--540 (in English);
\textbf{7A} (1987), 467--480; \textbf{8A} (1989), 61--94 (in Chinese).

\bibitem{Chen5b}
 Ding X.;  Chen G.-Q.;   Luo P. \,
Convergence of the fractional step Lax-Friedrichs scheme and Godunov scheme for the isentropic system of gas dynamics,
\textit{Commun. Math. Phys.}  \textbf{121} (1989), 63--84.

\bibitem{R.J.DiPerna2}
 DiPerna R. J. \,
 Convergence of the viscosity method for isentropic gas dynamics,
\textit{Commun. Math. Phys.}  \textbf{91} (1983), 1--30.

\bibitem{Duan-Li}
Duan Q.;   Li H.-L. \,
Global existence of weak solution for the compressible Navier-Stokes-Poisson system for gaseous stars,
\textit{J. Differential Equations}  \textbf{259} (2015), 5302--5330.

\bibitem{Ducomet-Feireisl-P-2004}
 Ducomet B.; Feireisl E.;  Petzeltova H.;  Stra$\check{\mbox{s}}$kraba  I. \,
Global in time weak solutions for compressible barotropic self-gravitating fluids,
\textit{Discrete Contin. Dyn. Syst.}  \textbf{11} (2004), 113--130.

\bibitem{Feireisl-2001}
Feireisl E.; Novotny  A.;   Petzeltova H. \,
On the existence of globally defined weak solutions to the Navier-Stokes equations of isentropic compressible fluids,
\textit{J. Math. Fluid Mech.}  \textbf{3} (2001), 358--392.

\bibitem{Fu-Lin-1998}
Fu C.-C.;   Lin S.-S. \,
On the critical mass of the collapse of a gaseous star in spherically symmetric and isentropic motion,
\textit{Japan J. Indust. Appl. Math.}   \textbf{15} (1998), 461--469.

\bibitem{D. Gilbarg}
Gilbarg D. \,
The existence and limit behavior of the
one-dimensional shock layer,
\textit{Amer. J. Math.}  \textbf{73} (1951), 256--274.

\bibitem{Goldreich-Weber-1980}
Goldreich P.;  Weber S. \,
Homologously collapsing stellar cores,
\textit{Astrophys. J.}  \textbf{238} (1980), 991--997.

\bibitem{GMWZ}
 Gu\`{e}s C.M.I.O.;  M\`{e}tivier G.;   Williams M.;   Zumbrun K. \,
Navier-Stokes regularization of multidimensional Euler shocks,
\textit{Ann. Sci. \'{E}cole Norm. Sup. (4)},  \textbf{39} (2006), 75--175.

\bibitem{Guo-1998}
Guo Y. \,
Smooth irrotational flows in the large to the Euler-Poisson system
in $\mathbb{R}^{3+1}$,
\textit{Commun. Math. Phys.}  \textbf{195} (1998), 249--265.

\bibitem{Guo-H-J}
Guo Y.;  Had$\check{\mbox{z}}$i\'c M.;  Jang J. \,
Continued gravitational collapse for Newtonian stars,
\textit{Arch. Ration. Mech. Anal.}  \textbf{239} (2021), 431--552.

\bibitem{Guo-H-J-1}
Guo Y.;  Had$\check{\mbox{z}}$i\'c M.;  Jang J. \,
Larson-Penston self-similar gravitational collapse,
\textit{Commun. Math. Phys.}  \textbf{386}  (2021), 1551--1601.

\bibitem{Guo-Ioescu-Pausader-2016}
Guo Y.;  Ionescu A.~D.;  Pausader B. \,
Global solutions of the Euler-Maxwell two-fluid system in 3D,
\textit{Ann. of Math. (2)},  \textbf{183} (2016), 377--498.


\bibitem{Guo-Pausader-2011}
 Guo Y.;   Pausader B. \,
Global smooth ion dynamics in the Euler-Poisson system,
\textit{Commun. Math. Phys.}  \textbf{303} (2011), 89--125.


\bibitem{Zhenhua Guo2}
 Guo Z. H.;  Jiu Q. S.;  Xin Z. \,
Spherically symmetric isentropic compressible flows with density-dependent viscosity
coefficients,
\textit{SIAM J. Math. Anal.}  \textbf{39} (2008), 1402--1427.

\bibitem{Hadzic-Jang-2018}
Had$\check{\mbox{z}}$i\'c M.;  Jang J. \,
Nonlinear stability of expanding star solutions of the radially symmetric
mass-critical Euler-Poisson system,
\textit{Comm. Pure Appl. Math.}  \textbf{71} (2018), 827--891.

\bibitem{Hadzic-Jang-2017}
Had$\check{\mbox{z}}$i\'c M.;  Jang J. \,
A class of global solutions to the Euler-Poisson system,
\textit{Commun. Math. Phys.}  \textbf{370} (2019),  475--505.

\bibitem{Huang-Wang}
 Huang F. M.;   Wang  Z. \,
Convergence of viscosity solutions for isothermal gas dynamics,
\textit{SIAM J. Math. Anal.}  \textbf{34} (2002), 595--610.

\bibitem{Hoff-1995}
Hoff D. \,
Strong convergence to global solutions for multidimensional flows
of compressible, viscous fluids with polytropic equations of state and discontinuous initial data,
\textit{Arch. Ration. Mech. Anal.}  \textbf{132} (1995), 1--14.

\bibitem{D. Hoff-Liu}
Hoff D.;   Liu T.-P. \,
The inviscid limit for the Navier-Stokes equations of compressible, isentropic flow with shock data,
\textit{Indiana Univ. Math. J.}  \textbf{38} (1989), 861--915.

\bibitem{Hugoniot}
Hugoniot P. H. \,
M\`{e}moire sur la propagation du mouvement dans les corps et sp\'{e}cialement dans les gaz parfaits,
2e Partie.
\textit{J. \'{E}cole Polytechnique Paris}, \textbf{58} (1889), 1--125.

\bibitem{Ionescu-Pausader-2013}
 Ionescu A. D.; Pausader  B. \,
 The Euler-Poisson system in 2D: global stability of the constant equilibrium solution,
\textit{ Int. Math. Res. Not. IMRN}  \textbf{4} (2013), 761--826.

\bibitem{Jang-2010}
Jang J. \,
Local well-posedness of dynamics of viscous gaseous stars,
\textit{ Arch. Ration. Mech. Anal.}  \textbf{195} (2010), 797--863.

\bibitem{JM-2015}
Jang J.; Masmoudi  N. \,
Well-posedness of compressible Euler equations in a physical vacuum,
\textit{ Comm. Pure Appl. Math.} \textbf{68} (2015), 61--111.

\bibitem{Jiang-Zhang-2001}
Jiang S.; Zhang  P. \,
On spherically symmetric solutions of the compressible isentropic Navier-Stokes equations,
\textit{ Commun. Math. Phys.} \textbf{215} (2001), 559--581.

\bibitem{Kobayashi-Suzuki-2008}
Kobayashi T.;   Suzuki T. \,
Weak solutions to the Navier-Stokes-Poisson equations,
\textit{Adv. Math. Sci. Appl.} \textbf{18} (2008), 141--168.

\bibitem{Kong-Li-2018}
 Kong H. H.;  Li  H.-L. \,
Free boundary value problem to 3D spherically symmetric compressible Navier-Stokes-Poisson equations,
\textit{Z. Angew. Math. Phys.} \textbf{68} (2017), Paper No. 21, 34 pp.

\bibitem{P. D. Lax}
Lax P. D. \,
Shock wave and entropy, \textit{Contributions to Functional Analysis},
ed. E. A. Zarantonello, pp. 603--634,
Academic Press, New York, 1971.

\bibitem{Ph. LeFloch}
LeFloch P.;  Westdickenberg  M. \,
Finite energy solutions to the isentropic Euler equations with geometric effects,
\textit{J. Math. Pures Appl.} \textbf{ 88} (2007), 386--429.

\bibitem{Lei-Gu}
Lei Z.;   Gu X. M. \,
Local well-posedness of the three dimensional compressible Euler-Poisson equations with physical vacuum,
\textit{J. Math. Pures Appl. (9)}, \textbf{105} (2016), 662--723.

\bibitem{Li-Xin-2015}
Li  J.;  Xin Z. \,
Global existence of weak solutions to the barotropic compressible Navier-Stokes flows with degenerate viscosities,
\textit{arXiv Preprint},  arXiv: 1504.06826v1.

\bibitem{Li-Wang-2006}
Li T.-H.;  Wang  D.-H. \,
Blow up phenomena of solutions to the Euler equations for compressible fluid flow,
\textit{J. Differential Equations} \textbf{22} (2006), 91--101.

\bibitem{Lieb-Loss}
 Lieb E. H.; Loss  M. \,
 {\em Analysis}, Second edition, Graduate Studies in Mathematics, \textbf{14},
 American Mathematical Society: Providence, 2001


\bibitem{Lin}
 Lin S.-S. \,
Stability of gaseous stars in spherically symmetric motions,
\textit{ SIAM J. Math. Anal.} \textbf{28} (1997), 539--569.

\bibitem{Lions-CNS-1998}
 Lions P. L. \,
\textit{Mathematical Topics in Fluid Dynamics 2: Compressible Models},
Oxford Science Publication, Oxford, 1998.

\bibitem{Lions P.-L.1}
 Lions P. L.; Perthame B.; Souganidis P. \,
Existence and stability of entropy solutions for the hyperbolic systems of isentropic gas dynamics in Eulerian and Lagrangian coordinates,
\textit{Comm. Pure Appl. Math.} \textbf{49} (1996), 599--638.

\bibitem{Lions P.-L.2}
Lions P. L.; Perthame B.; Tadmor E. \,
Kinetic formulation of the isentropic gas dynamics and p-systems,
\textit{Commun. Math. Phys.} \textbf{163} (1994), 415--431.

\bibitem{Luo-Xin-Zeng-2016-Adv}
Luo T.;  Xin Z.;   Zeng H. \,
On nonlinear asymptotic stability of the Lane-Emden solutions for the viscous gaseous star problem,
\textit{Adv. Math.} \textbf{291} (2016), 90--182.

\bibitem{Luo-Xin-Zeng-2016-CMP}
Luo T.;  Xin Z.;   Zeng H. \,
Nonlinear asymptotic stability of the Lane-Emden solutions for the viscous gaseous star problem
with degenerate density dependent viscosities,
\textit{Commun. Math. Phys.}  \textbf{347} (2016), 657--702.

\bibitem{Luo-Xin-Zeng-2014}
Luo T.;  Xin Z.;   Zeng H. \,
Well-posedness for the motion of physical vacuum of the three-dimensional compressible Euler equations with or without self-gravitation,
\textit{Arch. Ration. Mech. Anal.} \textbf{213} (2014), 763--831.

\bibitem{Okada-Makino-1993}
Okada M.;  Makino  T. \,
Free boundary value problems for the equation of spherically symmetrical motion of viscous gas,
\textit{Japan J. Ind. Appl. Math.}  \textbf{10} (1993), 219--235.

\bibitem{Makino-1986}
Makino T. \,
On a local existence theorem for the evolution equation of gaseous stars,
\textit{ Patterns and Waves}, pp. 459--479, \textit{ Stud. Math. Appl.} 18, North-Holland, Amsterdam, 1986.

\bibitem{Makino-1992}
Makino T. \,
Blowing up solutions of the Euler-Poisson equations for the evolution of gaseous stars,
\textit{Transport Theory Statist. Phys.} \textbf{21} (1992), 615--624.

\bibitem{Makino-1997}
Makino T. \,
On the spherically symmetric motion of self-gravitating isentropic gas surrounding
a solid ball, \textit{ Nonlinear Evolutionary Partial Differential Equations}, pp. 543--546,
AMS/IP Studies in Advanced Mathematics, Vol. 3, AMS, Providence, 1997.

\bibitem{Matsumura-1980}
Matsumura  A.;  Nishida T. \,
The initial value problem for the equations of motion of viscous and heat-conductive gases,
\textit{Publ. Res. Inst. Math. Sci.} \textbf{20} (1980), 67--104.

\bibitem{Mellet}
Mellet A.;  Vasseur A. \,
On the barotropic compressible Navier-Stokes equation,
\textit{Comm. Partial Diff. Eqs.} \textbf{32} (2007), 431--452.

\bibitem{MRRS-2020}
Merle F.;  Raphael P.;  Rodnianski I.;  Szeftel  J. \,
On the implosion of a compressible fluid I: Smooth self-similar inviscid profiles,
Ann. of Math. (2), {\bf 196} (2022),  567--778;
II: Singularity formation, Anna. of Math. (2), {\bf 196} (2022),  779--889. "


On the implosion of a three dimensional compressible fluid,
\textit{arXiv Preprint}, arXiv:1912.11009v2, 2020.

\bibitem{F. Murat}
Murat F. \,
Compacit\'e par compensation,
{\em Ann. Scuola Norm. Sup. Pisa Sci. Fis. Mat.}  \textbf{5} (1978), 489--507.

\bibitem{Rankine}
Rankine W. J. M. \,
On the thermodynamic theory of waves of finite longitudinal disturbance,
 \textit{Phi. Trans. Royal Soc. London}, \textbf{1870} (1960), 277--288.

\bibitem{Lord Rayleigh}
Rayleigh L.  (Strutt J. W.), \,
Aerial plane waves of finite amplitude,
\textit{ Proc. Royal Soc. London}, \textbf{84A} (1910), 247--284.

\bibitem{Rosseland}
Rosseland S. \,
\textit{ The Pulsation Theory of Variable Stars},
Dover Publications Inc., New York, 1964.

\bibitem{Schrecker}
Schrecker M. \,
Spherically symmetric solutions of the multidimensional compressible, isentropic Euler equations,
\textit{Trans. Amer. Math. Soc.}  \textbf{373} (2020), 727--746.

\bibitem{Serre}
Serre D. \,
Compensated integrability. Applications to the Vlasov-Poisson equation and other models in mathematical physics,
\textit{J. Math. Pures Appl. (9)} \textbf{127} (2019), 67--88.


\bibitem{Stokes}
Stokes G. G. \,
On a difficulty in the theory of sound. [\textit{Philos. Mag.} {\bf 33} (1848), 349--356;
\textit{Mathematical and Physical Papers}, Vol. II, Cambridge Univ. Press, Cambridge, 1883].
Classic Papers in Shock Compression Science, 71--79.
\textit{High-pressure Shock Compression of Condensed Matter}, Springer-Verlag, New York, 1998.

\bibitem{L. Tartar}
Tartar L. \, Compensated compactness and applications to partial
differential equations, \textit{Research Notes in Mathematics, Nonlinear Analysis and Mechanics}, Herriot-Watt Symposium, Vol. {4}, R.~J. Knops, ed., Pitman Press, 1979.

\bibitem{VY}
Vasseur A.;   Yu  C. \,
Existence of global weak solutions for 3D degenerate compressible Navier-Stokes equations,
\textit{Invent. Math.} \textbf{206} (2016), 935--974.


\bibitem{Whitham-1974}
 Whitham G. B. \, \textit{Linear and Nonlinear Waves}, Wiley, New York, 1974.

\bibitem{Xiao-2016}
Xiao J.-J. \,  Global weak entropy solutions to the Euler-Poisson system
with spherical symmetry, \textit{ Math. Models Methods Appl. Sci.}  \textbf{26} (2016), 1689--1734.

\bibitem{Xin-1993}
Xin Z. \,
Zero dissipation limit to rarefaction waves for the one-dimentional Navier-Stokes equations
of compressible isentropic gases, \textit{Comm. Pure Appl. Math.} \textbf{ 46} (1993), 621--665.

\bibitem{Zhang-Fang-2009}
Zhang T.;  Fang  D. Y. \,
Global behavior of spherically symmetric Navier-Stokes-Poisson system with degenerate viscosity coefficients,
\textit{Arch. Ration. Mech. Anal.} \textbf{191} (2009), 195--243.
\end{thebibliography}

\end{document}